\documentclass{amsart}
\usepackage{amssymb}
\usepackage{latexsym}

\date{\today}

\newcommand{\Z}{{\mathbb Z}}
\newcommand{\R}{{\mathbb R}}
\newcommand{\C}{{\mathbb C}}

\newcommand{\T}{{\mathbb T}}
\newcommand{\Q}{{\mathbb Q}}

\newcommand{\E}{{\mathbb E}}

\newcommand{\Y}{{\mathbb P}}

\newtheorem{theorem}{Theorem}[section]

\newtheorem{lemma}[theorem]{Lemma}
\newtheorem{example}[theorem]{Example}
\newtheorem{prop}[theorem]{Proposition}
\newtheorem{coro}[theorem]{Corollary}

\newtheorem{amoprob}{AMO-Problem}
\newtheorem*{ssprob}{Skew-Shift-Problem}
\renewcommand{\Im}{\mathrm{Im} \, }

\begin{document}

\title[Schr\"odinger Operators with Dynamically Defined Potentials]{Schr\"odinger Operators with Dynamically Defined Potentials: A Survey}

\dedicatory{Dedicated to the memory of Rosario {\rm (}Cherie{\rm )} Bautista Galvez}

\author[D.\ Damanik]{David Damanik}

\address{Department of Mathematics, Rice University, Houston, TX~77005, USA}

\thanks{The author was supported in part by NSF grants DMS--1067988 and DMS--1361625.}

\begin{abstract}
In this survey we discuss spectral and quantum dynamical properties of discrete one-dimensional Schr\"odinger operators whose potentials are obtained by real-valued sampling along the orbits of an ergodic invertible transformation. After an introductory part explaining basic spectral concepts and fundamental results, we present the general theory of such operators, and then provide an overview of known results for specific classes of potentials. Here we focus primarily on the cases of random and almost periodic potentials.
\end{abstract}

\maketitle

\tableofcontents

\section{Introduction}\label{s.1}

The spectral theory of Schr\"odinger operators with random or almost periodic potentials has been an area of very active study since the late 1970's. Essentially from the very beginning it has been understood and emphasized that these two classes of models share an important property, namely that the potentials can be generated dynamically. On the one hand, this makes a unified proof of basic spectral results possible, such as the almost sure constancy of the spectrum and the spectral type, since they hold as soon as the dynamical framework is fixed and an ergodic measure is chosen. On the other hand, by the very nature of the dynamical definition of the potentials, it comes as no surprise that tools from dynamics will enter the spectral analysis of these operators.

The field has made striking progress in the past ten years. This era was ushered in by Puig's proof that localization implies Cantor spectrum for the almost Mathieu operator via Aubry duality and reducibility \cite{P04}, and by the proof of zero-measure spectrum for the critical almost Mathieu operator with recurrent Diophantine frequencies by Avila and Krikorian \cite{AK06}. Both of these papers were released about ten years ago. Since then Avila has gone on to solve many of the major open problems for one-frequency quasi-periodic Schr\"odinger operators, and this work partly led to the Fields medal he was awarded in 2014. Moreover, through a series of papers in the past six years, the fine spectral properties of the Fibonacci Hamiltonian, which is the central quasicrystal model in one dimension, have been understood more or less completely, and this development was only possible after Gorodetski and his students had introduced new ideas from uniformly and partially hyperbolic dynamics in the spectral analysis of this operator.

This only reinforces the second point above. Many advances in the spectral theory of Schr\"odinger operators with dynamically defined potentials are only possible through the application of sophisticated techniques from dynamical systems, and this has come to a whole new level of fruition once researchers with a primary background in dynamics entered the field.

The purpose of this survey is to introduce the non-specialist reader to this field and some of the recent advances. We will assume that the reader has no detailed knowledge of the spectral analysis of Schr\"odinger operators.

A particular goal of ours is to explain why the recent results are interesting and for this reason we will begin with a section on quantum dynamics and spectral theory. Statements such as a certain operator having purely absolutely continuous spectrum are not very meaningful without the connection between such a spectral property and its consequences for the associated quantum evolution. We will explain this connection and survey the known results on how to derive bounds on the quantum evolution from statements about the operator driving the system.

Even though this paper is not short, to keep the length of it in check, very few ideas, arguments, and proofs can be presented. The proofs of most of the results stated here, and much more material, will be contained in the forthcoming monograph \cite{DF15}.

\section{Quantum Dynamics and Spectral Theory}\label{s.2}

In this section we discuss the dynamics of the time-dependent Schr\"odinger equation and ways of bounding the time-evolution. The classical way of doing this proceeds via separation of variables. That is, one studies the time-independent Schr\"odinger equation, and establishes connections between the solutions of this equation and the spectral properties of the associated Schr\"odinger operator. In the second step, one relates the latter to the behavior of the solutions of the time-dependent equation. The net result is the following: the slower the solutions of the time-independent equation decay or grow, the more continuous are the spectral measures of the Schr\"odinger operator in question, and the more rapid is the transport through the medium as modeled by the time-dependent Schr\"odinger equation. These connections are quantitative and allow one in principle to establish fine, and sometimes tight, estimates on the spreading behavior of a wavepacket. An important caveat is that these are in general only one-sided connections, but spectral information may sometimes be supplemented by further information in order to bound transport in the other direction as well.

\subsection{The Time-Dependent Schr\"odinger Equation}\label{ss.tdse}

Let us discuss the time-dependent Schr\"odinger equation. For the sake of simplicity we will focus on the special case we will be interested in. Namely, the evolution takes place in the Hilbert space
$$
\ell^2(\Z) = \Big\{ \psi : \Z \to \C : \sum_{n \in \Z} |\psi(n)|^2 < \infty \Big\},
$$
equipped with the inner product
$$
\langle \varphi, \psi \rangle = \sum_{n \in \Z} \overline{\varphi(n)} \psi(n),
$$
the norm
$$
\|\psi\| = \sqrt{\langle \psi , \psi \rangle} = \sqrt{\sum_{n \in \Z} |\psi(n)|^2},
$$
and the metric
$$
d(\varphi, \psi) = \| \varphi - \psi \|.
$$

Moreover, we consider a bounded self-adjoint operator $H$ in $\ell^2(\Z)$. That is, $H : \ell^2(\Z) \to \ell^2(\Z)$ is a linear map that is everywhere defined and obeys
$$
\|H\| := \sup_{\psi \not= 0} \frac{\|H \psi\|}{\|\psi\|} < \infty
$$
and
$$
\langle \varphi, H \psi \rangle = \langle H \varphi, \psi \rangle
$$
for every $\varphi, \psi \in \ell^2(\Z)$. The \emph{resolvent set} of $H$ is defined by
$$
\rho(H) = \{ E \in \C : (H-E)^{-1} \text{ exists and is bounded} \}.
$$
Here, $H-E$ is shorthand for $H - E \cdot I$, where $I$ denotes the identity operator $I(\psi) = \psi$. Existence of $(H-E)^{-1}$ means that $H-E$ is one-to-one and onto, and boundedness of $(H-E)^{-1}$ means that
$$
\sup_{\psi \not= 0} \frac{\|(H-E)^{-1} \psi\|}{\|\psi\|} < \infty.
$$
The \emph{spectrum} of $H$ is the set $\C \setminus \rho(H)$. All points $E$ in the spectrum of $H$ satisfy $|E| \le \|H\|$.

The self-adjointness of $H$ has a number of important consequences. First of all, the spectrum is real and boundedness of the inverse is automatic whenever $H - E$ is invertible, that is,
\begin{equation}\label{e.spectrumdef}
\sigma(H) = \{ E \in \R : (H-E)^{-1} \text{ does not exist} \}.
\end{equation}
Moreover, the functional calculus allows one to apply functions to $H$. Finally, the spectral theorem associates Borel measures with vectors in $\ell^2(\Z)$, which are called spectral measures. In particular, if $\mu_\psi$ denotes the spectral measure associated with $\psi \in \ell^2(\Z)$, then
\begin{equation}\label{e.spectraltheorem}
\langle \psi, g(H) \psi \rangle = \int g \, d\mu_\psi
\end{equation}
for every locally bounded measurable function $g : \R \to \C$. Every measure $\mu_\psi$ is supported by the spectrum, and conversely the spectrum is the smallest closed set that supports all spectral measures.

The \emph{time-independent Schr\"odinger equation} associated with $H$ is
\begin{equation}\label{e.timedepse}
i \partial_t \psi = H \psi, \quad \psi(\cdot,0) = \psi_0.
\end{equation}
The functional calculus allows one to write the solution of this equation in the form
\begin{equation}\label{e.timedepsesol}
\psi(\cdot,t) = e^{-itH} \psi_0.
\end{equation}

\subsection{The RAGE Theorem}\label{ss.rage}

Recall that the spectral theorem associates with each $\psi \in \ell^2(\Z)$ a Borel measure $\mu_\psi$ on $\R$ so that \eqref{e.spectraltheorem} holds. It follows in particular that $\mu_\psi(\R) = \|\psi\|^2$ (choose $g(E) \equiv 1$ in \eqref{e.spectraltheorem}). Let us consider the (unique) decomposition of $\mu_\psi$ into an absolutely continuous piece, a singular continuous piece, and a pure point piece. That is, we consider
\begin{equation}\label{e.measuredecompstandard}
\mu_\psi = \mu_{\psi,\mathrm{ac}} + \mu_{\psi,\mathrm{sc}} + \mu_{\psi,\mathrm{pp}},
\end{equation}
where $\mu_{\psi,\mathrm{ac}}$ gives zero weight to sets of zero Lebesgue measure, $\mu_{\psi,\mathrm{sc}}$ gives zero weight to individual points and is supported by some set of zero Lebesgue measure, and $\mu_{\psi,\mathrm{pp}}$ is supported by some countable set.

Given this measure decomposition, one can define the following subsets of $\ell^2(\Z)$:
\begin{align*}
\ell^2(\Z)_\mathrm{ac} & = \{ \psi \in \ell^2(\Z) : \mu_\psi = \mu_{\psi,\mathrm{ac}} \} , \\
\ell^2(\Z)_\mathrm{sc} & = \{ \psi \in \ell^2(\Z) : \mu_\psi = \mu_{\psi,\mathrm{sc}} \} , \\
\ell^2(\Z)_\mathrm{pp} & = \{ \psi \in \ell^2(\Z) : \mu_\psi = \mu_{\psi,\mathrm{pp}} \} .
\end{align*}
Each of these subsets turns out to be a closed subspace, and one has
$$
\ell^2(\Z) = \ell^2(\Z)_\mathrm{ac} \oplus \ell^2(\Z)_\mathrm{sc} \oplus \ell^2(\Z)_\mathrm{pp}.
$$
One also considers the continuous subspace
$$
\ell^2(\Z)_\mathrm{c} = \ell^2(\Z)_\mathrm{ac} \oplus \ell^2(\Z)_\mathrm{sc}
$$
and the singular subspace
$$
\ell^2(\Z)_\mathrm{s} = \ell^2(\Z)_\mathrm{sc} \oplus \ell^2(\Z)_\mathrm{pp}.
$$
Then, $\psi \in \ell^2(\Z)_\mathrm{c}$ (resp., $\psi \in \ell^2(\Z)_\mathrm{s}$) if and only if $\mu_\psi = \mu_{\psi,\mathrm{c}} := \mu_{\psi,\mathrm{ac}} + \mu_{\psi,\mathrm{sc}}$ (resp., $\mu_\psi = \mu_{\psi,\mathrm{s}} := \mu_{\psi,\mathrm{sc}} + \mu_{\psi,\mathrm{pp}}$).

Each of the subspaces above reduces the operator $H$, and hence we can consider the restriction of $H$ to it. The spectrum of the restriction of $H$ to $\ell^2(\Z)_\mathrm{ac}$ is denoted by $\sigma_\mathrm{ac}(H)$ and called the \emph{absolutely continuous} spectrum of $H$. The sets $\sigma_\mathrm{sc}(H)$, $\sigma_\mathrm{pp}(H)$, $\sigma_\mathrm{c}(H)$, and $\sigma_\mathrm{s}(H)$ are defined similarly. By convention, one of these sets is empty if and only if the corresponding subspace is trivial (i.e., consists only of the zero vector).

One says that $H$ \emph{has purely absolutely continuous spectrum} if $\ell^2(\Z)_\mathrm{ac} = \ell^2(\Z)$, and similarly for the other cases. Note that having purely absolutely continuous spectrum is not equivalent to $\sigma(H) = \sigma_\mathrm{ac}(H)$! If $H$ has purely absolutely continuous spectrum, it does follow that $\sigma(H) = \sigma_\mathrm{ac}(H)$ (and $\sigma_\mathrm{sc}(H) = \sigma_\mathrm{pp}(H) = \emptyset$), but conversely it is possible to have $\sigma(H) = \sigma_\mathrm{ac}(H)$ and $\ell^2(\Z)_\mathrm{ac} \not= \ell^2(\Z)$. Determining which of the subspaces $\ell^2(\Z)_\mathrm{ac}, \ell^2(\Z)_\mathrm{sc}, \ell^2(\Z)_\mathrm{pp}$ are non-trivial is referred to as determining the \emph{spectral type} of $H$.

\bigskip

The RAGE Theorem makes statements about the behavior of the solutions \eqref{e.timedepsesol} of the time-dependent Schr\"odinger equation \eqref{e.timedepse} in cases when the initial state $\psi_0$ belongs to one of the subspaces above.

\begin{theorem}[RAGE Theorem]\label{t.ragethm}
{\rm (a)} We have $\psi_0 \in \ell^2(\Z)_\mathrm{pp}$ if and only if for every $\varepsilon > 0$, there is $N \in \Z_+$ such that
$$
\sum_{|n| \ge N} | \langle \delta_n , e^{-itH} \psi_0 \rangle |^2 < \varepsilon \quad \text{ for every } t \in \R.
$$

{\rm (b)} We have $\psi_0 \in \ell^2(\Z)_\mathrm{c}$ if and only if for every $N \in \Z_+$,
$$
\lim_{T \to \infty} \frac{1}{2T} \int_{-T}^T \sum_{|n| \le N} | \langle \delta_n , e^{-itH} \psi_0 \rangle |^2 \, dt = 0.
$$

{\rm (c)} If $\psi_0 \in \ell^2(\Z)_\mathrm{ac}$, then for every $N \in \Z_+$,
$$
\lim_{|t| \to \infty} \sum_{|n| \le N} | \langle \delta_n , e^{-itH} \psi_0 \rangle |^2 = 0.
$$
\end{theorem}

In other words, if the spectral measure of the initial state is pure point, then the evolution is confined to a suitable finite set up to an arbitrarily small portion of the total weight; if the spectral measure of the initial state is continuous, then the time-averaged evolution leaves any finite set, and if the spectral measure of the initial state is absolutely continuous, then the evolution leaves any finite set even without any time-averaging.

This is the most basic instance of a general principle: the more continuous the spectral measure, the more the evolution spreads out in space with time. An implication of this kind can be made more quantitative. This is the objective of the following subsection.

\subsection{Hausdorff-Dimensional Properties of Spectral Measures}\label{ss.hdpsm}

While the measure decomposition \eqref{e.measuredecompstandard} is standard, the following refinement turns out to be useful as well. While the decomposition \eqref{e.measuredecompstandard} is obtained by decomposing the measure in question relative to Lebesgue measure and counting measure, one can interpolate between them by considering Hausdorff measures $h^\alpha$, $\alpha \in [0,1]$.

Recall that the $\alpha$-dimensional Hausdorff measure $h^\alpha$ is defined by
$$
h^\alpha (S) = \lim_{\delta \to 0} \; \inf_{\substack{\text{$\delta$-covers}\\ \text{of $S$}}} \; \sum |I_m|^\alpha,
$$
where $S \subseteq \R$ is a Borel set and a $\delta$-cover is a countable collection of intervals $I_m$ of length bounded by $\delta$ such that the union of these intervals contains the set in question. Note that $h^1$ coincides with Lebesgue measure and $h^0$ is the counting measure.

Given a finite Borel measure $\mu$ on $\R$ and $\alpha \in [0,1]$, the upper $\alpha$-derivative of $\mu$  is defined by
$$
D_\mu^\alpha(E) = \limsup_{\varepsilon \downarrow 0} \frac{\mu( ( E - \varepsilon, E + \varepsilon ) )}{(2\varepsilon)^\alpha}.
$$
Denote $T_f = \{ E : D^\alpha_\mu (E) < \infty \}$ and $T_\infty = \{ E : D^\alpha_\mu (E) = \infty \}$. Then \cite{rogers},

\begin{theorem}
We have $h^\alpha (T_\infty) = 0$ and $\mu(S \cap T_f) = 0$ for any $S$ with $h^\alpha(S) = 0$.
\end{theorem}

This suggests the following decomposition of $\mu$. Let $\mu_{\alpha c} (\cdot)= \mu(\cdot \cap T_f)$ and $\mu_{\alpha s} (\cdot)= \mu(\cdot \cap T_\infty)$. Then,
\begin{equation}\label{e.alphadecmoposition}
\mu = \mu_{\alpha c} + \mu_{\alpha s}.
\end{equation}
We say that $\mu$ is $\alpha$-continuous if $\mu_{\alpha s} = 0$ and $\alpha$-singular if $\mu_{\alpha c} = 0$. We also say that $\mu$ is zero-dimensional if $\mu_{\alpha c} = 0$ for every $\alpha > 0$.

The following result was proved by Last \cite{l}. Similar bounds were shown earlier under more restrictive assumptions by Guarneri \cite{G89} and Combes \cite{C93}.

\begin{theorem}\label{t.gclbound}
If $\mu_{\psi_0,\alpha c} \not= 0$, then for every $p > 0$, there is a constant $C = C(H,\psi_0,p)$ such
that
\begin{equation}\label{e.gclbound}
\frac{1}{2T} \int_{-T}^T \sum_{n \in \Z} |n|^p | \langle \delta_n , e^{-itH} \psi_0 \rangle |^2 \, dt \ge C T^{\alpha p}
\end{equation}
\end{theorem}

If the assumption $\mu_{\psi_0,\alpha c} \not= 0$ holds with some $\alpha > 0$, \eqref{e.gclbound} gives a lower bound on how fast (at least a part of) the evolution leaves a finite set. In this sense Theorem~\ref{t.gclbound} provides a quantitative counterpart to the RAGE Theorem.

It should be noted that these are strictly one-sided bounds. That is, growth of the left-hand side in \eqref{e.gclbound} does not imply any continuity properties for $\mu_{\psi_0}$.

There is a convenient way to capture the essence of the estimate \eqref{e.gclbound} via the introduction of so-called transport exponents. If we denote, for $p > 0$,
$$
X^p(t) = \sum_{n \in \Z} |n|^p | \langle \delta_n , e^{-itH} \psi_0 \rangle |^2
$$
and
$$
\tilde X^p(T) = \frac{1}{2T} \int_{-T}^T \sum_{n \in \Z} |n|^p | \langle \delta_n , e^{-itH} \psi_0 \rangle |^2 \, dt,
$$
describing the $p$-th moment of the position operator and a corresponding time-averaged quantity, then the associated transport exponents are given by
\begin{align*}
\beta^+(p) & = \limsup_{t \to \infty} \frac{\log X^p(t)}{p \log t}, \\
\beta^-(p) & = \liminf_{t \to \infty} \frac{\log X^p(t)}{p \log t}, \\
\tilde \beta^+(p) & = \limsup_{T \to \infty} \frac{\log \tilde X^p(T)}{p \log T}, \\
\tilde \beta^-(p) & = \liminf_{T \to \infty} \frac{\log \tilde X^p(T)}{p \log T}.
\end{align*}
It can be shown that each of these four functions of $p \in (0,\infty)$ takes values in $[0,1]$ and is non-decreasing in $p$; compare \cite{DT10}. Consequently, the limits
\begin{align*}
\alpha^+_u & = \lim_{p \to \infty} \beta^+(p), \quad \alpha^-_u = \lim_{p \to \infty} \beta^-(p) \\
\alpha^+_l & = \lim_{p \to 0} \beta^+(p), \quad \; \; \alpha^-_l = \lim_{p \to 0} \beta^-(p) \\
\tilde \alpha^+_u & = \lim_{p \to \infty} \tilde \beta^+(p), \quad \; \tilde \alpha^-_u = \lim_{p \to \infty} \tilde \beta^-(p) \\
\tilde \alpha^+_l & = \lim_{p \to 0} \tilde \beta^+(p), \quad \; \; \, \tilde \alpha^-_l = \lim_{p \to 0} \tilde \beta^-(p)
\end{align*}
exist and belong to $[0,1]$. Note that the dependence of these transport exponents on the initial state $\psi_0$ is left implicit.

The Guarneri-Combes-Last bound \eqref{e.gclbound} can then be succinctly stated as follows:
\begin{equation}\label{e.gclbound2}
\tilde \alpha^-_l \ge \dim_H \mu_{\psi_0},
\end{equation}
where $\dim_H \mu$ denotes the upper Hausdorff dimension of a measure $\mu$, given by $\dim_H \mu = \inf \{ \dim_H S : \mu(\R \setminus S) = 0 \}$. Replacing Hausdorff dimension with packing dimension, one can estimate the time-averaged ``$+$'' quantities. Indeed, Guarneri and Schulz-Baldes showed in \cite{GSB99} that
\begin{equation}\label{e.gsbbound}
\tilde \alpha^+_l \ge \dim_P \mu_{\psi_0},
\end{equation}

\subsection{The Schr\"odinger Operator}\label{ss.schrodingeroperator}

Up to this point it has only been assumed that $H$ is a bounded self-adjoint operator in $\ell^2(\Z)$. While all the results mentioned so far hold in this general case, the fact of the matter is that the operator $H$ appearing in the time-dependent Schr\"odinger equation \eqref{e.timedepse} is not some arbitrary operator, but rather a very specific operator, namely the one that arises by quantization from the Hamiltonian (the total energy function) in classical mechanics. Specifically, if $V : \Z \to \R$ is bounded, then $H$ acts on $\psi \in \ell^2(\Z)$ as follows:
\begin{equation}\label{e.oper}
[H \psi](n) = \psi(n+1) + \psi(n-1) + V(n) \psi(n).
\end{equation}
The function $V$ is called the potential, and it models the medium to which the quantum state $\psi$ is exposed. The quantization procedure associates with the potential energy function $V$ the operator that acts by multiplication with $V$. The term
\begin{equation}\label{e.laplace}
[\Delta \psi](n) = \psi(n+1) + \psi(n-1)
\end{equation}
in $H$ is the discrete Laplacian, and it arises from the kinetic energy by quantization.\footnote{It would be more accurate to also include the term $-2\psi(n)$ on the right-hand side of \eqref{e.laplace}, but it is a standard convention to drop this term and essentially subsume it in the energy. Similarly, it would also be more accurate to consider $-\Delta$ in \eqref{e.oper}, rather than $\Delta$, but this is another standard convention, which we follow here as well.}

Since $V$ is bounded and real-valued, it is easy to see that $H$ is a bounded self-adjoint operator. Thus, as soon as $V$ is fixed, we have an associated quantum evolution, given by \eqref{e.timedepse}, and spectral measures, obeying \eqref{e.spectraltheorem}. Our goal in this paper is to discuss this evolution and these measures in cases where $V$ is generated by some discrete-time dynamical system. We will describe this setting in detail in Subsection~\ref{ss.defandexamples}. Through the end of this section, however, we will consider the general case, where $V$ is merely assumed to be bounded and real-valued.

Thus, given some initial state $\psi_0$, we are interested in the time-evolution \eqref{e.timedepse}. As discussed in Subsections~\ref{ss.rage} and \ref{ss.hdpsm}, we would therefore like to understand the continuity properties of the spectral measure $\mu_{\psi_0}$. In many cases this type of analysis is carried out for all possible initial states $\psi_0$ simultaneously. More precisely, one tries to identify a universal measure with respect to which all spectral measures are absolutely continuous. In cases where there is a cyclic vector, one can just take the spectral measure associated with this vector. Here, a vector $\psi$ is called cyclic if the linear span of $\{ H^m \psi : m \ge 0 \}$ is dense in $\ell^2(\Z)$. Unfortunately, in our case at hand, a Schr\"odinger operator in $\ell^2(\Z)$, there is in general no cyclic vector.\footnote{If, on the other hand, one considers Schr\"odinger operators on the half line $\Z_+ = \{ 1, 2, 3, \ldots \}$, then the vector $\delta_1$ is cyclic and one can consider the universal spectral measure $\mu_{\delta_1}$.} There is, however, a canonical choice of a universal spectral measure. Indeed, if we consider
\begin{equation}\label{e.canonicalsp}
\mu_\mathrm{univ} = \mu_{\delta_0} + \mu_{\delta_1},
\end{equation}
then every spectral measure is absolutely continuous with respect to this measure. That is, if $B$ is a Borel set with $\mu_\mathrm{univ}(B) = 0$, then $\mu_\psi(B) = 0$ for every $\psi \in \ell^2(\Z)$.

With the universal spectral measure $\mu_\mathrm{univ}$ one can conveniently describe the sets $\sigma(H)$, $\sigma_\mathrm{ac}(H)$, $\sigma_\mathrm{sc}(H)$, $\sigma_\mathrm{pp}(H)$. Let us denote the topological support of a measure $\mu$ of $\mathrm{supp} \, \mu$, that is, $\mathrm{supp} \, \mu$ is the complement of the largest open set that has zero weight with respect to $\mu$.

\begin{theorem}\label{t.muunivandspectrum}
We have $\sigma(H) = \mathrm{supp} \, \mu_\mathrm{univ}$ and $\sigma_\mathrm{ac}(H) = \mathrm{supp} \, \mu_\mathrm{univ, \, ac}$, $\sigma_\mathrm{sc}(H) = \mathrm{supp} \, \mu_\mathrm{univ, \, sc}$, $\sigma_\mathrm{pp}(H) = \mathrm{supp} \, \mu_\mathrm{univ, \, pp}$.
\end{theorem}

There is an important alternative way to view $\mu_\mathrm{univ, \, pp}$. A minimal support of this measure is given by the set of eigenvalues of $H$. That is, the complement of the set of eigenvalues has zero weight with respect to $\mu_\mathrm{univ, \, pp}$, and each eigenvalue has positive weight with respect to $\mu_\mathrm{univ, \, pp}$. This has the following consequence:

\begin{coro}\label{c.sigmappandeigenvalues}
The set $\sigma_\mathrm{pp}(H)$ is equal to the closure of the set of eigenvalues of $H$, that is,
$$
\sigma_\mathrm{pp}(H) = \overline{ \left\{ E \in \R : \exists \psi \in \ell^2(\Z) \setminus \{ 0 \} \text{ such that } H\psi = E\psi \right\} }.
$$
\end{coro}

We end this subsection with a useful formula connecting solutions and Green functions. For $[n_1,n_2] = \{ n \in \Z : n_1 \le n \le n_2 \}$, denote by $H_{[n_1,n_2]}$ the restriction of $H$ to this interval, that is, $H_{[n_1,n_2]} = P_{[n_1,n_2]} H P_{[n_1,n_2]}^*$, where $P_{[n_1,n_2]} : \ell^2(\Z) \to \ell^2([n_1,n_2])$ is the canonical projection and $P_{[n_1,n_2]}^* : \ell^2([n_1,n_2]) \to \ell^2(\Z)$ is the canonical embedding.

Moreover, for
$E \not\in \sigma(H_{[n_1,n_2]})$ and $n,m \in [n_1,n_2]$, let
$$
G_{[n_1,n_2]}(n,m;E) := \langle \delta_n , \left( H_{[n_1,n_2]} - E \right)^{-1} \delta_m \rangle.
$$
Then, the following formula holds.

\begin{lemma}\label{l.eigengreenlocal}
Suppose $n \in [n_1,n_2] \subset \Z$ and $u$ is a solution of the difference equation $Hu=Eu$. If $E \not\in \sigma(H_{[n_1,n_2]})$ and $n \in [n_1,n_2]$, then
$$
u(n) = - G_{[n_1,n_2]}(n,n_1;E) u(n_1 - 1) - G_{[n_1,n_2]}(n,n_2;E) u(n_2 + 1).
$$
\end{lemma}

\subsection{Solutions of the Time-Independent Schr\"odinger Equation}\label{ss.tise}

Suppose $V : \Z \to \R$ is bounded and consider the bounded self-adjoint operator $H$ from \eqref{e.oper} and the associated universal spectral measure $\mu_\mathrm{univ}$ from \eqref{e.canonicalsp}. In this subsection we consider the difference equation
\begin{equation}\label{e.eve}
u(n+1) + u(n-1) + V(n) u(n) = E u(n)
\end{equation}
for $E \in \C$, and relate the behavior of these solutions to properties of the measure $\mu_\mathrm{univ}$. The equation \eqref{e.eve} is called the \emph{time-independent Schr\"odinger equation}.

While \eqref{e.eve} looks like the eigenvalue equation for the operator $H$, we emphasize that the solutions of \eqref{e.eve} we consider do not have to belong to $\ell^2(\Z)$. Thus, for each $E \in \C$, the solutions of \eqref{e.eve} form a two-dimensional vector space. Indeed, as soon as we fix two consecutive values of $u$, the whole solution is completely determined by \eqref{e.eve}. For example, suppose we fix $u(0)$ and $u(1)$, then any $u(n)$ is obtained by solving the difference equation ``from the origin to $n$.'' This can be formalized using transfer matrices as follows. If we set
\begin{equation}\label{e.onestepmats}
T(m;E) = \begin{pmatrix} E - V(m) & -1 \\ 1 & 0 \end{pmatrix}
\end{equation}
and
\begin{equation}\label{e.transmatrices}
A(n;E) =
\begin{cases}
T(n;E) \times \cdots \times T(1;E)             & n \ge 1 \\
I                                              & n = 0   \\
T(n+1;E)^{-1} \times \cdots \times T(0;E)^{-1} & n \le -1,
\end{cases}
\end{equation}
then $u$ solves \eqref{e.eve} for every $n \in \Z$ if and only if
\begin{equation}\label{e.meve}
\begin{pmatrix} u(n+1) \\ u(n) \end{pmatrix} = A(n;E) \begin{pmatrix}
u(1) \\ u(0) \end{pmatrix}
\end{equation}
for every $n \in \Z$. More generally, we denote the matrix that maps solution data from $m$ to $n$ by $A(n,m;E)$, that is,
\begin{equation}\label{e.meve2}
\begin{pmatrix} u(n+1) \\ u(n) \end{pmatrix} = A(n,m;E) \begin{pmatrix}
u(m+1) \\ u(m) \end{pmatrix}.
\end{equation}
This matrix is also given by a suitable product of one-step transfer matrices similar to \eqref{e.transmatrices}.

A convenient way to choose a basis of the solution space of \eqref{e.eve} is to prescribe a pair of initial conditions. For $\theta \in (-\frac{\pi}{2}, \frac{\pi}{2}]$, consider the pair $u_\theta, v_\theta$ of solutions of \eqref{e.eve} satisfying
\begin{equation}\label{e.evesolinicon}
\begin{pmatrix} u_\theta(1)  & v_\theta(1) \\ u_\theta(0) & v_\theta(0) \end{pmatrix} = \begin{pmatrix} \cos \theta & \sin \theta \\ - \sin \theta & \cos \theta \end{pmatrix}.
\end{equation}
Clearly, $u_\theta, v_\theta$ are linearly independent and hence form a basis of $\{ u : \Z \to \C : u \text{ solves \eqref{e.eve}} \}$. The relations \eqref{e.meve} and \eqref{e.evesolinicon} imply
\begin{equation}\label{e.transmatalternative}
\begin{pmatrix} u_\theta(n+1)  & v_\theta(n+1) \\ u_\theta(n) & v_\theta(n) \end{pmatrix} = A(n;E) \begin{pmatrix} \cos \theta & \sin \theta \\ - \sin \theta & \cos \theta \end{pmatrix}.
\end{equation}
Note that the determinant of the right-hand side is $1$, so that
\begin{equation}\label{e.wronskianconservation}
u_\theta(n+1) v_\theta(n) - v_\theta(n+1) u_\theta(n) = 1 \quad \text{ for every } n \in \Z.
\end{equation}

For $\theta = 0$, \eqref{e.transmatalternative} becomes
\begin{equation}\label{e.transmatalternative2}
A(n;E) = \begin{pmatrix} u_0(n+1)  & v_0(n+1) \\ u_0(n) & v_0(n) \end{pmatrix},
\end{equation}
which shows that the entries of the transfer matrices are given by suitable solutions of \eqref{e.eve} (namely the Dirichlet solution and the Neumann solution). In particular we can obtain estimates for the norms of transfer matrices if we have estimates for the norms of solutions. For example, if \eqref{e.eve} admits exponentially growing or decaying solutions, then the transfer matrix norms must grow exponentially (this is clear in the case of an exponentially growing solution; in the case of an exponentially decaying solution, apply the Wronskian conservation law \eqref{e.wronskianconservation}). There is a certain converse to this correspondence in the setting of exponential growth, which follows from the following abstract result proved by Ruelle in \cite{R79}.

\begin{theorem}\label{t.ruelle}
Suppose $T_n \in \mathrm{SL}(2,\R)$ obey
$$
\lim_{n \to \infty} \frac{1}{n} \| T_n \| = 0
$$
and
$$
\lim_{n \to \infty} \frac{1}{n} \log \| T_n \cdots T_1 \| = L > 0.
$$
Then there exists a one-dimensional subspace $V \subset \R^2$ such that
$$
\lim_{n \to \infty} \frac{1}{n} \log \| T_n \cdots T_1 v \| = - L \quad \text{ for } v \in V \setminus \{0\}
$$
and
$$
\lim_{n \to \infty} \frac{1}{n} \log \| T_n \cdots T_1 v \| = L \quad \text{ for } v \not\in V.
$$
\end{theorem}

Applied to the transfer matrices, Theorem~\ref{t.ruelle} gives that the condition $\lim_{n \to \infty} \frac{1}{n} \log \| A(n;E) \| = L > 0$ implies that there exists (up to a constant multiple) exactly one solution that decays exponentially, in fact at the same rate, at $\infty$, while all other solutions increase exponentially, also at the same rate. The same statement holds near $- \infty$. For dynamically defined potentials, the condition $\lim_{n \to \infty} \frac{1}{n} \log \| A(n;E) \| = L > 0$ will hold almost surely whenever the Lyapunov exponent $L(E)$ is positive; see our discussion below. Thus, the result just discussed is relevant in this situation.

\bigskip

An energy $E$ is an eigenvalue of $H$ if and only if \eqref{e.eve} admits a non-trivial $\ell^2$ solution $u$. There is a simple way of excluding the existence of $\ell^2$, and in fact decaying, solutions $u$ that relies on local (almost) repetitions of the potential and which goes back to Gordon \cite{G76}; see also \cite{D00, DP86, S87}. The history of this lemma is discussed in \cite{G16}. The most elementary statement of this kind is given in the following lemma.

\begin{lemma}\label{l.threeblockgordon}
Suppose the potential $V$ obeys $V(m+p) = V(m)$, $-p \le m \le p-1$. Then, every solution $u$ of \eqref{e.eve} satisfies
$$
\max \left\{ \left\| \begin{pmatrix} u(2p+1) \\ u(2p) \end{pmatrix} \right\|, \left\| \begin{pmatrix} u(p+1) \\ u(p) \end{pmatrix}  \right\| , \left\| \begin{pmatrix} u(-p+1) \\ u(-p) \end{pmatrix} \right\| \right\} \ge \frac{1}{2} \left\| \begin{pmatrix} u(1) \\ u(0) \end{pmatrix} \right\|.
$$
\end{lemma}

This lemma follows quickly from the Cayley-Hamilton theorem, which gives $A(p;E)^2 - [\mathrm{Tr} A(p;E)] A(p;E) + I = 0$. Applying this to either $(u(-p+1),u(-p))^T$ or $(u(1),u(0))^T$, depending on whether $|\mathrm{Tr} A(p;E)| > 1$ or $|\mathrm{Tr} A(p;E)| \le 1$, implies the lemma. Notice that as an immediate consequence, no non-trivial solution $u$ of \eqref{e.eve} can decay at both $\pm \infty$ when Lemma~\ref{l.threeblockgordon} applies for arbitrarily large $p$ (i.e., there is a sequence $p_k \to \infty$ for which one has the required three-block symmetries).

In fact, the same proof shows that if one is able to control the trace of the transfer matrix, one can work exclusively on one half-line:

\begin{lemma}\label{l.twoblockgordon}
Suppose the potential $V$ obeys $V(m+p) = V(m)$, $0 \le m \le p-1$. Then, every solution $u$ of \eqref{e.eve} satisfies
$$
\max \left\{ \left\| \begin{pmatrix} u(2p+1) \\ u(2p) \end{pmatrix} \right\|, \left\| \begin{pmatrix} u(p+1) \\ u(p) \end{pmatrix}  \right\| \right\} \ge \frac{1}{2 \max \{ | \mathrm{Tr} A(p;E) | , 1\}}  \left\| \begin{pmatrix} u(1) \\ u(0) \end{pmatrix} \right\|.
$$
\end{lemma}

This is useful in some situations, for example when studying the Fibonacci Hamiltonian and its generalizations, which are discussed in later sections.

Another useful remark is that one does not need exact repetitions, and one can instead allow errors that are (super-)exponentially small in the local period: A bounded potential $V : \Z \to \R$ is called a \emph{Gordon potential} if there are positive integers $p_k \to \infty$ such that
$$%\begin{equation}\label{f.gordondef2}
\forall \, C > 0 : \lim_{k \to \infty} \max_{1 \le n \le p_k} |V(n) - V(n \pm p_k)| C^{p_k} = 0.
$$%\end{equation}

\begin{lemma}\label{l.gordonpotential}
Suppose $V$ is a Gordon potential. Then, the operator $H$ has purely continuous spectrum. More precisely, for every $E \in \R$ and every solution $u$ of \eqref{e.eve}, we have
$$%\begin{equation}\label{f.gordonpotsolest}
\limsup_{|n| \to \infty} \left\| \begin{pmatrix} u(n+1) \\ u(n) \end{pmatrix} \right\| \ge \frac12 \left\| \begin{pmatrix} u(1) \\ u(0) \end{pmatrix} \right\|.
$$%\end{equation}
\end{lemma}

%\bigskip

Let us now discuss how to characterize the spectrum and the spectral type of $H$ in terms of solutions of \eqref{e.eve}. We say that $E$ is a generalized eigenvalue of $H$ if \eqref{e.eve} has a non-trivial solution, called the corresponding generalized eigenfunction, that satisfies
\begin{equation}\label{f.geneigenf}
|u(n)| \le C(1 + |n|)^\delta
\end{equation}
for suitable finite constants $C,\delta > 0$, and every $n \in \Z$.

\begin{theorem}\label{t.genef}
{\rm (a)} Every generalized eigenvalue of $H$ belongs to $\sigma(H)$ and hence is necessarily real.

{\rm (b)} Fix $\delta > \frac12$. Then, for $\mu_\mathrm{univ}$-almost every $E \in \R$, there exists a generalized eigenfunction satisfying
\eqref{f.geneigenf}.

{\rm (c)} The spectrum of $H$ is given by the closure of the set of generalized eigenvalues of $H$.
\end{theorem}

This theorem shows that the spectrum $\sigma(H)$ as a set is completely determined by the behavior of the solutions of \eqref{e.eve}. In fact the sets $\sigma_\mathrm{ac}(H)$, $\sigma_\mathrm{sc}(H)$, $\sigma_\mathrm{pp}(H)$ can also be described in terms of solutions. The case of $\sigma_\mathrm{pp}(H)$ is the easiest. Recall from Corollary~\ref{c.sigmappandeigenvalues} that it is simply given by the closure of the eigenvalues of $H$. Now, $E$ is an eigenvalue of $H$ if and only if \eqref{e.eve} has a non-trivial solution that is square-summable at both $\pm \infty$. Thus, square-summability is the way to discriminate between supports of $\mathrm{supp} \, \mu_\mathrm{univ, \, c}$ and $\mathrm{supp} \, \mu_\mathrm{univ, \, pp}$. This raises the following natural question: Is there a similar way to discriminate between supports of $\mathrm{supp} \, \mu_\mathrm{univ, \, ac}$ and $\mathrm{supp} \, \mu_\mathrm{univ, \, s}$?

It turns out that there is such a way, and the following important definition is due to Gilbert and Pearson. A non-trivial solution $u$ of \eqref{e.eve} is called \textit{subordinate at} $+ \infty$ if
$$
\lim_{L\to\infty} \frac{\|u\|_L}{\|v\|_L} = 0
$$
for any linearly independent solution $v$ of \eqref{e.eve}, where $\|\cdot\|_L$ denotes the norm of the solution over a lattice interval of length $L$. That is, for $L > 0$ we define
$$
\|u\|_L \equiv\left[\sum_{n=1}^{\lfloor L \rfloor}|u(n)|^2+ (L-\lfloor L \rfloor)|u(\lfloor L \rfloor+1)|^2\right]^{\frac{1}{2}},
$$
where $\lfloor L \rfloor$ denotes the integer part of $L$. Subordinacy at $-\infty$ is defined similarly, by considering solutions on $[L,0)$ for $L < 0$ and sending $L \to -\infty$. Finally, we say that a solution $u$ of \eqref{e.eve} is called \textit{subordinate} if it is subordinate at both $\pm \infty$.

The notion of subordinacy turns out to be the counterpart to square-summability and provides the desired split between the supports of $\mathrm{supp} \, \mu_\mathrm{univ, \, ac}$ and $\mathrm{supp} \, \mu_\mathrm{univ, \, s}$; see \cite{Gil89, GP87}.

\begin{theorem}\label{t.gilbert}
{\rm (a)} The singular part $\mu_\mathrm{univ, \, s}$ is supported by
$$
S = \{ E \in \R : \text{ \eqref{e.eve} has a subordinate solution} \}.
$$
That is, the complement of $S$ has zero weight with respect to $\mu_\mathrm{univ, \, s}$.

{\rm (b)} The set $N = N_+ \cup N_-$, where
$$
N_\pm = \{ E \in \R : \text{ \eqref{e.eve} has no solution that is subordinate at } \pm \infty \},
$$
is an essential support of $\mu_\mathrm{univ, \, ac}$. That is, $\mu_\mathrm{univ, \, ac} (\R \setminus N) = 0$, and for any measurable set $A$ with $\mu_\mathrm{univ, \, ac}(\R \setminus A) = 0$, we have $\mathrm{Leb}(N \setminus A) = 0$.
\end{theorem}

Combining this result with the earlier result for $\mu_\mathrm{univ, \, pp}$, we can state the following corollary.

\begin{coro}\label{c.solutionsupportsofspectra}
We have
\begin{align*}
\sigma_\mathrm{ac}(H) & = \overline{\{ E \in \R : \text{at $\infty$ or $-\infty$, \eqref{e.eve} has no subordinate solution} \}}^\mathrm{ess},  \\
\sigma_\mathrm{sc}(H) & \subseteq \overline{\{ E \in \R : \text{\eqref{e.eve} has a subordinate solution, which is not square-summable} \}} , \\
\sigma_\mathrm{pp}(H) & = \overline{\{ E \in \R : \text{\eqref{e.eve} has a non-trivial square-summable solution} \}}.
\end{align*}
\end{coro}

Note that we unfortunately cannot claim equality in the description of $\sigma_\mathrm{sc}(H)$. Of course in cases where $\{ E \in \R : \text{at $\infty$ or $-\infty$, \eqref{e.eve} has no subordinate solution} \}$ and $\{ E \in \R : \text{\eqref{e.eve} has a non-trivial square-summable solution} \}$ are both empty, we must have equality,\footnote{The right-hand side is always contained in $\sigma(H)$, and by the two sets above being empty, we must have $\sigma_\mathrm{sc}(H) = \sigma(H)$; this gives the reverse inclusion.} but in general there is no mechanism that deduces from the presence of subordinate non-square-summable solutions the presence of singular continuous spectrum.

An important special case of an energy $E$ in $N_\pm$ is where the transfer matrix $A(n;E)$ is bounded as a function of $n \in \Z_\pm$. This case occurs in many applications, so for the sake of easy reference, we formulate the following corollary; compare, for example, \cite{B91, JL99, S96, S92}.

\begin{coro}\label{c.boundedsolac}
Denote
$$
B_\pm = \{ E \in \R: \|A(n;E)\| \text{ is uniformly bounded for } n \in \Z_\pm \}.
$$
Then, $\mu_\mathrm{univ}$ is purely absolutely continuous on $B = B_+ \cup B_-$.
\end{coro}

On the one hand, it is not hard to see that $B_\pm \subseteq N_\pm$, so that Corollary~\ref{c.boundedsolac} follows from Theorem~\ref{t.gilbert}. On the other hand, Corollary~\ref{c.boundedsolac} is also an explicit consequence of Theorem~\ref{t.wholelinesubordinacy} below.

Absolute continuity of the spectral measure on some set is very often established through Corollary~\ref{c.boundedsolac}, that is, by showing boundedness of all solutions for the energies in question. As a consequence, it was quite tempting to conjecture that boundedness of solutions is not only sufficient, but also necessary for absolute continuity. This conjecture was often referred to as the Schr\"odinger conjecture. It was recently disproved by Avila in \cite{A14f}. We will say more about this work in a later section.

\bigskip

As discussed in Subsection~\ref{ss.hdpsm}, in order to prove quantitative transport bounds, we consider the decomposition \eqref{e.alphadecmoposition} of a spectral measure or the universal spectral measure with $\alpha$ as large as possible so that the $\alpha$-continuous piece of the measure is non-zero. Thus we ask if it is possible to study this decomposition via solutions of \eqref{e.eve}. It turns out that it is possible to introduce a notion of an $\alpha$-subordinate solution of \eqref{e.eve}, which coincides with the notion of a subordinate solution for $\alpha = 1$, such that the decomposition \eqref{e.alphadecmoposition} of the universal spectral measure can be related to those sets of energies where such solutions do or do not exist. This definition is due to Jitomirskaya and Last. For $\alpha \in [0,1]$, we say that a non-trivial solution $u$ of \eqref{e.eve} is $\alpha$\textit{-subordinate at} $+ \infty$ if
$$
\liminf_{L\to\infty} \frac{\|u\|_L^{2-\alpha}}{\|v\|_L^\alpha} = 0
$$
for any linearly independent solution $v$ of \eqref{e.eve}. Note that for $\alpha = 1$, we recover a weak form of the definition of subordinacy at $+ \infty$, while for $\alpha = 0$, we recover the definition of square-summability at $+ \infty$. Again, $\alpha$-subordinacy at $-\infty$ of a solution $u$ is defined similarly on the left half-line. In analogy to Theorem~\ref{t.gilbert} one would hope that
$$
S_\alpha = \{ E \in \R : \text{ \eqref{e.eve} has a solution that is $\alpha$-subordinate at both } \pm \infty \}
$$
is a support of $\mu_\mathrm{univ, \, \alpha s}$ and that $S_\alpha$ has zero weight with respect to $\mu_\mathrm{univ, \, \alpha c}$. Alas, this is at present unknown, even though it is a natural conjecture to make. What is known, on the other hand, is the following somewhat weaker statement \cite{DKL00}:

\begin{theorem}\label{t.wholelinesubordinacy}
Suppose $B \subset \R$ is a bounded Borel set. Assume that there are constants $\gamma_1, \gamma_2$ such that for every $E \in B$, there are constants $C_1(E), C_2(E)$ so that every solution $u$ of \eqref{e.eve} that is normalized in the sense that $|u(0)|^2 + |u(1)|^2 = 1$ obeys the estimate
\begin{equation}\label{e.powerlawbounds}
C_1(E) L^{\gamma_1} \le \|u\|_L \le C_2(E) L^{\gamma_2}
\end{equation}
for $L > 0$ sufficiently large. Set
$$
\alpha = \frac{2 \gamma_1}{\gamma_1 + \gamma_2}.
$$
Then the restriction of $\mu_\mathrm{univ}$ is purely $\alpha$-continuous, that is, $\mu_\mathrm{univ, \, \alpha s} (B) = 0$.
\end{theorem}

Let us add a few remarks. First of all, every $E$ in $B$ is a generalized eigenvalue and hence belongs to the spectrum of $H$. Moreover, as was the case in the description of $\mu_\mathrm{univ, \, ac}$ in Theorem~\ref{t.gilbert}.(b), information on one half-line is sufficient. Indeed, the conclusion of Theorem~\ref{t.wholelinesubordinacy} holds true if the assumptions about the solutions are phrased in terms of conditions on the left half-line. Finally, while the validity of the power-law bounds \eqref{e.powerlawbounds} is formally speaking stronger than the absence of solutions that are $\alpha$-subordinate at $+\infty$, this is the usual way in which this absence is established.

\subsection{Further Ways to Establish Transport Bounds}\label{ss.dtbounds}

While the combination of $\alpha$-subordinacy and the Guarneri-Combes-Last transport estimate provide a nice one-two punch, it is often quite difficult to actually establish the required power-law solution estimates \eqref{e.powerlawbounds}. Moreover, the results are strictly one-sided in the sense that transport may be fast even if $\mu_\mathrm{univ}$ is very singular. There are some extreme cases such as examples displaying almost ballistic transport and pure point spectral measures. In these cases fast transport cannot be established via the spectral continuity route and hence a different method is required.

A method to show transport bounds without resorting to spectral measures at all was suggested by Damanik and Tcheremchantsev. While this method has gone through several stages of evolution \cite{DST, DT03, DT05, DT07, DT08}, let us state here the simplest version of it.

The starting point is the following simple lemma.

\begin{lemma}\label{l.kkllemma}
$$
\int_0^\infty e^{-2t/T} | \langle \delta_n , e^{-itH} \psi_0 \rangle |^2 \, dt = \frac{1}{2\pi} \int_\R | \langle \delta_n , (H - E - \tfrac{i}{T})^{-1} \psi_0 \rangle |^2 \, dE.
$$
\end{lemma}

The proof uses the spectral theorem twice to write the inner products as integrals over suitable spectral measures. This formula suggests replacing the Ces\`aro time-averages considered in Subsection~\ref{ss.hdpsm} with the following average,
$$
\frac{2}{T} \int_0^\infty e^{-2t/T} | \langle \delta_n , e^{-itH} \psi_0 \rangle |^2 \, dt,
$$
which by Lemma~\ref{l.kkllemma} is equal to
$$
\frac{1}{\pi T} \int_\R | \langle \delta_n , (H - E - \tfrac{i}{T})^{-1} \psi_0 \rangle |^2 \, dE.
$$
Thus, the modified time-averaged moments of the position operator are
$$
\tilde X^p(T) = \frac{2}{T} \int_0^\infty e^{-2t/T} \sum_{n \in \Z} |n|^p | \langle \delta_n , e^{-itH} \psi_0 \rangle |^2 \, dt,
$$
and one can now define the resulting modified time-averaged transport exponents $\tilde \beta^\pm(p)$ as before. The fact of the matter is that they are actually not modified at all; see \cite[Section~2.6]{DT10}. That is, the values of $\tilde \beta^\pm(p)$ coincide in the two cases, and hence when studying these transport exponents, one can choose the underlying way of time-averaging that is more convenient.

Why is this change of perspective useful? Consider the special case $\psi_0 = \delta_0$ for simplicity. Notice that in this case
$$
\langle \delta_n , (H - E - \tfrac{i}{T})^{-1} \psi_0 \rangle = \langle \delta_n , (H - E - \tfrac{i}{T})^{-1} \delta_0 \rangle
$$
is simply a Green's function entry. There are powerful tools one can use to estimate these quantities, especially in the one-dimensional case we are interested in. In fact, the Green's function can be expressed in terms of solutions of $Hu = (E + \tfrac{i}{T})u$, which in turn can be expressed with the help of transfer matrices. Thus, we can study time-averaged transport directly by estimating transfer matrices! Note, however, that this needs to happen at energies with non-trivial imaginary part.

Pursuing this further, Damanik and Tcheremchantsev established the following result in \cite{DT03}.

\begin{theorem}
Suppose that for some $K, C, \alpha > 0$, the following condition holds: For any $N>0$ large enough, there exists a nonempty Borel set $A(N)\subset \R$ such that
$\mathcal{E}(N) \subset [-K,K]$ and
\begin{equation}\label{maincond}
\|A(n,m;E)\| \le C N^{\alpha} \ \  \forall E \in \mathcal{E}(N),\ \forall \ n,m: |n|\le N, |m| \le N.
\end{equation}
Let $N(T) = T^{1/(1+\alpha)}$ and let $\mathcal{N}(T)$ be the $1/T$-neighborhood of the set $\mathcal{E}(N(T))$:
$$
\mathcal{N}(T)=\{ E \in \R : \exists E' \in \mathcal{E}(N(T)), |E-E'| \le 1/T \}.
$$
Then for all $T > 1$ large enough, the following bound holds:
$$
\sum_{|n| \ge N(T)/2} \frac{2}{T} \int_0^\infty e^{-2t/T} | \langle \delta_n , e^{-itH} \delta_0 \rangle |^2 \, dt \ge \frac{\hat{C}}{T} |\mathcal{N}(T)| N^{1-2\alpha} (T),
$$
where $\hat{C}$ is some uniform positive constant and $|\cdot|$ denotes Lebesgue measure.

In particular, for any $p > 0$, one has the following bound for the time-averaged moments of
the position operator:
$$
\tilde X^p (T) \ge \frac{\hat{C}}{T}|\mathcal{N}(T)| N^{p+1-2 \alpha} (T).
$$
\end{theorem}

This result is useful even if the set of energies for which one can control transfer matrix norms consists of a single element:

\begin{coro}
If
$$
\|A(n,m; E_0)\| \le C(E_0)(|n|+|m|)^{\alpha}
$$
for some $E_0 \in \R$, uniformly in $n,m \in\Z$, then
$$
\tilde \alpha^-_u \ge \frac{1}{1+\alpha}.
$$
\end{coro}

Contrary to relying on spectral continuity properties, this approach can also be used to show dynamical upper bounds. Indeed, the following theorem was shown in \cite{DT07}.

\begin{theorem}
Suppose $H$ is as in \eqref{e.oper}, and $K \ge 4$ is such that $\sigma (H) \subseteq [-K+1,K-1]$. Suppose that, for some $C \in (0,\infty)$ and $\alpha \in (0,1)$, we have
\begin{equation}\label{assumeright}
\int_{-K}^K \left( \max_{1 \le n \le C T^\alpha} \left\| A \left( n;E+ \tfrac{i}{T}
\right) \right\|^2 \right)^{-1} dE = O(T^{-m})
\end{equation}
and
\begin{equation}\label{assumeleft}
\int_{-K}^K \left( \max_{1 \le -n \le C T^\alpha} \left\| A \left( n;E+ \tfrac{i}{T}
\right) \right\|^2 \right)^{-1} dE = O(T^{-m})
\end{equation}
for every $m \ge 1$. Then
\begin{equation}\label{apubound}
\tilde \alpha_u^+ \le \alpha.
\end{equation}
In particular,
\begin{equation}\label{bppbound}
\tilde \beta^+ (p) \le \alpha \quad \text{ for every } p > 0.
\end{equation}
\end{theorem}

The paper \cite{DT08} shows how to obtain analogous upper bounds for non-time-averaged transport exponents in terms of transfer matrix estimates.

\section{Schr\"odinger Operators with Dynamically Defined Potentials}\label{s.3}

\subsection{Basic Definitions and Examples}\label{ss.defandexamples}

Suppose we are given a probability measure space $(\Omega,\mathcal{B},\mu)$. We will usually leave the $\sigma$-algebra $\mathcal{B}$ implicit and just write $(\Omega,\mu)$. Integration with respect to $\mu$ will be denoted by $\E(\cdot)$, that is, if $f \in L^1(\Omega,d\mu)$, then
$$
\E(f) = \int f(\omega) \, d\mu(\omega).
$$

Suppose further that $T: \Omega \to \Omega$ is an \emph{invertible measure-preserving transformation}. That is, $T$ is a one-to-one and onto map so that for every $B \in \mathcal{B}$, we have $T B, T^{-1} B \in \mathcal{B}$ and $\mu(B) = \mu (T B) = \mu( T^{-1} B)$. Conversely, given $(\Omega,\mathcal{B},T)$, a probability measure $\mu$ with the invariance property above is called an \emph{invariant probability measure} for $T$, or just $T$-\emph{invariant}.

\begin{example}\label{x.torusshift}
Translation on a torus: $\Omega$ is the $d$-dimensional torus $\T^d = \R^d / \Z^d$, $\mathcal{B}$ is the Borel $\sigma$-algebra, $\mu$ is normalized Lebesgue measure, denoted by $\mathrm{Leb}$, $T$ is given by a translation, that is,
$$
T(\omega_1 , \ldots, \omega_d) = (\omega_1 + \alpha_1 , \ldots , \omega_d + \alpha_d),
$$
where $\alpha_1 , \ldots , \alpha_d$ are real numbers, which can and will be chosen in the interval $[0,1)$. The interesting case is where $1, \alpha_1 , \ldots , \alpha_d$ are linearly independent over the rational numbers, and we will assume this unless noted otherwise.
\end{example}

\begin{example}\label{x.skewshift}
Skew-shift on a torus: $\Omega$ is the $2$-torus $\T^2$, $\mathcal{B}$ is the Borel $\sigma$-algebra, $\mu$ is $\mathrm{Leb}$, and $T$ is given by
$$
T(\omega_1,\omega_2) = (\omega_1 + \alpha , \omega_1 + \omega_2),
$$
where $\alpha \in (0,1)$ is irrational.
\end{example}

\begin{example}\label{x.hyptoraut}
Hyperbolic toral automorphism: $\Omega = \T^2$, $\mathcal{B}$ is the Borel $\sigma$-algebra, $\mu$ is $\mathrm{Leb}$, and  $T$ is given by
$$
T(\omega_1,\omega_2) = (2\omega_1 + \omega_2 , \omega_1 + \omega_2).
$$
\end{example}

\begin{example}\label{x.sequenceshift}
Shift on a sequence space: Fix some compact interval $I \subset \R$ with the induced topology. Consider the infinite product $\Omega = I^\Z$ with the product topology and the Borel $\sigma$-algebra $\mathcal{B}$. The shift transformation $T : \Omega \to \Omega$ is given by
$$
(T \omega)_n = \omega_{n+1}.
$$
There are many $T$-invariant measures $\mu$. An important class is obtained by taking $\mu = \rho^\Z$, where $\rho$ is a Borel probability measure on $I$.
\end{example}

\begin{example}\label{x.symbolicshift}
Shift on a symbolic sequence space: This is a slight variation of the previous example. Fix a finite {\rm (}or at most countable{\rm )} set $\mathcal{A}$, called the \emph{alphabet}, equipped with the discrete topology. Consider the infinite product $\Omega = \mathcal{A}^\Z$ with the product topology and the Borel $\sigma$-algebra $\mathcal{B}$. The shift transformation $T : \Omega \to \Omega$ is again given by
$$
(T \omega)_n = \omega_{n+1}.
$$
As before, there are many $T$-invariant measures $\mu$ and an important class of examples is obtained by taking $\mu = \rho^\Z$, where $\rho$ is a probability measure on $\mathcal{A}$. Other interesting examples are given by Markov measures.
\end{example}

We say that $(\Omega,\mu,T)$ is \emph{ergodic} if in addition every invariant function is constant. More precisely, this means
that if $f$ is a measurable function on $\Omega$ and $f(\omega) = f(T \omega) = f(T^{-1} \omega)$ for $\mu$-almost every $\omega \in \Omega$, then there is a constant $f_*$ and a set $\Omega_* \subseteq \Omega$ of full $\mu$-measure so that $f(\omega) = f_*$ for every $\omega \in \Omega_*$. Equivalently, every measurable set $E \subset \Omega$ with $T^{-1} E = E$ must satisfy $\mu(E) = 0$ or $1$. It turns out that all the examples listed above are ergodic.

\bigskip

We are now ready to define the central object of interest in this paper, namely Schr\"odinger operators with dynamically defined potentials. Suppose $(\Omega,\mu,T)$ is ergodic and $f : \Omega \to \R$ is measurable and bounded. Define potentials,
\begin{equation}\label{e.ergpotential}
V_\omega (n) = f (T^n \omega), \quad \omega \in \Omega, \; n \in \Z,
\end{equation}
and Schr\"odinger operators on $\mathcal{H} = \ell^2(\Z)$,
\begin{equation}\label{e.ergoper}
[H_\omega \psi](n) = \psi(n+1) + \psi(n-1) + V_\omega(n) \psi(n).
\end{equation}
The family $\{ H_\omega \}_{\omega \in \Omega}$ is called an \emph{ergodic family of Schr\"odinger operators}. Our goal is to study the spectral properties of the operators $H_\omega$. The canonical spectral measure of $H_\omega$ will be denoted by $\mu_\omega$, that is,
$$
\int_\R \frac{d\mu_\omega(E')}{E' - E} = \langle \delta_0 , (H_\omega - E)^{-1} \delta_0 \rangle + \langle \delta_1 , (H_\omega - E)^{-1} \delta_1 \rangle
$$
for $E \in \C_+$. The Lebesgue decomposition of $\mu_\omega$ will be denoted by
$$
\mu_\omega = \mu_{\omega, \mathrm{ac}} + \mu_{\omega, \mathrm{sc}} + \mu_{\omega, \mathrm{pp}}.
$$

\subsection{Invariance of the Spectrum and the Spectral Type}

A central result of the theory says that the spectrum of $H_\omega$ is non-random in the sense that this set is actually independent of $\omega$ for a full-measure set of $\omega$'s. This theorem was shown by Pastur in 1980; see \cite{P80}.

\begin{theorem}\label{t.pastursetofev}
Given an ergodic family $\{ H_\omega \}_{\omega \in \Omega}$, there exists a set $\Sigma \subseteq \R$ such that for $\mu$-almost every $\omega$, $\sigma(H_\omega) = \Sigma$ and $\sigma_{{\rm disc}}(H_\omega) = \emptyset$. Moreover, for every $E$, $\mu \left( \{ \omega : E \text{ is an eigenvalue of } H_\omega) \} \right) = 0$.
\end{theorem}

The idea behind the proof is simple. By ergodicity, any $T$-invariant measurable function is almost everywhere constant. One may, for example, consider the function that associates, for some fixed open energy interval $J \subset \R$, the dimension of the range of the associated spectral projection, that is, $\mathrm{Tr} \, \chi_J(H_\omega)$. This dimension is zero if and only if $J$ does not intersect the spectrum of $H_\omega$. Thus, almost everywhere constancy implies that either $J$ almost surely does intersect the spectrum, or that it almost surely does not. Varying $J$ in a countable way (choose rational endpoints, for example) then allows one to conclude the almost sure constancy of the spectrum.

It is possible to modify this way of reasoning somewhat to focus on the partial spectra. This allowed Kunz and Souillard \cite{KS80} to prove the following result, also in 1980.

\begin{theorem}\label{t.kunzsouillard}
Given an ergodic family $\{ H_\omega \}_{\omega \in \Omega}$, there exist sets $\Sigma_{{\rm ac}}, \Sigma_{{\rm ac}}, \Sigma_{{\rm pp}} \subseteq \R$ such that for $\mu$-almost every $\omega$, $\sigma_\bullet(H_\omega) = \Sigma_\bullet$, $\bullet \in \{ \mathrm{ac,sc,pp} \}$.
\end{theorem}

As a consequence, whenever an ergodic family of Schr\"odinger operators $\{ H_\omega \}_{\omega \in \Omega}$ is specified, it is one of the most basic goals to identify the sets $\Sigma, \Sigma_{{\rm ac}}, \Sigma_{{\rm sc}}, \Sigma_{{\rm pp}}$.

In general, one cannot claim more than mere full-measure statements. For instance, in Examples~\ref{x.hyptoraut}--\ref{x.symbolicshift} above, it is quite easy to see that there are many outliers, for which the spectrum and/or the spectral parts differ from the typical behavior.

On the other hand, sometimes it is possible to go beyond that. This usually works by approximation and hence requires additional structure from $\Omega$. Suppose in addition to our general assumptions that $\Omega$ is in fact a compact metric space, $T$ is a homeomorphism, and the sampling function $f : \Omega \to \R$ is continuous. We say that $(\Omega,T)$ is \emph{minimal} if the \emph{orbit} $O(\omega) = \{ T^n \omega : n \in \Z \}$ is dense in $\Omega$ for every $\omega \in \Omega$. If $(\Omega,T)$ is minimal, then for each pair $\omega, \omega' \in \Omega$, $\omega$ can be approximated by a sequence chosen from the orbit of $\omega'$. For the associated potentials, this means pointwise convergence, and for the associated operators, this means strong operator convergence. The net result is that $\sigma(H_\omega) \subseteq \sigma(H_{\omega'})$. Reversing roles, we obtain the following result.

\begin{prop}\label{p.constantspec}
Suppose that $\Omega$ is a compact metric space, $T$ is a homeomorphism, and the sampling function $f : \Omega \to \R$ is continuous. If $(\Omega,T)$ is minimal, then there exists a set $\Sigma \subseteq \R$ such that for every $\omega \in \Omega$, $\sigma(H_\omega) = \Sigma$.
\end{prop}

A similar approximation argument works for the absolutely continuous spectrum, as shown by Last and Simon \cite{LS99}. The details regarding the necessary semi-continuity statement are far more involved however.

\begin{theorem}\label{t.lastsimon}
Suppose that $\Omega$ is a compact metric space, $T$ is a homeomorphism, and the sampling function $f : \Omega \to \R$ is continuous. If $(\Omega,T)$ is minimal, then there exists a set $\Sigma_{{\rm ac}} \subseteq \R$ such that for every $\omega \in \Omega$, $\sigma_{{\rm ac}}(H_\omega) = \Sigma_{{\rm ac}}$.
\end{theorem}

On the other hand, there is no such result for the singular continuous spectrum and the point spectrum. The simplest counterexample is given by the super-critical almost Mathieu operator with Diophantine frequency; compare Subsection~\ref{ss.mitransition}.

Minimality may be replaced by unique ergodicity. We still assume that $\Omega$ is a compact metric space, $T$ is a homeomorphism, and the sampling function $f : \Omega \to \R$ is continuous. We say that $(\Omega,T)$ is \emph{uniquely ergodic} if there is exactly one ergodic Borel probability measure $\mu$. This is equivalent to there being exactly one $T$-invariant Borel probability measure $\mu$ (this measure must then necessarily be ergodic). Moreover, $(\Omega,T)$ is called \emph{strictly ergodic} if it is both minimal and uniquely ergodic. Finally, we sometimes just say that $T$ is minimal, uniquely ergodic, or strictly ergodic when we are referring to the respective property of the topological dynamical system $(\Omega, T)$.

The pair of results above holds in the uniquely ergodic situation, as shown by Kotani \cite{K97}.

\begin{theorem}\label{t.kotaniue}
Suppose that $\Omega$ is a compact metric space, $T$ is a homeomorphism, and the sampling function $f : \Omega \to \R$ is continuous. If $(\Omega,T)$ is uniquely ergodic with unique invariant measure $\mu$, then there exist sets $\Sigma, \Sigma_{{\rm ac}} \subseteq \R$ such that for every $\omega \in \mathrm{supp} \, \mu$, we have $\sigma(H_\omega) = \Sigma$ and $\sigma_{{\rm ac}}(H_\omega) = \Sigma_{{\rm ac}}$.
\end{theorem}

For our key examples, the following statements hold true regarding minimality and unique ergodicity. A torus translation is minimal if and only if it is uniquely ergodic, which in turns holds if and only if $1, \alpha_1 , \ldots , \alpha_d$ are linearly independent over the rational numbers. The skew-shift is minimal if and only if it is uniquely ergodic, which in turns holds if and only if $\alpha$ is irrational. Examples~\ref{x.hyptoraut}--\ref{x.symbolicshift} are neither minimal, nor uniquely ergodic (assuming they are non-degenerate; i.e., the shift spaces do not consist of only a single element).

\subsection{Lyapunov Exponents and the Integrated Density of States}

\subsubsection{The Cocycles Generating the Transfer Matrices}

Recall that for any potential $V : \Z \to \R$, we associate transfer matrices via \eqref{e.onestepmats}--\eqref{e.transmatrices}. In our present setting, the potential depends on the parameter $\omega \in \Omega$, and we will therefore denote the transfer matrices associated with $V_\omega$ by $T_\omega$ and $A_\omega$. Thus, the solutions to
\begin{equation}\label{e.ergeve}
u(n+1) + u(n-1) + V_\omega(n) u(n) = E u(n)
\end{equation}
obey
\begin{equation}\label{e.ergmeve}
\begin{pmatrix} u(n+1) \\ u(n) \end{pmatrix} = A_\omega(n;E) \begin{pmatrix}
u(1) \\ u(0) \end{pmatrix}.
\end{equation}

Since the potentials $\{ V_\omega \}$ are dynamically defined, it is not surprising that the transfer matrices are dynamically defined as well. Concretely, consider for $E \in \C$, the following skew-product:
$$
(T,A_E) : \Omega \times \C^2 \to \Omega \times \C^2, \quad (\omega,v) \mapsto (T \omega , A_E(T \omega) v),
$$
where
$$
A_E(\omega) = \begin{pmatrix} E - f(\omega) & -1 \\ 1 & 0 \end{pmatrix}.
$$
The iterates of $(T,A_E)$ may be written in the form $(T,A_E)^n = (T^n, A_E^n)$ with a suitable choice of matrix function $\omega \mapsto A_E^n(\omega)$. In fact, it is easy to check that
$$
A_E^n(\omega) = A_{\omega}(n;E).
$$
In other words, the iteration of the map $(T,A_E)$ generates the transfer matrices $A_{\omega}(n;E)$ in the second component.

The cocycle $A_E$ is said to be \emph{uniformly hyperbolic} if there are $C > 0$ and $\lambda > 1$ such that $\|A_E^n(\omega)\| \ge C \lambda^{|n|}$ for every $\omega \in \Omega$ and $n \in \Z$. We write
$$
\mathcal{UH} = \{ E \in \C : A_E \text{ is uniformly hyperbolic} \} \quad \text{and} \quad \mathcal{UH}_\R = \mathcal{UH} \cap \R.
$$
All non-real $E$ belong to $\mathcal{UH}$, that is,
\begin{equation}\label{e.nonrealisuh}
\C \setminus \R \subseteq \mathcal{UH}.
\end{equation}

\subsubsection{Lyapunov Exponents}

It is readily seen that the so-called cocycle condition holds, $A_E^{m+n}(\omega) = A_E^m(T^n \omega) A^n_E(\omega)$. Since norms are submultiplicative, this shows that $f_n(\omega,E) = \log \| A^n_E(\omega) \|$ satisfies the subadditivity condition $f_{n+m}(\omega,E) \le f_n(\omega,E) + f_m(T^n \omega ,E)$. Kingman's Subadditive Ergodic Theorem therefore implies the following:

\begin{prop}\label{p.leexistence}
For every $E \in \C$, there is a number $L(E) \in [0,\infty)$, called the \textit{Lyapunov exponent}, so that
\begin{align*}
L(E) & = \inf_{n \ge 1} \frac1n \, \E( \log \| A^n_E(\omega) \| ) \\
& = \lim_{n \to \infty} \frac1n \, \E( \log \| A^n_E(\omega) \| ) \\
& = \lim_{n \to \infty} \frac{1}{n} \log \| A^n_E(\omega) \| \; \text{ for } \mu-\text{almost every } \omega \in \Omega.
\end{align*}
\end{prop}

Clearly, $L(E) > 0$ for every $E \in \mathcal{UH}$. The converse is in general not true, and hence we denote the set of energies at which $A_E$ is non-uniformly hyperbolic by
\begin{equation}\label{e.nuhdef}
\mathcal{NUH} = \{ E \in \C : L(E) > 0 \text{ and } E \not\in \mathcal{UH} \}.
\end{equation}
Note that by \eqref{e.nonrealisuh}, we have $\mathcal{NUH} \subseteq \R$. Finally, we also set
\begin{equation}\label{e.zdef}
\mathcal{Z} = \{ E \in \C : L(E) = 0 \},
\end{equation}
so that
$$
\R = \mathcal{UH}_\R \sqcup \mathcal{NUH} \sqcup \mathcal{Z}.
$$

Here is how this partition relates to the spectrum.

\begin{theorem}\label{t.spectrumandenergypartition}
In general, we have
\begin{equation}\label{e.lespecincl}
\mathcal{Z} \subseteq \Sigma \subseteq \mathcal{NUH} \sqcup \mathcal{Z}.
\end{equation}
Moreover, if $\Omega$ is a compact metric space, $T$ is a homeomorphism, $f$ is continuous, and the $T$-orbit of $\omega \in \Omega$ is dense, then
\begin{equation}\label{e.johnson1}
\sigma(H_\omega) = \mathcal{NUH} \sqcup \mathcal{Z}.
\end{equation}
In particular, if $(\Omega,T)$ is minimal and $f$ is continuous, then
\begin{equation}\label{e.johnson2}
\sigma(H_\omega) = \Sigma = \mathcal{NUH} \sqcup \mathcal{Z}
\end{equation}
for every $\omega \in \Omega$.
\end{theorem}

This result is generally referred to as Johnson's theorem; compare \cite{J86}. Thus, in the minimal situation, we have
\begin{equation}\label{e.johnson3}
\R \setminus \Sigma = \mathcal{UH}_\R.
\end{equation}
This suggests a way of proving that the spectrum is ``small,'' namely by showing that ``many'' or ``most'' energies $E$ belong to $\mathcal{UH}$. Dating back to the early 1980's, the appearance of Cantor spectra, a phenomenon earlier thought to be exotic, has been a topic of intense study. Originally observed for almost periodic potentials, it turned out to be much more prevalent. Let us describe a result that shows that the appearance of Cantor spectra is generic in a suitable sense.

Recall that a compact subset $C$ of $\R$ is called a Cantor set if it contains no isolated points and no intervals. A spectrum $\Sigma$ arising in the ergodic setting does not contain isolated points due to Theorem~\ref{t.pastursetofev}. Thus, in order to establish that $\Sigma$ is a Cantor set, one needs to show that it contains no intervals. If the identity \eqref{e.johnson3} holds, this is equivalent to proving that ``uniform hyperbolity is dense,'' that is, for a dense set of energies $E \in \R$, the associated cocycle $A_E$ is uniformly hyperbolic. Such a result was shown under suitable assumptions by Avila, Bochi, and Damanik \cite{ABD09}:

\begin{theorem}\label{t.abdthm}
Suppose $\Omega$ is a compact metric space and $T: \Omega \to \Omega$ is a strictly ergodic homeomorphism that fibers over an almost periodic dynamical system. This means that there exists an infinite compact abelian group $G$, some $\alpha \in G$, and an onto continuous map $h : \Omega \to G$ such that $h(T(\omega)) = h(\omega) + \alpha$ for every $\omega \in \Omega$. Then, for every $E \in \R$, the set
$$
\mathcal{UH}_E = \{ f \in C(\Omega,\R) : A_E \text{ is uniformly hyperbolic} \}
$$
is open and dense. In particular, the set
$$
\mathcal{CS} = \{ f \in C(\Omega,\R) : \Sigma \text{ is a Cantor set} \}
$$
is residual.
\end{theorem}

This shows that Cantor spectrum is generic for base transformations that are much more general than almost periodic ones. It suffices that they contain an almost periodic factor. In particular this result applies to the skew-shift, for which the result is quite surprising. Skew-shift models were expected to not have Cantor spectra, but this turned out to be wrong at least $C^0$-generically.

In fact, the real obstruction to generic Cantor spectra can be formulated in terms of the Schwartzman asymptotic cycle \cite{S57}, with which one can describe the possible gap labels (the possible values the integrated density of states can take in gaps of $\Sigma$); compare, for example, \cite{B92, J86}. Whenever the possible gap labels are dense, Cantor spectrum will be generic as shown by Avila, Bochi, and Damanik in \cite{ABD12}.

\subsubsection{The Integrated Density of States}

Define the probability measure $\nu$ on $\R$ by
\begin{equation}\label{e.idsdef}
\int g(E) \, d\nu(E) = \E \left( \langle \delta_0 , g(H_\omega) \delta_0 \rangle \right)
\end{equation}
for bounded measurable $g$. The measure $\nu$ is called the \emph{density of states measure} associated with the family $\{ H_\omega \}_{\omega \in \Omega}$. Note that by definition $\nu$ is the $\mu$-average of the spectral measure corresponding to the pair $(H_\omega,\delta_0)$, but also one-half the $\mu$-average of the canonical spectral measure associated with $H_\omega$. The function $N$ defined by
$$
N(E) = \int \chi_{(-\infty,E]} (E') \, d\nu(E')
$$
is called the \textit{integrated density of states}.

The following result of Avron and Simon \cite{AS83} follows quickly from the definition of $\nu$:

\begin{theorem}\label{t.as83}
The almost sure spectrum is given by the points of increase of $N$, that is, $\Sigma = \mathrm{supp} \, \nu$.
\end{theorem}

Here, as usual, $\mathrm{supp} \, \nu$ denotes the topological support of the measure $\nu$.

Some of the statements in Theorem~\ref{t.pastursetofev} are reflected in the integrated density of states as shown by Delyon and Souillard \cite{DS84}:

\begin{theorem}
The integrated density of states is continuous.
\end{theorem}

Here is a different approach to the density of states measure. Denote the restriction of $H_\omega$ to $[1,n]$ with Dirichlet boundary conditions by $H_{\omega}^{(n)}$. For $\omega \in \Omega$ and $n \ge 1$, define probability measures $\nu_{\omega,n}$ by placing uniformly distributed point masses at the eigenvalues $E^{(n)}_\omega (1) < \cdots < E^{(n)}_\omega (n)$ of $H_\omega^{(n)}$, that is,
$$
\int g(E) \, d\nu_{\omega,n}(E) = \frac{1}{n} \, \sum_{j = 1}^n g(E^{(n)}_\omega (j)).
$$
Then, for $\mu$-almost every $\omega \in \Omega$, the measures $\nu_{\omega,n}$ converge weakly to $\nu$ as $n \to \infty$.

\subsubsection{The Thouless Formula}

Avron and Simon \cite{AS83} proved the following formula connecting the density of states measure and the Lyapunov exponent (see also Craig-Simon \cite{CS83} for an alternative proof).

\begin{theorem}\label{t.thouless}
For every $E \in \C$, we have
\begin{equation}\label{e.thouless}
L(E) = \int \log | E' - E | \, d\nu(E').
\end{equation}
\end{theorem}

This formula is called the Thouless formula and it says that the Lyapunov exponent is the negative of the logarithmic potential of the density of states measure. Using this interpretation, the following result of Simon \cite{S07} (which is essentially already in \cite{ST92}) is not too difficult to deduce.

\begin{theorem}\label{t.dos.eq}
If $L$ vanishes identically on $\Sigma$, then $\nu$ is the equilibrium measure of the compact set $\Sigma$.
\end{theorem}

The equilibrium measure is the unique probability measure supported on $\Sigma$ that minimizes the logarithmic energy
$$
\mathcal{E}(\rho) = - \iint \log | x - y | \, d\rho(x) \, d\rho(y)
$$
among such measures. By a standard result in logarithmic potential theory, the existence and the uniqueness of the minimizer follow as soon as at least one measure with finite logarithmic energy exists. In the case at hand, the Thouless formula implies that the density of states measure $\nu$ has finite logarithmic energy (since the spectrum $\Sigma$ is compact and the Lyapunov exponent is bounded on it).

The Thouless formula also implies the following general regularity result for the integrated density of states as shown by Craig and Simon \cite{CS83}.

\begin{theorem}\label{t.loghoelderids}
The integrated density of states is $\log$-H\"older continuous, that is, there is some uniform constant $C$ such that for real $E_1,E_2$ with $|E_1 - E_2| < 1/2$,
$$
|N(E_1) - N(E_2)| \le C \left( \log \left( |E_1 - E_2|^{-1} \right) \right)^{-1}.
$$
\end{theorem}

In this general setting this bound is optimal; compare Craig \cite{C83} and Gan-Kr\"uger \cite{GK11}. The regularity statement can often be improved for specific cases. We will describe some results of this kind in later sections.

\subsection{Kotani Theory}\label{ss.kotani}

Given a set $A \subseteq \R$, the \textit{essential closure} of $A$ is defined as follows:
$$
\overline{A}^\mathrm{ess} = \{ E \in \R : | (E-\varepsilon,E+\varepsilon) \cap A| > 0 \text{ for every } \varepsilon > 0 \}.
$$
Here, $| \cdot |$ denotes Lebesgue measure on $\R$. Note that $\overline{A}^\mathrm{ess} = \emptyset$ if and only if $|A| = 0$.

Recall that $\mathcal{Z}$ denotes the set $\{ E \in \R : L(E) = 0 \}$.

\begin{theorem}\label{t.ipkthm}
$\Sigma_\mathrm{ac} = \overline{ \mathcal{Z} }^\mathrm{ess}$.
\end{theorem}

The inclusion ``$\subseteq$'' was proved by Ishii \cite{I73} and Pastur \cite{P80}. The other inclusion was proved by Kotani \cite{K84} and is a much deeper result. In fact, the Ishii-Pastur half of the result is really an immediate consequence of the general theory of one-dimensional Schr\"odinger operators. See, for example, \cite{B97, DS83, LS99}. Moreover, the Ishii-Pastur half of the result can be strengthened considerably as shown by Simon in \cite{S07}. Not only are the spectral measures purely singular on $\mathcal{NUH}$ for $\mu$-almost every $\omega \in \Omega$, they must be purely zero-dimensional there! (In fact, an even stronger statement is true, for $\mu$-almost every $\omega \in \Omega$, the restriction of the spectral measures to $\mathcal{NUH}$ admits a support of capacity zero; see \cite{S07}.)

Denote the spectral measure associated with $H_\omega$ and $\delta_0$ by $\nu_\omega$. In particular, the density of states measure $\nu$ is the $\mu$-average of the measures $\nu_\omega$. Consider the absolutely continuous parts of these measures and their Radon-Nikodym derivatives. Kotani \cite{K97} has shown that they are related as follows.

\begin{theorem}\label{t.kot97}
For almost every $E \in \mathcal{Z}$,
\begin{equation}\label{acpartids}
\frac{d\nu^{(\mathrm{ac})}}{dE} (E) = \E \left( \frac{d\nu_\omega^{(\mathrm{ac})}}{dE} (E) \right).
\end{equation}
\end{theorem}

This result has a useful consequence \cite{K97}:

\begin{coro}\label{c.pureaccoro}
The spectrum of $H_\omega$ is purely absolutely continuous for $\mu$-almost every $\omega \in \Omega$ if and only if the density of states measure is purely absolutely continuous and the Lyapunov exponent vanishes almost everywhere with respect to it.
\end{coro}

Let us pass from a measurable setting to a topological setting. To fix a universal topology, we consider spaces of sequences, on which the topology will be given by pointwise convergence. Given an ergodic dynamical system $(\Omega,\mu,T)$ and a measurable bounded sampling function $f : \Omega \to \R$ defining potentials $V_\omega(n) = f(T^n \omega)$ as before, we associate the following dynamical system $(R^\Z,\tilde \mu,S)$: $R$ is a compact interval that contains the range of $f$, $\tilde \mu$ is the Borel measure on $R^\Z$ induced by $\mu$ via $\Phi(\omega) = V_\omega$ (i.e., $\tilde \mu(A) = \mu(\Phi^{-1}(A))$), and $S$ is the standard shift transformation on $R^\Z$. Clearly, the topological support $\mathrm{supp} \, \tilde \mu$ is closed and $S$-invariant.

For an $S$-ergodic Borel measure $\tilde \mu$ on $R^\Z$, let $\Sigma_\mathrm{ac}(\tilde \mu) \subseteq \R$ denote the almost sure absolutely continuous spectrum, that is, $\sigma_\mathrm{ac} (\Delta + V) = \Sigma_\mathrm{ac}(\nu)$ for $\tilde \mu$ almost every $V$. If $\tilde \mu$ comes from $(\Omega,\mu,T,f)$, then $\Sigma_\mathrm{ac}(\tilde \mu)$ coincides with the set $\Sigma_\mathrm{ac}$ introduced earlier. The support theorem \cite{K85} says that $\Sigma_\mathrm{ac}(\tilde \mu)$ is monotonically decreasing in the support of $\tilde \mu$.

\begin{theorem}\label{t.suppthm}
For every $V \in \mathrm{supp} \, \tilde \mu$, we have $\sigma_\mathrm{ac} (\Delta + V) \supseteq \Sigma_\mathrm{ac}(\tilde \mu)$, and hence
$$
\Sigma_\mathrm{ac}(\tilde \mu) = \bigcap_{V \in \mathrm{supp} \tilde \mu} \!\!\sigma_\mathrm{ac}(\Delta + V).
$$
In particular, $\mathrm{supp} \, \tilde \mu_1 \subseteq \mathrm{supp} \, \tilde \mu_2$ implies that $\Sigma_\mathrm{ac} (\tilde \mu_1) \supseteq \Sigma_\mathrm{ac}(\tilde \mu_2)$.
\end{theorem}

Here is a typical application of the support theorem:

\begin{coro}\label{c.ppgaps}
Let $\mathrm{Per}_{\tilde \mu}$ be the set of $V \in \mathrm{supp} \, \tilde \mu$ that are periodic, that is, $S^p V = V$ for some $p \in \Z_+$. Then,
$$
\Sigma_\mathrm{ac}(\tilde \mu) \subseteq \bigcap_{V \in \mathrm{Per}_{\tilde \mu}} \!\!\sigma(\Delta + V).
$$
\end{coro}

If there are sufficiently many gaps in the spectra of these periodic operators, one can show in this way that $\Sigma_\mathrm{ac}(\nu)$ is empty.

The following result shows that ergodic Schr\"odinger operators with non-empty absolutely continuous spectrum are deterministic.

\begin{theorem}\label{t.kottopthm}
Assume that $\mathrm{Leb} \, (\mathcal{Z}) > 0$. Then,\\
{\rm (a)} Each $V \in \mathrm{supp} \, \tilde \mu$ is determined completely {\rm (}among all elements of $\mathrm{supp} \, \tilde \mu${\rm )} by $V_- = V|_{\Z_-}$ {\rm (}resp., $V_+ = V|_{\Z_+}${\rm )}.\\
{\rm (b)} If we let
$$
(\mathrm{supp} \, \tilde \mu)_\pm = \{ V_\pm : V \in \mathrm{supp} \,
\tilde \mu \},
$$
then the mappings
$$
(\mathrm{supp} \, \tilde \mu)_\pm \ni V_\pm \mapsto V_\mp \in (\mathrm{supp} \, \tilde \mu)_\mp
$$
are continuous with respect to pointwise convergence.
\end{theorem}

Negating this statement, one obtains a criterion for purely singular spectrum. Namely, call $(\Omega,\mu,T,f)$ \textit{topologically deterministic} if there exist continuous mappings $E_\pm : (\mathrm{supp} \, \tilde \mu)_\pm \to (\mathrm{supp} \, \tilde \mu)_\mp$ that are formal inverses of one another and obey $V^\#_- \in \mathrm{supp} \, \tilde \mu$ for every $V_- \in (\mathrm{supp} \, \nu)_-$, where
$$
V^\#_-(n) = \begin{cases} V_-(n) & n \le 0, \\ E_-(V_-)(n) & n \ge 1. \end{cases}
$$
This also implies $V^\#_+ \in \mathrm{supp} \, \tilde \mu$ for every $V_+ \in (\mathrm{supp} \, \tilde \mu)_+$, where
$$
V^\#_+(n) = \begin{cases} V_+(n) & n \ge 1, \\ E_+(V_+)(n) & n \le 0. \end{cases}
$$
Otherwise, $(\Omega,\mu,T,f)$ is said to be \textit{topologically non-deterministic}.

\begin{coro}\label{topnd}
If $(\Omega,\mu,T,f)$ is \textit{topologically non-deterministic}, we have that $\mathrm{Leb} \, (\mathcal{Z}) = 0$, and therefore $\Sigma_\mathrm{ac} = \emptyset$.
\end{coro}

This readily applies to the random case, but also to certain models with weak correlations. Here is another, less obvious, application, which turns out to have far-reaching consequences.

\begin{theorem}\label{t.kotthmfv}
Suppose that $(\Omega, T, \mu)$ is ergodic, $f : \Omega \to \R$ takes finitely many values, and the resulting potentials $V_\omega$ are $\mu$-almost surely not periodic. Then, $\mathrm{Leb} \, (\mathcal{Z}) = 0$, and therefore $\Sigma_\mathrm{ac} = \emptyset$.
\end{theorem}

This was shown by Kotani in \cite{K89}. The proof is actually quite short, given Theorem~\ref{t.kottopthm} above. The almost everywhere positivity of the Lyapunov exponent for non-periodic ergodic potentials taking finitely many values is the basis for quite extensive work on subshift potentials, some of which will be discussed in Section~\ref{s.7}.

The fact that the potentials take finitely many values is actually not that crucial. It suffices that the sampling function has a discontinuity that can be exploited. This of course needs the topological situation to be present from the outset. Consider the case where $\Omega$ is a compact metric space, $T$ is a homeomorphism, and $\mu$ is an ergodic Borel probability measure. We say that $l\in\R$ is an \textit{essential limit} of $f$ at $\omega_0$ if there exists a sequence $\{\Omega_k\}$ of sets each of positive measure such that for any sequence $\{\omega_k\}$ with $\omega_k \in \Omega_k$, both $\omega_k\to\omega_0$ and $f(\omega_k) \to l$. If $f$ has more than one essential limit at $\omega_0$, we say that $f$ is \textit{essentially discontinuous} at this point. Damanik and Killip \cite{DK05} showed the following:

\begin{theorem}\label{t.damkilthm1}
Suppose $\Omega$ is a compact metric space, $T : \Omega \to \Omega$ a homeomorphism, and
$\mu$ an ergodic Borel probability measure. If there is an $\omega_0 \in \Omega$ such
that $f$ is essentially discontinuous at $\omega_0$ but continuous at all points $T^n
\omega_0$, $n < 0$, then $\mathrm{Leb} \, (\mathcal{Z}) = 0$, and hence $\Sigma_\mathrm{ac} = \emptyset$.
\end{theorem}

Let us now turn to the case of continuous sampling functions $f$. The proof of Theorem~\ref{t.damkilthm1} certainly breaks down and it is not clear where some sort of non-determinism should come from in the quasi-periodic case, for example. Of course, absence of absolutely continuous spectrum does not hold for a general continuous $f$. Thus, the following result by Avila and Damanik \cite{AD05} is somewhat surprising:

\begin{theorem}\label{t.avdamthm1}
Suppose $\Omega$ is a compact metric space, $T : \Omega \to \Omega$ a homeomorphism, and
$\mu$ a non-atomic ergodic Borel probability measure. Then, there is a residual set of
functions $f$ in $C(\Omega)$ such that $\Sigma_{{\rm ac}}(f) = \emptyset$.
\end{theorem}

Recall that a subset of $C(\Omega)$ is called residual if it contains a countable intersection of dense open sets. A residual set is locally uncountable.

One would expect some absolutely continuous spectrum for weak perturbations with sufficiently nice potentials; especially in the one-frequency quasi-periodic case. However, the proof of Theorem~\ref{t.avdamthm1} can easily be adapted to yield the following result, also contained in \cite{AD05}, which shows that continuity of the sampling function is not sufficient to ensure the existence of absolutely continuous spectrum for weakly coupled quasi-periodic potentials.

\begin{theorem}\label{avdamthm2}
Suppose $\Omega$ is a compact metric space, $T : \Omega \to \Omega$ a homeomorphism, and $\mu$ a non-atomic ergodic Borel probability measure. Then, there is a residual set of functions $f$ in $C(\Omega)$ such that $\Sigma_{{\rm ac}}(\lambda f) = \emptyset$ for almost every $\lambda > 0$.
\end{theorem}

As we have seen, there are many situations in which we have $\Sigma_\mathrm{ac} = \emptyset$. Of course, there are also cases where $\Sigma_\mathrm{ac} \not= \emptyset$, the most obvious being the periodic case. There are also some aperiodic cases as we will see later when we discuss specific classes of potentials. Nevertheless, for a long time all known examples with $\Sigma_\mathrm{ac} \not= \emptyset$ were almost periodic. This has led a number of people to conjecture that $\Sigma_\mathrm{ac} \not= \emptyset$ in fact implies almost periodicity of the potentials. Two of them were Kotani and Last, and hence this conjecture was sometimes called the Kotani-Last conjecture. It was recently explicitly stated in \cite{D07b, J07, S07}. However, the conjecture turned out to be wrong. It was disproved (in the form stated here and in \cite{D07b, J07, S07}) by Avila in \cite{A14f}.\footnote{There are related results for Jacobi matrices by Volberg and Yuditskii \cite{VY14} and for continuum Schr\"odinger operators by Avila \cite{A14f}, Damanik and Yuditskii \cite{DY14}, and You and Zhou \cite{YZ16}.} We will say more about this work in a later section.

\medskip

Let us also mention the surveys \cite{D07b, K97} of Kotani theory and its applications, where the interested reader can find further related material.

Moreover, it is a perhaps surprising, but certainly amazing, fact that much of Kotani theory has a deterministic counterpart, without any need of a dynamical definition of the potentials and an underlying ergodic measure, as shown by Remling in \cite{R11}.

\section{Random Potentials}\label{s.4}

In this section we discuss the case of random potentials. Random potentials arise in the setting of Example~\ref{x.sequenceshift}, where $\Omega = I^\Z$ with a compact interval $I \subset \R$, $T : \Omega \to \Omega$ is given by the shift transformation $(T \omega)_n = \omega_{n+1}$, $\mu = \rho^\Z$, where $\rho$ is a Borel probability measure on $I$, and the sampling function is given by the evaluation at the origin, $f(\omega) = \omega_0$.\footnote{This paper focused on the one-dimensional case. We wish to point out, however, that random potentials have been studied in great depth in higher dimensions as well. One understands the behavior near the edges of the spectrum very well. The methods employed in the analysis of random operators in dimensions greater than one are quite different; the two most prominent ones are based on a multi-scale analysis or the fractional moment method.} That is, the elements $\omega$ of $\Omega$ themselves serve as the potentials of the operators \eqref{e.ergoper}. Of course, we will assume that $\rho$ is non-degenerate in the sense that $\mathrm{supp} \, \rho$ contains more than one element. In this case, the resulting family $\{ H_\omega \}_{\omega \in \Omega}$ of operators is referred to as the \emph{Anderson model}. The special case where $\mathrm{supp} \, \rho$ contains precisely two elements is called the \emph{Bernoulli-Anderson model}. This is the case with the least amount of randomness, and as a consequence the proofs of the expected results for the Anderson model are the most difficult in the Bernoulli case. In fact, for nicer single-site distributions $\rho$, the proofs can be simpler by orders of magnitude.

\subsection{The Spectrum}

The spectrum of the Anderson model has a very simple description (where the sum below denotes the sum set $A + B = \{ a + b : a \in A, \; b \in B \}$):

\begin{theorem}\label{t.randomspectrum}
For the Anderson model, we have
\begin{equation}\label{e.randomspectrum}
\Sigma = [-2,2] + \mathrm{supp} \, \rho.
\end{equation}
\end{theorem}

The proof of this result is not difficult; let us sketch it. First of all, $\mu$-almost all elements $\omega \in \Omega$ will be such that the range of $V_\omega$ is dense in $\mathrm{supp} \, \rho$. In other words, for these $\omega$'s, we have $\sigma(V_\omega) = \mathrm{supp} \, \rho$. Since the norm of the Laplacian is bounded by (in fact, it is equal to) $2$, it follows for these $\omega$'s that $\sigma(H_\omega) \subseteq [-2,2] + \mathrm{supp} \, \rho$. This establishes the inclusion ``$\subseteq$'' in \eqref{e.randomspectrum}. Conversely, $\mu$-almost surely there are for each $E \in \mathrm{supp} \, \rho$ long stretches where $V_\omega$ only takes values very close to $E$. By considering suitable trial functions for the Laplacian, one can derive from this that $[-2,2] + E$ must be contained in the almost sure spectrum. Since this is true for every $E \in \mathrm{supp} \, \rho$, the inclusion ``$\supseteq$'' in \eqref{e.randomspectrum} follows.

Theorem~\ref{t.randomspectrum} shows that the spectrum of an Anderson model cannot be arbitrary. In fact, it will always be given by a finite union of compact intervals since the spectrum is bounded and each of its connected components has length at least $4$. Conversely, any compact set with the property that each of its connected components has length at least $4$ arises as the almost sure spectrum of a suitably chosen Anderson model.

\subsection{Various Notions of Anderson Localization}

So, what are the expected results? One says that the Anderson model exhibits \emph{Anderson localization}. There are usually two different statements that are referred to, a spectral statement and a (quantum) dynamical statement. \emph{Spectral Anderson localization} is the assertion that for $\mu$-almost every $\omega \in \Omega$, the operator $H_\omega$ has pure point spectrum with exponentially decaying eigenfunctions. More precisely,  for $\mu$-almost every $\omega \in \Omega$, there are $E_k(\omega) \in \R$ and $u_k(\omega) \in \ell^2(\Z)$ such that $H_\omega u_k(\omega) = E_k(\omega) u_k(\omega)$ for every $k$, $\{ u_k(\omega) \}_k$ form a basis of $\ell^2(\Z)$, and
\begin{equation}\label{e.expdecsol}
|u_k(n;\omega)| \le C_{k;\omega} e^{-\gamma_{k;\omega} |n|}
\end{equation}
with suitable constants $C_{k;\omega}, \gamma_{k;\omega} > 0$. \emph{Dynamical Anderson localization} is a less well-defined notion, but it typically means at least that for $\mu$-almost every $\omega \in \Omega$, we have
\begin{equation}\label{e.dynamicallocalization}
\sup_t \sum_{n \in \Z} |n|^p | \langle \delta_n , e^{-itH_\omega} \delta_0 \rangle |^2 < \infty
\end{equation}
for every $p > 0$. There are stronger statements that can be proved in some cases, such as replacing the $\mu$-almost everywhere statement by an expectation $\E(\cdot)$, or by claiming explicit (semi-)uniform exponential decay of $| \langle \delta_n , e^{-itH_\omega} \delta_0 \rangle |$. But in any event, dynamical Anderson localization refers to the absence of transport in a random medium.

The two notions of Anderson localization are related, though not equivalent. Indeed, dynamical localization in a suitable formulation implies spectral localization, while the converse does not hold in general. For an example with a good amount of randomness, for which the implication ``spectral localization $\Rightarrow$ dynamical localization'' fails fairly spectacularly, one can consider the so-called random dimer model. Starting with the Bernoulli-Anderson model, with $\mathrm{supp} \, \rho = \{ 0, \lambda \}$ say, the random dimer model results from doubling up all the sites. That is, the operator $H_\omega$ has the potential $V_\omega$ with
$$
V_\omega(2n-1) = V_\omega(2n) = \omega_n
$$
for every $n \in \Z$. This model can be realized in our framework by using Example~\ref{x.symbolicshift} with a suitable Markov measure; compare the discussion in \cite{AD14}.

When considering the transfer matrices $A_\omega(n;E)$ associated with the random dimer model, it is natural to group the factors in pairs, which are
$$
\begin{pmatrix} E & -1 \\ 1 & 0 \end{pmatrix}^2 \quad \text{ and } \quad \begin{pmatrix} E - \lambda & -1 \\ 1 & 0 \end{pmatrix}^2.
$$
In particular, for the energy $E = 0$, we have the basic building blocks
$$
\begin{pmatrix} 0 & -1 \\ 1 & 0 \end{pmatrix}^2 = \begin{pmatrix} -1 & 0 \\ 0 & -1 \end{pmatrix} \quad \text{ and } \quad \begin{pmatrix} - \lambda & -1 \\ 1 & 0 \end{pmatrix}^2.
$$
In particular, up to a sign, the matrix $A_\omega(2n;0)$ will be given by a power of the matrix
$$
\begin{pmatrix} - \lambda & -1 \\ 1 & 0 \end{pmatrix},
$$
which means that $\|A_\omega(2n;0)\|$ remains bounded in $n$, provided that $|\lambda| < 2$. This of course implies that $\|A_\omega(n;0)\|$ is bounded as well. The methods of Subsection~\ref{ss.dtbounds} therefore imply quasi-ballistic transport for every $\omega \in \Omega$! In particular, dynamical localization fails in a rather extreme way (see \cite{JS07, JSS03} for a precise description of the quantum dynamics of this model). On the other hand, the method of proof outlined below applies to the random dimer model and implies that spectral localization holds for this model; compare \cite{dBG00}. Moreover, it is a general result of Simon that pure point spectrum implies the absence of genuine ballistic transport \cite{S90}.

There are ways to supplement the requirement that there be an orthonormal basis of exponentially decaying eigenvectors in such a way that dynamical localization in a suitable formulation is indeed a consequence. Such connections were first established in the paper \cite{DJLS96} by del Rio, Jitomirskaya, Last, and Simon. For example, if for some $\omega \in \Omega$, one has semi-uniformly localized eigenfunctions (SULE) in the sense that there are $\alpha > 0$ and $\{ n_m \} \subseteq \Z$ such that for each $\delta > 0$, there is $C_\delta$ so that the eigenvectors $u_m$ obey
$$
|u_m(n)| \le C_\delta e^{\delta |n_m| - \alpha | n - n_m |},
$$
then it follows that semi-uniform dynamical localization (SUDL) holds, that is,
$$
\sup_{t \in \R} \left| \left\langle \delta_n, e^{-itH_\omega} \delta_m \right\rangle \right| \le C_\delta e^{\delta |m| - \alpha | n - m |}.
$$

\subsection{Positivity of the Lyapunov Exponent}\label{ss.furstenberg}

The first step in a proof of Anderson localization is the proof of positivity for the Lyapunov exponent. Indeed, as discussed earlier, the positivity of the Lyapunov exponent is a necessary condition for exponential decay of solutions of \eqref{e.ergeve}, which in turn is necessary for there to exist exponentially decaying eigenvectors of $H_\omega$.

A general theorem of F\"urstenberg about products of random matrices is tailor-made for this particular goal. In fact, the general result applies easily and in full generality to yield the positivity of the Lyapunov exponent at all energies for every realization of the Anderson model. Let us first state F\"urstenberg's theorem and then show how it may be applied to the Anderson model.

Let $\tilde \rho$ be a probability measure on $\mathrm{SL}(2,\R)$ which satisfies
\begin{equation}\label{logintegrable}
\int \log \| M \| \, d\tilde \rho(M) < \infty.
\end{equation}
Let us consider i.i.d.\ matrices $T_1,T_2,\ldots$, each distributed according to $\tilde \rho$. Write $M_n = T_n \cdots T_1$. We are interested in the Lyapunov exponent $L \ge 0$, given by
$$
L = \lim_{n \to \infty} \frac1n \log \| M_n \| , \quad \tilde \rho^{\Z_+}-\text{a.s.}
$$
We are interested in conditions that ensure $L > 0$. To motivate the result below, let us give some examples with $L = 0$:
\begin{itemize}

\item If $\tilde \rho$ is supported in $\mathrm{SO}(2,\R)$, then $L = 0$.

\item If
$$
\tilde \rho \left\{ \left( \begin{array}{cc} 2 & 0 \\ 0 & 1/2 \end{array} \right) \right\} = \frac12 \quad \text{ and } \quad \rho \left\{ \left( \begin{array}{cc} 1/2 & 0 \\ 0 & 2 \end{array} \right) \right\} = \frac12,
$$
then $L = 0$: We have that
$$
M_n = \left( \begin{array}{cc} m_n & 0 \\ 0 & m_n^{-1} \end{array} \right),
$$
where $\log m_n = a_1 + \cdots + a_n$ and $\{a_j\}$ are i.i.d.\ random variables taking values $\pm \log 2$, each with probability $1/2$. Thus, $\log \| M_n \| = | a_1 + \cdots + a_n |$ and the strong law of large numbers gives $\frac1n \log \| M_n \| \to 0$ almost surely.

\item If $p \in (0,1)$ and
$$
\tilde \rho \left\{ \left( \begin{array}{cc} 2 & 0 \\ 0 & 1/2 \end{array} \right) \right\} = p \quad \text{ and } \quad \tilde \rho \left\{ \left( \begin{array}{rc} 0 & 1 \\ -1 & 0 \end{array} \right) \right\} = 1-p,
$$
then $L = 0$.

\end{itemize}

F\"urstenberg's Theorem shows that this list is essentially exhaustive in the sense that the two mechanisms above, no growth of norms or a finite (cardinality $=2$) invariant set of directions, are the only ones that can preclude a positive Lyapunov exponent.

Call two non-zero vectors $v_1,v_2$ in $\R^2$ equivalent if $v_2 = \lambda v_1$ for some $\lambda \in \R$. The set of equivalence classes is denoted by $\R\Y^1$. Since every $M \in \mathrm{SL}(2,\R)$ is invertible, it induces a mapping from $\R\Y^1$ to $\R\Y^1$ in the obvious way.

\begin{theorem}\label{t.furst}
Let $\tilde \rho$ be a probability measure on $\mathrm{SL}(2,\R)$ which satisfies \eqref{logintegrable}. Denote by $G_{\tilde \rho}$ the smallest closed subgroup of $\mathrm{SL}(2,\R)$ which contains $\mathrm{supp} \, \tilde \rho$.

Assume
\begin{itemize}
\item[(i)] $G_{\tilde \rho}$ is not compact.

\item[(ii)] There is no set $L \subseteq \R\Y^1$ of cardinality $1$ or $2$ such that $M(L) = L$ for all $M \in G_{\tilde \rho}$.
\end{itemize}
Then, $L > 0$.
\end{theorem}

This is a special case of a much more general result proved by F\"urstenberg in \cite{F63}. Note that the assumptions are monotonic in the support of the measure in the sense that if Theorem~\ref{t.furst} applies to $\tilde \rho$, then it applies to measures whose support contains the support of $\tilde \rho$.

\bigskip

Let us now apply F\"urstenberg's Theorem to the Anderson model. Recall that $\mathrm{supp} \, \rho$ has cardinality $\ge 2$, and $f : \Omega \to \R$ is given by $f(\omega) = \omega_0$.

For every $E \in \R$ fixed, the measure $\rho$ on the interval $J$ induces the measure $\tilde \rho$ on $\mathrm{SL}(2,\R)$ by push-forward via
$$
v \mapsto \left( \begin{array}{cr} E - v & -1 \\ 1 & 0 \end{array} \right).
$$
The definitions are such that the Lyapunov exponent $L$ associated with this $\tilde \rho$ is equal to $L(E)$ defined earlier.

Let us check that F\"urstenberg's Theorem applies. Since $\mathrm{supp} \, \rho$ has cardinality at least two, $\mathrm{supp} \, \tilde \rho$ has cardinality at least two, and hence $G_{\tilde \rho}$ contains at least two distinct elements of the form
$$
M_x = \left( \begin{array}{cr} x & -1 \\ 1 & 0 \end{array} \right),
$$
for example, $M_a$ and $M_b$ with $a \not= b$. Note that
$$
M^{(1)} = M_a M_b^{-1} = \left( \begin{array}{cc} 1 & a-b \\ 0 & 1 \end{array} \right) \in G_{\tilde \rho}.
$$
Taking powers of the matrix $M^{(1)}$, we see that $G_{\tilde \rho}$ is not compact.

Consider the equivalence class of $e_1 = (1,0)^T$ in $\R\Y^1$. Then $M^{(1)} e_1 = e_1$ and for every $v \in \R\Y^1$, $(M^{(1)})^n v$ converges to $e_1$. Thus, if there is a finite invariant set of directions $L$, it must be equal to $\{e_1\}$. However,
$$
M^{(2)} = M_a^{-1} M_b = \left( \begin{array}{cc} 1 & 0 \\ a-b & 1 \end{array} \right) \in G_{\tilde \rho}
$$
and $M^{(2)} e_1 \not= e_1$; contradiction. Thus, the conditions (i) and (ii) of Theorem~\ref{t.furst} hold and, consequently, $L = L(E) > 0$.

Since $E \in \R$ was arbitrary and all we needed was $\# \mathrm{supp} \, \rho \ge 2$, we obtain the following consequence that holds in full generality.

\begin{theorem}\label{t.amposle}
For the Anderson model, we have $L(E) > 0$ for every $E \in \R$.
\end{theorem}

\subsection{Spectral Localization via Spectral Averaging}\label{ss.specaver}

Theorem~\ref{t.amposle} suggests that we are already tantalizingly close to being able to deduce one of our primary goals, namely spectral Anderson localization. Indeed, combined with Proposition~\ref{p.leexistence}, Theorem~\ref{t.amposle} implies that for every $E \in \R$,
$$
\lim_{|n| \to \infty} \frac{1}{|n|} \log \| A^n_E(\omega) \| = L(E) > 0 \; \text{ for } \mu-\text{almost every } \omega \in \Omega.
$$
Thus, Theorem~\ref{t.ruelle} is applicable and shows that at both $\pm \infty$, the solutions of \begin{equation}\label{e.ergeve.e}
u(n+1) + u(n-1) + V_\omega(n) u(n) = E u(n)
\end{equation}
either decay or increase exponentially. By Theorem~\ref{t.genef}, we can focus our attention on those energies $E$ for which polynomially bounded solutions exist. In the presence of our exponential dichotomy, this means that we can exclude exponential growth, and hence the generalized eigenfunctions are in fact exponentially decaying. This implies both that spectrally every energy is an eigenvalue, and that the corresponding eigenvectors decay exponentially!

So why is this not already a complete proof? The cheat here lies in a change of quantifiers. We have passed from a ``for every energy $E$ and almost every $\omega$'' statement to a statement of the form ``for almost every $\omega$ and every energy $E$.'' A more honest application of Fubini only allows us to conclude a ``for almost every $\omega$ and Lebesgue almost every energy $E$'' statement, and the exclusion of a set of zero Lebesgue measure still makes it possible that we lost some singular continuous spectrum in the process.

All is not lost, however. If we could force spectral measures away from the set of zero Lebesgue measure that needs to be excluded, the argument above still works and then allows us to conclude as desired. The mechanism that can be employed to force spectral measures away from sets of zero Lebesgue measure is called \emph{spectral averaging}. In the context of the Anderson model it can be applied whenever the single-site distribution has a non-trivial absolutely continuous component. In this subsection we explain how this works and how one may deduce spectral Anderson localization for such ``nice'' single-site distributions in a rather elegant way. For a paper pioneering spectral averaging methods in proofs of spectral localization, see \cite{SW86} by Simon and Wolff. See also \cite{S95} for an introductory paper on rank-one perturbations and applications to spectral averaging.

The key input is the spectral averaging formula from the theory of rank-one perturbations, which we now recall. Suppose that $A$ is a bounded self-adjoint operator on $\ell^2(\Z)$ and $\phi \in \ell^2(\Z) \setminus \{ 0 \}$. For $\lambda \in \R$, we consider the operator
$$
A_\lambda = A + \lambda \langle \phi, \cdot \rangle \phi,
$$
which is a self-adjoint rank one perturbation of $A$. Denote the spectral measure associated with $A_\lambda$ and $\phi$ by $\mu_\lambda$. Then, we have
\begin{equation}\label{f.specavertp}
\int \left[ d\mu_\lambda(E) \right] \, d\lambda = dE
\end{equation}
in the sense that if $g \in L^1(\R,dE)$, then $g \in L^1(\R,d\mu_\lambda)$ for Lebesgue almost every $\lambda$, $\int g(E) \, d\mu_\lambda(E) \in L^1(\R,d\lambda)$, and
$$
\int \left( \int g(E) \, d\mu_\lambda(E) \right) \, d\lambda = \int g(E) \, dE.
$$

Now let us return to our discussion of the Anderson model. Assume that $\rho_\mathrm{ac} \not= 0$. By the argument described above, we have
\begin{equation}\label{f.noweight}
\mathrm{Leb} (\R \setminus \mathcal{E}_\omega) = 0
\end{equation}
for $\mu$-almost every $\omega \in \Omega$, where
$$
\mathcal{E}_\omega = \{ E \in \R : L(E) > 0 , \, \exists \text{ solutions } u_\pm \text{ with } |u_\pm(n)| \sim e^{-L(E) |n|} \text{ as } n \to \pm \infty \}.
$$
Note that the sets $\mathcal{E}_\omega$ are invariant with respect to a modification of $V_\omega$ on a finite set! We will perform such a modification, within the family $\{V_\omega\}$, on the set $\{0,1\}$ because the pair $\{ \delta_0, \delta_1 \}$ is cyclic for each operator $H_\omega$.

Denote the set of $\omega$'s for which \eqref{f.noweight} holds by $\Omega_0$. By invariance and ergodicity, it follows that
\begin{equation}\label{f.omega0fm}
\mu(\Omega_0) = 1.
\end{equation}
For $\omega \in \Omega_0$, consider the operators
$$
H_{\omega,\lambda_0,\lambda_1} = H_\omega + \lambda_0 \langle \delta_0 , \cdot \rangle \delta_0 + \lambda_1 \langle \delta_1 , \cdot \rangle \delta_1,
$$
where $\lambda_0, \lambda_1 \in \R$. For every fixed $\lambda_0$, it follows from \eqref{f.specavertp} and \eqref{f.noweight} that the spectral measure of the pair $(H_{\omega,\lambda_0,\lambda_1} , \delta_1)$ gives zero weight to the set $\R \setminus \mathcal{E}_\omega$ for Lebesgue almost
every $\lambda_1 \in \R$. Similarly, for every fixed $\lambda_1$, the spectral measure of the pair $(H_{\omega,\lambda_0,\lambda_1}
, \delta_0)$ gives zero weight to the set $\R \setminus \mathcal{E}_\omega$ for Lebesgue almost every $\lambda_0 \in \R$. As a consequence, we find that for Lebesgue almost every $(\lambda_0,\lambda_1) \in \R^2$, the universal spectral measure of $H_{\omega,\lambda_0,\lambda_1}$ (the sum of the spectral measures of $\delta_0$ and $\delta_1$) gives zero weight to the set $\R \setminus \mathcal{E}_\omega$. Write $G_\omega$ for this set of ``good'' pairs $(\lambda_0,\lambda_1)$, so that
\begin{equation}\label{f.gomegafm}
\mathrm{Leb} (\R^2 \setminus G_\omega) = 0.
\end{equation}
Let
$$
\Omega_1 = \{ \omega + \lambda_0 \delta_0 + \lambda_1 \delta_1 : \omega \in \Omega_0 , \; (\lambda_0,\lambda_1) \in G_\omega \}.
$$
Since $\rho_\mathrm{ac} \not= 0$, it follows that from \eqref{f.omega0fm} and \eqref{f.gomegafm} that
$$
\mu(\Omega_1) > 0.
$$
Thus, by assumption on $\rho$, with positive $\rho \times \rho$ probability, it follows from \eqref{f.noweight} that the whole-line spectral measure (corresponding to the sum of the $\delta_0$ and $\delta_1$ spectral measures) assigns no weight to $\R \setminus \mathcal{E}_\omega$ and hence, with positive $\mu$ probability, the operator $H_\omega$ is spectrally localized. Since localization is a shift-invariant event, the operator $H_\omega$ must in fact be spectrally localized for $\mu$-almost every $\omega$.

This establishes the following result:

\begin{theorem}\label{t.kotlocwl}
For the Anderson model with a single-site distribution $\rho$ that has a non-trivial absolutely continuous component, we have that the family $\{ H_\omega \}_{\omega \in \Omega}$ is spectrally localized.
\end{theorem}

\subsection{Spectral and Dynamical Localization via Multi-Scale Analysis}\label{ss.speclocalization}

While the proof of Theorem~\ref{t.kotlocwl} is elegant and relatively short, with more work one can establish spectral localization for the Anderson model in full generality.

\begin{theorem}\label{t.speclocandmod}
For the Anderson model with a single-site distribution $\rho$ whose support is bounded and contains at least two elements, we have that the family $\{ H_\omega \}_{\omega \in \Omega}$ is spectrally localized.
\end{theorem}

This is a special case of a result of Carmona, Klein, and Martinelli \cite{CKM87}. The assumption that $\mathrm{supp} \, \rho$ contains at least two elements is clearly necessary (as otherwise $\mu$-almost all potentials $V_\omega$ are constant and hence $H_\omega$ almost surely has purely absolutely continuous spectrum. The assumption that $\mathrm{supp} \, \rho$ is bounded is not necessary. In fact, spectral localization is proved in \cite{CKM87} under the weaker assumption that $\rho$ has some finite moment. Since we focus in this paper on the case of bounded ergodic potentials, we impose the corresponding condition on $\rho$ in Theorem~\ref{t.speclocandmod}.

The proof of Theorem~\ref{t.speclocandmod} is based on multi-scale analysis. This is a method for proving localization that was originally introduced by Fr\"ohlich and Spencer \cite{FS83} and then developed further in many papers (e.g., in \cite{GK01, vDK89}, among many others).

The purpose of a multi-scale analysis is to inductively prove decay estimates for the resolvent of finite-volume restrictions of the operator that hold with large probability. To make this inductive procedure work, one needs two ingredients, an initial length-scale estimate and a Wegner estimate. The former is used to establish the base case, and the latter is used in the induction step.

The initial length-scale estimate can be established, in the one-dimensional case we consider here, as a consequence of the positivity of the Lyapunov exponent, that is, the result provided by Theorem~\ref{t.amposle}, which holds in complete generality. As we saw above, Theorem~\ref{t.amposle} is proved for the Bernoulli case ($\# \mathrm{supp} \, \rho = 2$) and then derived for the general case using the monotonicity of the argument in the support of $\rho$. The Wegner estimate, on the other hand, is quite easy to prove for nice single-site distributions $\rho$, but it is quite difficult to establish in more singular cases (of which the Bernoulli case is the most singular one). Thus, this is precisely the point \cite{CKM87} had to address, and the authors accomplish this by deriving a Wegner-type estimate from H\"older continuity of the integrated density of states, which in turn follows from H\"older continuity of the Lyapunov exponent via the Thouless formula, with the latter property being a consequence of the F\"urstenberg approach to positive Lyapunov exponents; see Le Page \cite{L84}. That is, a finite-volume statement (the Wegner-type estimate) is derived from an infinite-volume statement (the H\"older continuity of the integrated density of states), which in turn is used to prove an infinite-volume statement (spectral localization; as a consequence of a successful multi-scale induction).

\bigskip

It is a realization of Germinet and de Bi\`evre \cite{GD98} that if one can successfully carry out a multi-scale analysis, then one not only gets spectral localization as a consequence, but also dynamical localization. This connection was developed further in several papers; see, for example, \cite{DS01, GK01}. For the model at hand we may therefore state the following.

\begin{theorem}\label{t.speclocandmod}
For the Anderson model with a single-site distribution $\rho$ whose support is bounded and contains at least two elements, we have that the family $\{ H_\omega \}_{\omega \in \Omega}$ is dynamically localized.
\end{theorem}

This result certainly includes the statement that \eqref{e.dynamicallocalization} holds $\mu$-almost surely, but one can, for example, also make statements about $\mu$-expectations and replace $p$-th moments by larger (sub-exponential) functions.

\section{Almost Periodic Potentials}\label{s.5}

In this section we discuss the class of almost periodic potentials. A bounded $V : \Z \to \R$ is called \emph{almost periodic} if the set of its translates has compact closure in $\ell^\infty(\Z)$. That is, on $\ell^\infty(\Z)$, consider the shift transformation $S : \ell^\infty(\Z) \to \ell^\infty(\Z)$ given by $(S(W))(n) = W(n+1)$. The shift orbit of $V$ is $O(V) = \{ S^m(V) : m \in \Z \}$, and $V$ is almost periodic if and only if the closure of $O(V)$ in $\ell^\infty(\Z)$ is compact. It turns out that Schr\"odinger operators with almost periodic potentials may be studied within the framework of ergodic Schr\"odinger operators. To this end, we need appropriate choices of $\Omega,T,\mu$, and $f$. This will be discussed in Subsection~\ref{ss.aphull}. Once this has been realized, the general results from the theory of ergodic Schr\"odinger operators become applicable to the almost periodic case. There are three important subclasses of almost periodic potentials, and they will be discussed in Subsections~\ref{ss.periodic}--\ref{ss.quasiper}. These are periodic potentials, limit-periodic potentials, and quasi-periodic potentials. All three classes describe physically relevant models, and each of them has a rich mathematical theory. In this section we will survey some of the most important results for them. It will be impossible to be comprehensive, and we will provide the reader with pointers for further reading.

\subsection{The Hull}\label{ss.aphull}

Suppose $V$ is almost periodic. Let us denote the closure of $O(V)$ in $\ell^\infty(\Z)$ by $\Omega(V)$. The set $\Omega(V)$ is called the \emph{hull} of $V$. We want to equip $\Omega(V)$ with an abelian group structure. Since the dense subset $O(V)$ of $\Omega(V)$ already carries a natural group structure, we wish to extend it to the closure by continuity. Of course, if $V$ is a periodic point of $S$, both orbit and hull are easily seen to be isomorphic to $\Z_p = \Z/p\Z$, where $p$ is the (minimal) period. For simplicity, let us exclude this degenerate case from the discussion of the extension of the group structure and assume that $V$ is not a periodic point of $S$. In this case, the group structure on $O(V)$ is that of $\Z$.

Ordinarily, we would denote the abelian group structure by $+$, but to avoid confusion with the operation of adding functions, we will denote it by $\ast$. Thus, $O(V)$ is a group under the operation
$$
S^{k_1} V \ast S^{k_2} V = S^{k_1 + k_2} V,
$$
and $(O(V),\ast) \simeq (\Z,+)$ by our non-periodicity assumption. For $W_1 = \lim_{j \to \infty} S^{k_j} V$ and $W_2 = \lim_{j \to \infty} S^{\ell_j} V$, we wish to define $W_1 \ast W_2$, and the most natural choice is to set
\begin{equation}\label{e.hullgroupop}
W_1 \ast W_2 = \left( \lim_{j \to \infty} S^{k_j} V \right) \ast \left( \lim_{j \to \infty} S^{\ell_j} V \right) = \lim_{j \to \infty} S^{k_j + \ell_j} V.
\end{equation}
It is not hard to see that this indeed converges and is in fact well-defined. As a consequence, the operation \eqref{e.hullgroupop} equips the compact space $\Omega(V)$ with an abelian group structure.

In particular, we can choose normalized Haar measure on $\Omega(V)$ as our probability measure, and it will be invariant with respect to the transformation in question, which is the restriction of the shift transformation $S$ to $\Omega(V)$ (since by \eqref{e.hullgroupop}, the shift on $\Omega(V)$ is just given by the action of $S(V)$, i.e.\ $S(W) = W \ast S(V)$ for every $W \in \Omega (V)$). Finally the sampling function to be used is the evaluation at the origin, $f(W) = W(0)$, $W \in \Omega(V)$. In this way, every $W \in \Omega(V)$ has the desired representation \eqref{e.ergpotential}, which in this case becomes $W(n) = f(S^n(W))$.

To summarize, every almost periodic potential may be realized as one element in a canonically chosen ergodic family of potentials $\{ V_\omega \}_{\omega \in \Omega}$. (The periodic case excluded in the discussion above is trivial.) In this family, $\Omega$ is a compact abelian group, $T$ is a minimal translation (i.e., by the action of a fixed group element, and so that all orbits are dense), and $\mu$ is normalized Haar measure.

Conversely, one can start with the latter scenario and generate almost periodic potentials in this way.

\begin{prop} \label{p.abeliangroup.char.alm.per}
A potential $V \in \ell^{\infty}(\Z)$ is almost periodic if and only if it can be represented via
\begin{equation} \label{f.abeliangroup.char.alm.per}
V(n) = V_{\omega}(n) = f(T^n \omega),
\end{equation}
where $ \Omega $ is a compact abelian group, $\omega \in \Omega$, $f: \Omega \to \R$ is continuous, and $ T:\Omega \to \Omega $ is a minimal translation, say $T = \cdot + \alpha$.
\end{prop}

This point of view is sometimes useful, especially in the discussion of the specific subclasses below.

\subsection{Periodic Potentials}\label{ss.periodic}

A potential is \emph{periodic} if and only if it is a fixed point of the shift transformation $S$. That is, there is some $p \in \Z_+$ such that $S^p(V) = V$. In other words, $V(n + p) = V(n)$ for every $n \in \Z$. We assume that $p \in \Z_+$ is minimal with this property; in this case it is called the \emph{minimal period} of $V$, $O(V) = \Omega(V)$ is isomorphic to $\Z_p$, and $\mu$ assigns the weight $1/p$ to each element of the hull.

Since every $V_\omega$ is a translate of any other $V_{\omega'}$ in this case, the associated operators are unitarily equivalent, and the constancy of the spectrum and the spectral parts is immediate. It is therefore sufficient to describe the spectral properties of a single operator, and we discuss those of the Schr\"odinger operator with the initial potential $V$.

The following theorem summarizes the most important spectral results for this operator. A key role is played by the \emph{monodromy matrix} $A(p;E)$ and the \emph{discriminant}, which is defined by
\begin{equation}\label{e.discriminant}
\Delta(E) = \mathrm{Tr} \, A(p;E).
\end{equation}
From the explicit form of $A(p;E)$, we see that $\Delta$ is a monic real polynomial of degree $p$.

\begin{theorem}\label{t.perspth}
Suppose $V : \Z \to \R$ is periodic with minimal period $p$.

{\rm (a)} If $\Delta(E) \in (-2,2)$, then $\|A(n;E)\|$ is bounded. If $\Delta(E) \in \{-2,2\}$, then $\|A(n;E)\|$ is linearly bounded. If $\Delta(E) \not\in [-2,2]$, then $\|A(n;E)\|$ grows exponentially.

{\rm (b)} If $D \in [-2,2]$, then all solutions of $\Delta(\cdot) = D$ are real. If $D \in (-2,2)$, then all roots of $\Delta(\cdot) = D$ are simple.

{\rm (c)} The spectrum of $H$ is given by
\begin{equation}\label{f.shdzinm22}
\sigma(H) = \{ z : \Delta(E) \in [-2,2] \}.
\end{equation}
It consists of $p$ compact intervals, $B_1, \ldots , B_p$, called \textit{bands}, which are obtained by taking the closure of the $p$ mutually disjoint open intervals whose union is $\Delta^{-1} \left( (-2,2) \right)$. Thus, there are $m \le p-1$ bounded open intervals that separate bands, called \textit{open gaps}, and $p-1-m$ points, where two bands overlap, called \textit{closed gaps}.

{\rm (e)} The canonical spectral measure of $H$ is purely absolutely continuous. An essential support of this measure is given by $\{ E : \Delta(E) \in (-2,2) \}$. For energies $E$ in this set, all solutions of \eqref{e.eve} are bounded.
\end{theorem}

\subsection{Limit-Periodic Potentials}\label{ss.limitper}

We say that $V \in \ell^\infty(\Z)$ is \emph{limit-periodic} if it belongs to the $\ell^\infty$-closure of the set of periodic points of $S$, that is, there exist a sequence $\{V_j\}$ in $\ell^\infty(\Z)$ and a sequence $\{p_j\}$ in $\Z_+$ such that $S^{p_j}V_j = V_j$ for every $j$ and $\lim_{j \to \infty} \|V - V_j\|_\infty = 0$.

Every periodic $V$ is limit-periodic, every limit-periodic $V$ is almost periodic, and a non-periodic almost periodic $V$ is limit-periodic if and only if its hull is totally disconnected. This leads us naturally to the following definition: A \emph{Cantor group} is an abelian topological group which is compact, totally disconnected, and perfect (i.e., it has no isolated points). Then the fundamental structure theorem for limit-periodic potentials is the following.

\begin{prop} \label{p.cantorgroup.char.lim.per}
A potential $V \in \ell^{\infty}(\Z)$ is limit-periodic if and only if it can be represented via
\begin{equation} \label{f.cantorgroup.char.lim.per}
V(n) = V_{\omega}(n) = f(T^n \omega),
\end{equation}
where $\Omega$ is a Cantor group, $\omega \in \Omega $, $f: \Omega \to \R $ is continuous, and $T : \Omega \to \Omega$ is a minimal translation, say $T = \cdot  + \alpha$.
\end{prop}

This point of view is particularly useful here because it allows us to separate the base dynamics from the sampling function, so that we can keep the former fixed and vary the latter. The theorems below show that, as the sampling function is varied, various kinds of spectral behavior can be observed. In particular, all basic spectral types are possible (with the minor caveat that the pure point spectrum result does need some additional assumptions on the Cantor group).

The first result, which is due to Avila \cite{A09}, shows that Cantor spectrum is generic in the sense that it holds, for fixed base dynamics, for a dense $G_\delta$ set of sampling functions. Moreover, the generic type of Cantor spectrum in the limit-periodic setting comes with zero Lebesgue measure.

\begin{theorem}\label{t.lp.cantorspec.generic}
Suppose that $\Omega$ is a Cantor group and $T : \Omega \to \Omega$ is a minimal translation. Then, there exists a dense $G_{\delta}$ subset $\mathcal{C} \subseteq C(\Omega,\R) $ so that for all $f \in \mathcal{C}$, the spectrum $\Sigma$ associated with the potentials \eqref{f.cantorgroup.char.lim.per} is a Cantor set of zero Lebesgue measure.
\end{theorem}

This is in some sense  contrary to the expectations that were prevalent in the early days of the study of limit-periodic operators. In fact, due to the limit-periodicity of the potentials, and hence the very strong sense in which these operators are approximated by periodic ones, which in turn have band spectrum and purely absolutely continuous spectral measures, much of the early effort had focused on proving purely absolutely continuous spectrum. The Cantor structure of the spectrum was an objective as well, but note that such Cantor sets must have positive Lebesgue measure if the spectral measures are absolutely continuous! The following theorem describes this scenario in our setting. It was shown in this form by Damanik and Gan in \cite{DGa11a}, but the result is in the spirit of results of Avron-Simon \cite{AS81}, Chulaevsky \cite{C81} and others from the 1980's, and its proof follows the line of reasoning from those papers quite closely.

\begin{theorem}\label{t.lp.acspec.dense}
Suppose that $\Omega$ is a Cantor group and $T : \Omega \to \Omega$ is a minimal translation. Then, there is a dense set $\mathcal{A} \subseteq C(\Omega,\R)$ such that for every $f \in \mathcal{A}$ and $\omega \in \Omega$, the spectrum of $H_{\omega}$ is a Cantor set of positive Lebesgue measure and $H_{\omega}$ has purely absolutely continuous spectrum.
\end{theorem}

One can even ensure, in the same generality, that the Cantor spectrum is homogeneous in the sense of Carleson \cite{C83b}; see \cite{F14}.

While absolute continuity does happen for a dense set of sampling functions, Theorem~\ref{t.lp.cantorspec.generic} implies that it cannot occur on a generic set. Again, observe that a zero-measure set cannot support any absolutely continuous measures. Indeed, the generic spectral type turns out to be singular continuous; compare \cite{A09, DGa11a}.

\begin{theorem}\label{t.lp.scspec.generic}
Suppose that $\Omega$ is a Cantor group and $T : \Omega \to \Omega$ is a minimal translation. Then, there exists a dense $G_{\delta}$ set $\mathcal{S} \subseteq C(\Omega,\R)$ such that for every $f \in \mathcal{S}$ and every $\omega \in \Omega$, the spectrum of $H_{\omega}$ is a Cantor set of zero Lebesgue measure and $H_{\omega}$ has purely singular continuous spectrum.
\end{theorem}

There is another dense set of sampling functions, where interesting spectral phenomena occur; compare \cite{A09, DG10}.

\begin{theorem}\label{t.lp.scspec.posle}
Suppose that $\Omega$ is a Cantor group and $T : \Omega \to \Omega$ is a minimal translation. Then, there exists a dense subset $\mathcal{H} \subseteq C(\Omega,\R) $ such that for all $\omega \in \Omega$, the spectrum of $H_{\omega}$ is a Cantor set having zero Hausdorff dimension and $H_{\omega} $ has purely singular continuous spectrum. Moreover, the Lyapunov exponent $L(E)$ is a positive continuous function of $E$.
\end{theorem}

Indeed, all statements in Theorems~\ref{t.lp.scspec.generic} and \ref{t.lp.scspec.posle} except for singular continuity were shown by Avila in \cite{A09}, and a proof of singular continuity was added by Damanik and Gan in \cite{DG10, DGa11a}. The main reason why sampling functions in $\mathcal{H}$ are of interest is that they provide counterexamples to a conjecture of Simon, who had conjectured that positive Lyapunov exponents imply positive-measure spectrum. As Theorem~\ref{t.lp.scspec.posle} shows this is as far from the truth as possible. In fact, positive Lyapunov exponents do not even imply that the spectrum has positive Hausdorff dimension.

As explained by Gan in \cite{G10}, Cantor groups that have minimal translations are procyclic groups. We can classify such Cantor groups by studying their frequency integer sets. Every Cantor group with a minimal translation has a unique maximal frequency integer set $S = \{n_k\} \subseteq \Z_+$ with the property that $n_{k+1}/n_k$ is prime for every $k$; see, for example, \cite{G10}. We say that condition $A$ holds if there exists some integer $m \ge 2$ such that for every $k$, we have $n_k < n_{k+1} \leq n^m_k,$ that is, $\log n_{k+1}/\log n_k$ is uniformly bounded. Cantor groups admitting a minimal translation and satisfying condition $A$ are easily seen to exist. For them, it is possible to show that the third basic spectral type may occur; compare \cite{DGa11}.

\begin{theorem}\label{t.lp.ppspec}
Suppose that $\Omega$ is a Cantor group and $T : \Omega \to \Omega$ is a minimal translation. Suppose further that condition $A$ holds. Then, there exists $f \in C(\Omega,\R)$ such that for every $\omega \in \Omega$, the spectrum of $H_{\omega}$ is pure point, and all eigenvectors decay exponentially. In fact, the exponential decay is uniform for all $\omega$'s and all energies.
\end{theorem}

This result is surprising for several reasons. First, as pointed out above, the early works were aiming for absolutely continuous spectrum, and this was motivated by limit-periodic potentials being well approximated by periodic potentials. Pure point spectrum is as far away from absolutely continuous spectrum as possible. Second, the strength of the localization result is startling. Such a scenario is often called uniform localization, and it was not clear if such a strong localization statement can ever hold. That it arises in the limit-periodic world is indeed quite surprising. All other operator families that are known to be localized (random potentials, strongly coupled quasi-periodic potentials or skew-shift potentials) are either not known to be uniformly localized or known to be not uniformly localized. Motivated by \cite{DGa11} and extending \cite{Jit96}, Han \cite{Han15} showed that phase uniformity is a general phenomenon in the context of uniform localization.

Even though the spectral measures of limit-periodic operators are generically singular continuous and hence a study of their transport exponents is potentially interesting, only few works have studied transport exponents in this context; see, for example, \cite{CO14, DLY14}.

In light of the results above, the following open problems (listed already in \cite{DGa11}) arise naturally.

\medskip

\textit{Problem 1.} Is it true that for $f$ from a suitable dense subset of $C(\Omega,\R)$, $H_\omega$ has pure point spectrum for (Haar-) almost every $\omega \in \Omega$?

\medskip

We already know that for generic $f \in C(\Omega,\R)$, $H_\omega$ has purely singular continuous spectrum for every $\omega \in \Omega$, and also that for $f$ from a suitable dense subset of $C(\Omega,\R)$, $H_\omega$ has purely absolutely continuous spectrum for every $\omega \in \Omega$. Thus, an affirmative answer to Problem~1 would clarify the effect of the choice of $f$ on the spectral type. Since the methods leading to Theorem~\ref{t.lp.ppspec} are essentially restricted to large potentials, one should not expect them to yield an answer to Problem~1 and one should in fact pursue methods involving some randomness aspect.

\medskip

\textit{Problem 2.} Is the spectral type of $H_\omega$ always the same for every $\omega \in \Omega$?

\medskip

For quasi-periodic potentials, this is known not to be the case; see below. However, the mutual approximation by translates for two given elements in the hull is stronger in the limit-periodic case than in the quasi-periodic case, so it is not clear if similar counterexamples to uniform spectral types exist in the limit-periodic world.

Another related problem is the following:

\medskip

\textit{Problem 3.} Is the spectral type of $H_\omega$ always pure?

\medskip

Again, in the quasi-periodic world, this is known not to be the case: there are examples that have both absolutely continuous spectrum and point spectrum (cf., e.g., \cite{B06, B02b, FK02}).

Returning to the issue of point spectrum, one interesting aspect of a result stated (in the continuum case), but not proved, by Molchanov and Chulaevsky in \cite{MC84} is the coexistence of pure point spectrum with the absence of non-uniform hyperbolicity. That is, in their examples, the Lyapunov exponent vanishes on the spectrum and yet the spectral measures are pure point. This is the only known example of this kind and it would therefore be of interest to have a complete published proof of a result exhibiting this phenomenon. Especially since our study is carried out in a different framework, we ask within this framework the following question:

\medskip

\textit{Problem 4.} For which $f \in C(\Omega,\R)$ does the Lyapunov exponent vanish throughout the spectrum and yet $H_\omega$ has pure point spectrum for (almost) every $\omega \in \Omega$?

\medskip

Given the existing ideas, it is conceivable that Problems~1 and 4 are closely related and may be answered by the same construction. If this is the case, it will then still be of interest to show for a dense set of $f$'s that there is almost sure pure point spectrum with \emph{positive} Lyapunov exponents.

\subsection{Quasi-Periodic Potentials}\label{ss.quasiper}

Quasi-periodic potentials are generated by a minimal translation on a finite-dimensional torus and a continuous sampling function. That is, a quasi-periodic potential is of the form
\begin{equation}\label{e.qppotential}
V(n) = f(\omega + n \alpha),
\end{equation}
where $\alpha, \omega \in \T^d$ and $f \in C(\T^d,\R)$. The vector $\alpha = (\alpha_1, \ldots, \alpha_d) \in \T^d$ is assumed to be such that
\begin{equation}\label{e.alphaassumption}
k_j \in \Z, \, 1 \le j \le d, \quad \sum_{j = 1}^d k_j \alpha_j = 0 \in \T \quad \Rightarrow \quad k_j = 0, \, 1 \le j \le d.
\end{equation}
In other words, the entries of $\alpha$ together with $1$ are linearly independent over the rational numbers. This condition is equivalent to the translation by $\alpha$ on $\T^d$ being minimal.

The spectral theory of quasi-periodic Schr\"odinger operators has been extensively studied. In this subsection we will focus on some highlights, but won't attempt to give a comprehensive survey of the relevant literature. The main reason is that there is another contemporary survey by Jitomirskaya and Marx \cite{JM14} that focuses exclusively on the quasi-periodic case, and we refer the reader to that paper for more information.

\subsubsection{The Spectral Type}

The spectral properties of quasi-periodic Schr\"odinger operators are affected by the regularity of the sampling function $f$. For example, in the low-regularity regime, having purely singular continuous spectrum is typical, while in the strong-regularity regime, the absence of singular continuous spectrum is typical.

Let us make these statements more precise. We begin with the low-regularity setting. That is, nothing beyond continuity is assumed. Specializing Theorem~\ref{t.avdamthm1} to the case at hand, we obtain the following.

\begin{theorem}\label{t.qpgenericsingularspectrum}
Suppose $\alpha \in \T^d$ obeys \eqref{e.alphaassumption}. Then, there is a residual set $\mathcal{F}_\mathrm{s} \subseteq C(\T^d,\R)$ such that for every $f \in \mathcal{F}_\mathrm{s}$ and every $\omega \in \T^d$, the Schr\"odinger operator with potential \eqref{e.qppotential} has purely singular spectrum.
\end{theorem}

Recall that this is really a consequence of Kotani theory and holds in much greater generality. A complementary result was obtained by Boshernitzan and Damanik in \cite{BD08}.

\begin{theorem}\label{t.qpgenericcontinuousspectrum}
Suppose $\alpha \in \T^d$ obeys \eqref{e.alphaassumption}. Then, there is a residual set $\mathcal{F}_\mathrm{c} \subseteq C(\T^d,\R)$ such that for every $f \in \mathcal{F}_\mathrm{c}$ and Lebesgue almost every $\omega \in \T^d$, the Schr\"odinger operator with potential \eqref{e.qppotential} has purely continuous spectrum.
\end{theorem}

This result, too, holds in greater generality, albeit not quite as broadly as the previous theorem. See \cite{BD08} for the scope of their method, which includes in particular the skew-shift case for which the generic absence of eigenvalues was surprising at the time.

Combining Theorems~\ref{t.qpgenericsingularspectrum} and \ref{t.qpgenericcontinuousspectrum}, we obtain the generic singular continuity result mentioned above.

\begin{coro}\label{c.qpgenericsingularcontinuousspectrum}
Suppose $\alpha \in \T^d$ obeys \eqref{e.alphaassumption}. Then, there is a residual set $\mathcal{F}_\mathrm{sc} \subseteq C(\T^d,\R)$ such that for every $f \in \mathcal{F}_\mathrm{sc}$ and Lebesgue almost every $\omega \in \T^d$, the Schr\"odinger operator with potential \eqref{e.qppotential} has purely singular continuous spectrum.
\end{coro}

Let us now discuss the spectral type for analytic sampling functions. We consider mainly the one-frequency case, that is, $\alpha \in \T$. It turns out that the partition of the spectrum as $\Sigma = \mathcal{Z} \sqcup \mathcal{NUH}$ is particularly relevant. That is, pure point spectrum is typical in $\mathcal{NUH}$, whereas purely absolutely continuous spectrum is typical in $\mathcal{Z}$. Here is a result of Bourgain and Goldstein \cite{BG00} on spectral localization in $\mathcal{NUH}$. Recall that the Lyapunov exponents $L(E)$ depends not only on the energy $E$, but also on the sampling function $f$ and the frequency $\alpha$.

\begin{theorem}\label{t.bourgaingold}
Assume that $f$ is a $1$-periodic real-analytic function and that the Lyapunov exponent is strictly positive for any $\alpha \in \T \setminus \Q$ and any $E \in \R$. Then, for Lebesgue almost all $(\alpha,\omega) \in \T^2$, the Schr\"odinger operator with potential \eqref{e.qppotential} has pure point spectrum with exponentially decaying eigenfunctions.
\end{theorem}

Bourgain and Jitomirskaya showed in \cite{BJ00} that in the setting of the previous theorem, dynamical localization holds as well.

\begin{theorem}\label{t.qpdynloc}
Assume that $f$ is a $1$-periodic real-analytic function and that the Lyapunov exponent is strictly positive for any $\alpha \in \T \setminus \Q$ and any $E \in \R$. Then, for Lebesgue almost all $(\alpha,\omega) \in \T^2$, the Schr\"odinger operator with potential \eqref{e.qppotential} is dynamically localized in the sense \eqref{e.dynamicallocalization}, that is, we have
$$
\sup_t \sum_{n \in \Z} |n|^p | \langle \delta_n , e^{-itH} \delta_0 \rangle |^2 < \infty
$$
for every $p > 0$.
\end{theorem}

What about the input to Theorem~\ref{t.bourgaingold}? This is provided by a theorem due to Sorets and Spencer \cite{SS91}, which was already known at the time \cite{BG00} was published.

\begin{theorem}\label{t.soretsspencer}
Assume that $g$ is a non-constant $1$-periodic real-analytic function. Then, there exists $\lambda_0 > 0$ such that the following holds for $f = \lambda g$ with $\lambda > \lambda_0$. The Lyapunov exponent associated with the Schr\"odinger operator with potential \eqref{e.qppotential} is strictly positive for any $\alpha \in \T \setminus \Q$ and any $E \in \R$.
\end{theorem}

Thus, combining Theorems~\ref{t.bourgaingold}--\ref{t.soretsspencer}, we see that for analytic sampling functions, (spectral and dynamical) localization occurs at sufficiently large coupling for almost all frequencies. It is not an artifact that a zero-measure set of frequencies has to be excluded. Indeed, for frequencies that are sufficiently well approximated by rational numbers, spectral localization (and hence also dynamical localization) fails due to a result of Gordon \cite{G76}, which was highlighted by Avron and Simon \cite{AS83}:

\begin{theorem}\label{t.godronavronsimon}
Assume that $f \in C(\T,\R)$ and $\alpha \in \T \setminus \Q$ is such that for suitable rational numbers $\{ \frac{p_k}{q_k} \}_{k \ge 1}$, we have
\begin{equation}\label{e.liouville}
\left| \alpha - \frac{p_k}{q_k} \right| \le k^{-q_k}
\end{equation}
for $k \ge 1$. Then, for every $\omega \in \T$, the Schr\"odinger operator with potential \eqref{e.qppotential} has purely continuous spectrum.
\end{theorem}

Namely, it is not hard to see that under the assumption of Theorem~\ref{t.godronavronsimon}, the potential $V$ defined by \eqref{e.qppotential} is a Gordon potential for every $\omega \in \T$. Lemma~\ref{l.gordonpotential} then yields the conclusion.

Irrational $\alpha$ obeying \eqref{e.liouville} form a specific explicit class of \emph{Liouville numbers}. It is easy to see that this set is a dense $G_\delta$ set of zero Lebesgue measure. Irrational numbers that are not well approximated by rational numbers are called \emph{Diophantine}. As with Liouville numbers, there are several ways of imposing a Diophantine condition, some of which will lead to sets of full Lebesgue measure. It is an interesting open problem to extend the Bourgain-Goldstein localization result to an explicit full-measure set of Diophantine frequencies. That is, is there a full-measure set of Diophantine $\alpha$ for which spectral localization holds for Lebesgue almost all $\omega \in \T$, assuming that $f$ is a $1$-periodic real-analytic function for which the Lyapunov exponent is strictly positive for every $E \in \R$?

While localization is typical for analytic one-frequency quasi-periodic Schr\"odinger operators in the large coupling regime, purely absolutely continuous spectrum occurs in the weak-coupling regime, as shown by Bourgain and Jitomirskaya \cite{BJ02b}. This result actually does have an explicit Diophantine condition that is imposed on $\alpha$. Let us denote the distance from $0$ in $\T$ by $\| \cdot \|_{\T}$.

\begin{theorem}\label{t.bourgainjitoac}
Assume that $g$ is a $1$-periodic real-analytic function. Then, there exists $\lambda_1 > 0$ such that the following holds for $f = \lambda g$ with $0 < \lambda < \lambda_1$. If $\alpha \in \T \setminus \Q$ is Diophantine in the sense that
$$
\exists c > 0, \, r > 1 \text{ such that } \|n\alpha/2\|_{\T} > \frac{c}{|n|^r} \text{ for every } n \in \Z \setminus \{ 0 \},
$$
then for Lebesgue almost all $\omega \in \T$, the Schr\"odinger operator with potential \eqref{e.qppotential} has purely absolutely continuous spectrum.
\end{theorem}

This clarifies the typical (in the frequency and the phase) behavior for analytic one-frequency quasi-periodic Schr\"odinger operators in the regime of large and small coupling. In general, there is a significant gap between the two regimes, that is, the numbers $\lambda_0$ and $\lambda_1$ in the theorems above will be far apart. An exception is given by the almost Mathieu case, $g(\omega) = 2 \cos(2\pi \omega)$, which will be discussed in Section~\ref{s.6}. In this special case, the statements above actually hold with $\lambda_0 = \lambda_1 = 1$.

This gap was filled to a large extent when Avila developed his global theory of analytic quasi-periodic one-frequency Schr\"odinger operators in \cite{A14c, A14d, A14e, A14g}.\footnote{In fact, this work was presented in Avila's Fields Medalist lecture at the 2014 ICM in Seoul.} Recall that we can view the spectrum as a disjoint union of sets of energies, $\Sigma = \mathcal{NUH} \sqcup \mathcal{Z}$; see Theorem~\ref{t.spectrumandenergypartition}. At least for Diophantine frequencies, we also know that localization phenomena occur in $\mathcal{NUH}$; this follows from localized (in the energy parameter) versions of Theorems~\ref{t.bourgaingold} and \ref{t.qpdynloc}. Thus, we would like to understand the spectral type in $\mathcal{Z}$. By Theorem~\ref{t.ipkthm} we have $\Sigma_\mathrm{ac} = \overline{ \mathcal{Z} }^\mathrm{ess}$, but this leaves the question open of whether $\mathcal{Z}$ can locally have portions of zero Lebesgue measure or whether there can be any additional singular spectrum even when $\mathcal{Z}$ has (everywhere) positive measure. As we will see when we discuss the almost Mathieu operator, both phenomena can actually occur in the context of analytic quasi-periodic one-frequency operators. These issues were addressed by Avila using a further decomposition of $\mathcal{Z}$ into two subsets, namely the subcritical energies and the critical energies. These notions are defined via cocycle behavior, and more concretely by what happens when the phase $\omega \in \T$ is complexified. Namely, an energy in $\mathcal{Z}$ is subcritical if the Lyapunov exponent remains zero for sufficiently small perturbations of the phase in the imaginary direction, and critical otherwise. Fixing the frequency $\alpha$, Avila showed in \cite{A14d} that for a typical analytic $f$,\footnote{Here, ``typical'' is meant in the measure theoretical sense of prevalence.} there are no critical energies, and hence the spectrum splits into a localized regime and a subcritical regime. By showing that subcriticality implies almost-reducibility (which for some time was referred to as the ``almost reducibility conjecture'' (cf.~\cite{AJ10}), at least until it was proved in \cite{A14e, A14g}), Avila was then able to show that the subcritical regime in fact must be purely absolutely continuous. As a net result, one obtains that for a typical analytic one-frequency potential, there is no singular continuous spectrum and the decomposition $\Sigma = \mathcal{NUH} \sqcup \mathcal{Z}$ corresponds precisely to the decomposition into a localized regime and an absolutely continuous regime.

\bigskip

Much of the work on quasi-periodic potentials has focused on cases of extremal regularity, that is, analytic sampling functions and merely continuous sampling functions. For some work on quasi-periodic potentials of intermediate regularity, we refer to reader to \cite{Bj05, FJZ10, JM14b, JN11, K05, WZ14, WZ14b}) and references therein.

Similarly, while we have limited our discussion of the analytic category above to the one-frequency case, for results on the multi-frequency case we refer the reader to \cite{Bo05, E92, GS01, HA09, JM14} and references therein.

\subsubsection{Cantor Spectrum}

Another topic of wide interest is Cantor spectrum. For continuous sampling functions, this spectral phenomenon also turns out to be generic. Indeed, specializing Theorem~\ref{t.abdthm} by Avila, Bochi, and Damanik \cite{ABD09} to the case at hand, we obtain:

\begin{theorem}\label{t.qpgenericcantorspectrum}
Suppose $\alpha \in \T^d$ obeys \eqref{e.alphaassumption}. Then, there is a residual set $\mathcal{F}_\mathrm{cantor} \subseteq C(\T^d,\R)$ such that for every $f \in \mathcal{F}_\mathrm{cantor}$ and every $\omega \in \T^d$, the spectrum of the Schr\"odinger operator with potential \eqref{e.qppotential}  is a Cantor set.
\end{theorem}

For analytic sampling functions and in the regime of positive Lyapunov exponents, Cantor spectrum is typical as well as shown by Goldstein and Schlag \cite{GS11}.

\begin{theorem}
Denote
$$
\mathrm{Dioph} = \Big\{ \alpha \in \T : \exists c > 0, \, r > 1 \text{ such that } \|n \alpha\|_{\T} \ge \frac{c}{n (\log n)^{r}} \, \forall n > 1 \Big\}.
$$
Assume that $f$ is a $1$-periodic real-analytic function and that the Lyapunov exponent is strictly positive for any $\alpha \in (\alpha_1,\alpha_2)$ and any $E \in (E_1,E_2)$. Then there exists a set $B \subset \T$ of Hausdorff dimension zero such that for any $\alpha \in \mathrm{Dioph} \setminus B$, the intersection of $(E_1,E_2)$ with the spectrum of the Schr\"odinger operator with potential \eqref{e.qppotential} is a Cantor set.
\end{theorem}

On the other hand, Avila and Jitomirskaya showed the following result in the analytic non-perturbative small coupling regime \cite{AJ10}.

\begin{theorem}\label{t.avilajitocantor}
For typical {\rm (}i.e., outside a suitable set of infinite codimension{\rm )} $1$-periodic real-analytic $g$, there exists $\lambda_2 > 0$ such that the following holds for $f = \lambda g$ with $0 < \lambda < \lambda_2$. If $\alpha \in \T \setminus \Q$ is Diophantine in the sense that
$$
\exists c > 0, \, r > 1 \text{ such that } \|n\alpha/2\|_{\T} > \frac{c}{|n|^r} \text{ for every } n \in \Z \setminus \{ 0 \},
$$
then for Lebesgue almost all $\omega \in \T$, the Schr\"odinger operator with potential \eqref{e.qppotential} has Cantor spectrum.
\end{theorem}

In fact the stronger statement that all gaps allowed by the gap labeling theorem are open is shown.

For intermediate regularity, see \cite{Si87, WZ14b} for results on Cantor spectra.

\section{The Almost Mathieu Operator}\label{s.6}

We will discuss the almost Mathieu operator
$$
[H^{\lambda,\alpha}_\omega \psi](n) = \psi(n+1) + \psi(n-1) + 2 \lambda \cos (2 \pi (\omega + n \alpha)) \psi(n).
$$
This special case of a quasi-periodic Schr\"odinger operator deserves a separate section for a number of reasons. It it the single case of a quasi-periodic Schr\"odinger operator that has been more or less completely analyzed. The wealth of the results and the sheer number of papers devoted to this operator are quite impressive. Moreover, even in this single family, one can already see, for example, that all possible spectral types may arise in the quasi-periodic context and one can also see the mechanisms behind these phenomena. Related to this, the study of the almost Mathieu case has informed the study of the general quasi-periodic case. Many of the known results in the general (predominantly analytic) setting are extensions of results earlier obtained for the almost Mathieu case.

We see that the almost Mathieu operator is a one-frequency quasi-periodic Schr\"odinger operator, where the sampling function is given by the trigonometric polynomial
$$
f(\omega) = 2 \lambda \cos (2 \pi \omega).
$$
As above, we consider $\alpha$ and $\omega$ as elements of $\T = \R / \Z$. It is easy to see that $H^{\lambda,\alpha}_\omega = H^{-\lambda,\alpha}_{\omega + \frac{1}{2}}$ and hence we may focus on the case $\lambda > 0$. If $\alpha$ is irrational, then the spectrum of $H^{\lambda,\alpha}_\omega$ is independent of $\omega$ and may be denoted by $\Sigma^{\lambda,\alpha}$. In fact, when one talks about the almost Mathieu operator, it is implicitly assumed that $\alpha$ is irrational. However, it is sometimes useful to consider rational approximations of $\alpha$ and hence periodic approximations of the quasi-periodic operator. In the general case, we will set
\begin{equation}\label{e.amospecunion}
\Sigma^{\lambda,\alpha} = \bigcup_{\omega \in \T} \sigma(H^{\lambda,\alpha}_\omega),
\end{equation}
and this definition agrees with the one above in the irrational case.

Much of the development of the theory of the almost Mathieu operator has been driven by three conjectures, which have been around since the late 1970's/early 1980's \cite{AA80, S82}.

\begin{amoprob}[Measure of the Spectrum]\label{prob.amo1}
For every $\lambda > 0$ and every irrational $\alpha \in \T$, we have
$$
\mathrm{Leb}(\Sigma^{\lambda,\alpha}) = 4 | 1 - \lambda |.
$$
\end{amoprob}

\begin{amoprob}[Metal-Insulator Transition]\label{prob.amo2}
Suppose $\lambda > 0$, $\alpha \in \T$ is irrational, and $\omega \in \T$. Then,
\begin{itemize}

\item $H^{\lambda,\alpha}_\omega$ has purely absolutely continuous spectrum if $\lambda < 1$,

\item $H^{\lambda,\alpha}_\omega$ has purely singular continuous spectrum if $\lambda = 1$,

\item $H^{\lambda,\alpha}_\omega$ is spectrally localized if $\lambda > 1$.

\end{itemize}
\end{amoprob}

\begin{amoprob}[Ten Martini Problem]\label{prob.amo3}
For every $\lambda > 0$ and every irrational $\alpha \in \T$, $\Sigma^{\lambda,\alpha}_\omega$ is a Cantor set.
\end{amoprob}

\subsection{The Main Results}

The three theorems below concern the Lebesgue measure of the spectrum, the metal-insulator transition, and the ten Martini problem. They are stated in the generality in which they are currently known and summarize the results of many authors, obtained over the course of about three decades.

\begin{theorem}[Measure of the Spectrum]\label{t.amomeasureofspectrum}
For every $\lambda > 0$ and every irrational $\alpha \in \T$, we have
\begin{equation}\label{e.amoprobmeasure}
\mathrm{Leb}(\Sigma^{\lambda,\alpha}) = 4 | 1 - \lambda |.
\end{equation}
\end{theorem}

This shows that the original conjecture holds in complete generality. This result was established in the papers \cite{AK06, AMS90, JK02, JL98, L93, L94}.

\begin{theorem}[Metal-Insulator Transition]\label{t.amospectraltype}
{\rm (a)} If $\lambda < 1$, then for every $\alpha$ and every $\omega$, the spectrum is purely absolutely continuous.\\
{\rm (b)} If $\lambda = 1$, then for every irrational $\alpha$ and all but countably many $\omega$, the spectrum is purely singular continuous.\\
{\rm (c)} If $\lambda > 1$, then for almost every $\alpha$ and almost every $\omega$, the spectrum is pure point and the eigenfunctions decay exponentially. \\
{\rm (d)} If $\lambda > 1$, then for generic $\alpha$ and every $\omega$, the spectrum is purely singular continuous.\\
{\rm (e)} If $\lambda > 1$, then for every irrational $\alpha$ and generic $\omega$, the spectrum is purely singular continuous.
\end{theorem}

This shows that the original conjecture holds in a full measure sense, but fails in a generic sense for $\lambda > 1$. The situation at $\lambda = 1$ is not completely resolved yet, and it is still expected that one always has purely singular continuous spectrum in this case. The theorem above combines results from \cite{A14a, A14b, AD08, AK06, AS83, G76, GJLS97, J99, JS94, L94}.

\begin{theorem}[Ten Martini Problem]\label{t.amotenmartini}
For every $\lambda > 0$ and every irrational $\alpha \in \T$, $\Sigma^{\lambda,\alpha}$ is a Cantor set.
\end{theorem}

This shows that also in this case the original conjecture holds in complete generality. The relevant papers are \cite{AJ09, AK06, BS82, CEY90, J99, L94, P04}.

In the following subsections we will present some of the main ideas that go into the proof of these theorems.

\subsection{Aubry Duality}

Consider the Hilbert space $L^2(\T \times \Z)$ and the operator $H^{\lambda,\alpha} : L^2(\T \times \Z) \to L^2(\T \times \Z)$
given by
$$
[H^{\lambda,\alpha} \varphi](\omega,n) = \varphi(\omega,n+1) + \varphi(\omega,n-1) + 2 \lambda \cos (2 \pi (\omega + n \alpha))
\varphi(\omega,n).
$$
Introduce the duality transform $\mathcal{A} : L^2(\T \times \Z) \to L^2(\T \times \Z)$, which is given by
$$
[\mathcal{A} \varphi](\omega,n) = \sum_{m \in \Z} \int_\T e^{-2\pi
i (\omega + n \alpha)m} e^{-2\pi i n \eta} \varphi(\eta,m) \,
d\eta.
$$
This definition assumes initially that $\varphi$ is such that the sum in $m$ converges, but note that in terms of the Fourier
transform on $L^2(\T \times \Z)$, we have $[\mathcal{A} \varphi](\omega,n) = \hat \varphi (n, \omega + n \alpha)$, which may be used to extend the definition to all of $L^2(\T \times \Z)$ and shows that $\mathcal{A}$ is unitary.

A first consequence of Aubry duality is a formula relating the spectra of $H^{\lambda,\alpha}_\omega$ and $H^{\lambda^{-1},\alpha}_\omega$, as shown by Avron and Simon in \cite{AS83}.

\begin{theorem}\label{t.dualspectrum}
We have $\Sigma^{\lambda,\alpha} = \lambda \Sigma^{\lambda^{-1},\alpha}$.
\end{theorem}

Moreover, Gordon-Jitomirskaya-Last-Simon stated the following theorem in \cite{GJLS97}, which relates the type of the spectral measures at $\lambda$ and the dual coupling $\lambda^{-1}$.

\begin{theorem}\label{t.gjlsdual}
Suppose $\lambda > 0$ and $\alpha \in \T$ is irrational.\\
{\rm (a)} We have $H^{\lambda,\alpha} \mathcal{A} = \lambda \mathcal{A} H^{\lambda^{-1},\alpha}$.
\\
{\rm (b)} If $H^{\lambda,\alpha}_\omega$ has pure point spectrum for almost every $\omega \in \T$, then $H^{\lambda^{-1},\alpha}_\omega$ has purely absolutely continuous spectrum for almost every $\omega \in \T$.
\\
{\rm (c)} If $H^{\lambda,\alpha}_\omega$ has some point spectrum for almost every $\omega \in \T$, then $H^{\lambda^{-1},\alpha}_\omega$ has some absolutely continuous spectrum for almost every $\omega \in \T$.
\end{theorem}

For many years this theorem had been the basis of employing Aubry duality to relate point spectrum and absolutely continuous spectrum. However, it turned out that the proof given in \cite{GJLS97} does not actually establish Theorem~\ref{t.gjlsdual}. Namely, the proof relies on a statement, for which it quotes the paper \cite{DS83} by Deift and Simon, which is not actually proved in \cite{DS83}.  Instead, one may rely on duality arguments given in \cite{A14b, AJ10, JM12b} in order to show that Theorem~\ref{t.gjlsdual} actually does hold as formulated.

\subsection{The Herman Estimate}

Herman proved the following lower bound for the Lyapunov exponent in \cite{H83}.

\begin{theorem}\label{t.hermanestimate}
If $\lambda > 0$ and $\alpha \in \T$ is irrational, then the Lyapunov exponent associated with the almost Mathieu operator satisfies $L(E) \ge \log \lambda$ for every $E$.
\end{theorem}

The key idea of the proof is to complexify the phase $\omega$ and to employ subharmonicity in this new variable. Since the argument is so elegant and short, let us give it here. Setting $w = e^{2 \pi i\omega}$, we see that
$$
2 \lambda \cos(2\pi (\omega + m \alpha)) = \lambda \left( e^{2 \pi i \alpha m} w
  + e^{-2 \pi i \alpha m} w^{-1} \right) .
$$
Thus, the one-step transfer matrices have the form
$$
T_\omega(m; E ) = \left( \begin{array}{cr} E - \lambda \left( e^{2 \pi i \alpha m} w + e^{-2 \pi i \alpha m} w^{-1} \right) & -1 \\ 1 & 0 \end{array} \right)
$$
If we define
$$
N_n(w) = w^n A_\omega(n;E) = (w T_\omega(n; E)) \cdots (w T_\omega(1; E)),
$$
initially on $|w| = 1$, we see that $N_n$ extends to an entire function and hence $w \mapsto \log \| N_n(w) \|$ is subharmonic. Thus,
$$
\int_0^1 \log \| N_n(e^{2 \pi i\omega}) \| \, d\omega \ge \log \| N_n(0) \| = n \log \lambda.
$$
Moreover, $\|N_n(e^{2 \pi i\omega})\| = \|A_\omega(n;E)\|$. Thus,
$$
L(E) = \lim_{n \to \infty} \frac1n \int_\T \log \| A_\omega(n;E) \| \, d\omega = \lim_{n \to \infty} \frac1n \int_\T \log \| N_n(e^{2 \pi i\omega}) \| \, d\omega \ge \log \lambda ,
$$
and Theorem~\ref{t.hermanestimate} follows.

The argument above extends readily to trigonometric polynomials, and also to multi-frequency models. In this sense the Herman estimate was the precursor to the result by Sorets-Spencer (Theorem~\ref{t.soretsspencer}) and its extension to the multi-frequency case by Bourgain \cite{Bo05}.

\medskip

Since the Lyapunov exponent is non-negative, the Herman estimate may be rewritten as $L(E) \ge \max \{ \log \lambda , 0 \}$. This estimate in turn is sharp, as shown by Bourgain and Jitomirskaya in \cite{BJ02}:

\begin{theorem}
If $\lambda > 0$ and $\alpha$ is irrational, then $L(E) = \max \{ \log \lambda , 0 \}$ for every $E \in
\Sigma^{\lambda,\alpha}$.
\end{theorem}

\subsection{The Measure of the Spectrum}\label{ss.amomeasure}

Helffer and Sj\"ostrand proved the following in \cite{HS89}.

\begin{theorem}
There is a constant $A < \infty$ such that for every $\alpha \in \T$ irrational with continued fraction coefficients obeying $a_k \ge A$ for every $k \in \Z_+$, $\mathrm{Leb}(\Sigma^{1,\alpha}) = 0$.
\end{theorem}

In fact, they obtain a detailed description of the quantitative self-similarity properties of the spectrum. It would be very interesting to extend this work to non-critical coupling.

A different approach to studying the measure of the spectrum is based on periodic approximations, obtained by replacing the irrational frequency of the quasiperiodic potential with suitable rational numbers nearby. To make this approach effective, one needs good quantitative information about the periodic operators and a suitable quantitative continuity statement. The starting point is the paper \cite{AMS90} by Avron, van Mouche, and Simon, which is devoted to both issues, namely a study of rational frequencies and a quantitative continuity result.

In addition to the union of spectra \eqref{e.amospecunion}, consider also the intersection of spectra,
\begin{equation}\label{e.amospecintersection}
\sigma^{\lambda,\alpha} = \bigcap_{\omega \in \T} \sigma(H^{\lambda,\alpha}_\omega).
\end{equation}
Then, the following is shown in \cite{AMS90} (see also \cite{JM12} for interesting follow-up work).

\begin{theorem}\label{t.avmsthm}
Suppose $\alpha \in \T$ is rational. Write $\alpha = p/q$ with $(p,q) = 1$. \\[1mm]
{\rm (a)} We have
$$
\mathrm{Leb}\left(\sigma^{\lambda,\alpha}\right) = \begin{cases} 4 | 1 - \lambda | & \text{if } 0 < \lambda < 1, \\ 0 & \text{if } \lambda \ge 1. \end{cases}
$$
{\rm (b)} For $0 < \lambda < 1$, we have
$$
\mathrm{Leb}\left(\sigma^{\lambda,\alpha}\right) \le \mathrm{Leb}\left(\Sigma^{\lambda,\alpha}\right) \le \mathrm{Leb}\left(\sigma^{\lambda,\alpha}\right) + 4\pi \lambda^{q/2}.
$$
\end{theorem}

Notice that if $(p_k,q_k) = 1$ and $q_k \to \infty$, then Theorem~\ref{t.avmsthm} yields
\begin{equation}\label{e.avmsconv}
\lim_{k \to \infty} \mathrm{Leb}\left(\Sigma^{\lambda,\frac{p_k}{q_k}}\right) = 4 | 1 - \lambda |
\end{equation}
if $0 < \lambda < 1$. Thus, if $\alpha \in \T$ is irrational and $\frac{p_k}{q_k}$ are the continued fraction approximants, then \eqref{e.avmsconv} suggests strongly that \eqref{e.amoprobmeasure} holds for $0 < \lambda < 1$. But then \eqref{e.amoprobmeasure} will also hold for $\lambda > 1$ due to Theorem~\ref{t.dualspectrum}.

We see that proving \eqref{e.amoprobmeasure} using \eqref{e.avmsconv} requires a continuity result for $\alpha \mapsto \Sigma^{\lambda,\alpha}$. The authors of \cite{AMS90} also proved $\frac12$-H\"older continuity of this map with respect to the Hausdorff metric.

\begin{theorem}\label{t.avmsthm2}
Suppose $\lambda > 0$. Then, there exists a constant $C$ such that
$$
\mathrm{dist}_H \left( \Sigma^{\lambda,\alpha}, \Sigma^{\lambda,\alpha'} \right) \le C \left| \alpha - \alpha' \right|^{1/2}.
$$
\end{theorem}

In fact, this result holds for general one-frequency quasi-periodic Schr\"odinger operators with a $C^1$ sampling function $f$, and the constant $C$ may be chosen as $6 \|f'\|^{1/2}_\infty$.

While this seemingly put Avron, van Mouche, and Simon very close to proving \eqref{e.amoprobmeasure} for every $\lambda \not= 1$, it was only shown by Last a few years later how to derive the desired conclusion under a suitable additional assumption. Namely, Last proved the following in \cite{L93}.

\begin{theorem}\label{t.lastmeasure}
Suppose $\alpha \in \T$ is irrational with an unbounded continued fraction expansion. Equivalently, there are rational numbers $\frac{p_k}{q_k}$ such that
$$
\lim_{k \to \infty} q_k^2 \left| \alpha - \frac{p_k}{q_k} \right| = 0.
$$
Then, for every $\lambda \in (0,\infty) \setminus \{ 1 \}$, we have $\mathrm{Leb}(\Sigma^{\lambda,\alpha}) = 4 | 1 - \lambda |$.
\end{theorem}

The assumption in Theorem~\ref{t.lastmeasure} holds for Lebesgue almost every $\alpha$.

The work of Avron, van Mouche, and Simon seems to exclude the case of critical coupling, $\lambda = 1$. However, by approximating the critical coupling from below with non-critical values, Last \cite{L94} was nevertheless able to show the following estimates.

\begin{theorem}\label{t.lastavmsthm}
Suppose $\alpha = p/q \in \T$ is rational with $(p,q) = 1$. Then,
$$
\frac{2(\sqrt{5} + 1)}{q} \le \mathrm{Leb}\left(\Sigma^{1,\alpha}\right) \le \frac{8e}{q}.
$$
\end{theorem}

As before, rational approximation then leads to the desired result for Lebesgue almost every irrational $\alpha \in \T$, as also shown in \cite{L94}.

\begin{theorem}\label{t.lastzeromeasure}
Suppose $\alpha \in \T$ is irrational with an unbounded continued fraction expansion. Then, we have $\mathrm{Leb}(\Sigma^{1,\alpha}) = 0$.
\end{theorem}

The limitation in the frequencies covered by Theorems~\ref{t.lastmeasure} and \ref{t.lastzeromeasure} comes directly from the quantitative version of the continuity of spectra stated in Theorem~\ref{t.avmsthm2}. To cover more or even all irrational frequencies, an improved continuity statement was necessary. Such an improvement was obtained by Jitomirskaya and Last in \cite{JL98} for $\lambda > 14.5$ and by Jitomirskaya and Krasovsky in \cite{JK02} for $\lambda > 1$ (really as a consequence of $L(E) > 0$, which by the Herman estimate holds when $\lambda > 1$). Recall that if \eqref{e.amoprobmeasure} holds for some $\lambda$, it also holds for $\lambda^{-1}$. As a consequence, the following theorem was obtained in \cite{JK02}.

\begin{theorem}\label{t.jitokrasmeasure}
Suppose $\alpha \in \T$ is irrational. Then, for every $\lambda \in (0,\infty) \setminus \{ 1 \}$, we have $\mathrm{Leb}(\Sigma^{\lambda,\alpha}) = 4 | 1 - \lambda |$.
\end{theorem}

This resolves AMO-Problem~\ref{prob.amo1}, except for critical coupling and the zero-measure set of frequencies with bounded continued fraction expansion. Resolving this issue completely became one of the problems on Barry Simon's list of Schr\"odinger operator problems for the 21st century \cite{S00}. The complete solution to AMO-Problem~\ref{prob.amo1}, and hence Theorem~\ref{t.amomeasureofspectrum}, was finally obtained by Avila and Krikorian in \cite{AK06}. They considered critical coupling, $\lambda = 1$, and frequencies $\alpha \in \T$ that are recurrent Diophantine, that is, that are such that infinitely many of their iterates of the Gauss map (which truncates the continued fraction expansion, that is, $[0;a_1,a_2,a_3,\ldots]$ is sent to $[0;a_2,a_3,a_4\ldots]$) satisfy a suitable fixed Diophantine condition. Avila and Krikorian \cite{AK06} showed the following:

\begin{theorem}\label{t.avilakrikorianmeasure}
Suppose $\alpha \in \T$ is recurrent Diophantine. Then, we have $\mathrm{Leb}(\Sigma^{1,\alpha}) = 0$.
\end{theorem}

The set of recurrent Diophantine $\alpha \in \T$ has full Lebesgue measure and it includes all numbers with bounded continued fraction. In other words, the union of the frequencies covered by Last and the frequencies covered by Avila and Krikorian is equal to all irrational $\alpha \in \T$. Thus, combining Theorems~\ref{t.lastzeromeasure}--\ref{t.avilakrikorianmeasure}, Theorem~\ref{t.amomeasureofspectrum} follows.

\subsection{The Metal-Insulator Transition}\label{ss.mitransition}

Recall that AMO-Problem~\ref{prob.amo2} claims purely absolutely continuous spectrum for subcritical coupling $\lambda < 1$, purely singular continuous spectrum for critical coupling $\lambda = 1$, and spectral localization for supercritical coupling $\lambda > 1$.

The first result relevant to this problem was actually a negative one. Namely, specializing the results of Gordon \cite{G76} and Avron and Simon \cite{AS83} to the almost Mathieu case, Theorem~\ref{t.godronavronsimon} becomes:

\begin{theorem}\label{t.godronavronsimon2}
Assume $\alpha \in \T$ is such that for suitable rational numbers $\{ \frac{p_k}{q_k} \}_{k \ge 1}$, we have
$$%\begin{equation}\label{e.liouville2}
\left| \alpha - \frac{p_k}{q_k} \right| \le k^{-q_k}
$$%\end{equation}
for $k \ge 1$. Then, for every $\lambda > 0$ and every $\omega \in \T$, $H^{\lambda,\alpha}_\omega$ has purely continuous spectrum.
\end{theorem}

In particular, there is no spectral localization in the supercritical regime if $\alpha$ is Liouville in the sense above. Recall that the set of Liouville numbers has zero Lebesgue measure, but it is large in the sense that it is a dense $G_\delta$ set.

But even for typical frequencies, a correction to the expected result is necessary, as shown by Jitomirskaya and Simon in \cite{JS94}.

\begin{theorem}\label{t.jitosimon}
Assume $\lambda > 0$ and $\alpha \in \T$ is irrational. Then, for $\omega$'s from a dense $G_\delta$ subset of $\T$, $H^{\lambda,\alpha}_\omega$ has purely continuous spectrum.
\end{theorem}

The proof relies on the fact that the cosine function is even, and hence the potential of the almost Mathieu operator has long stretches on which it is almost reflection symmetric. This has consequences for the generalized eigenfunctions. If these stretches of almost symmetry are suitably located (which does happen for a generic set of phases $\omega \in \T$), one can in this way exclude the presence of square-summable solutions to the difference equation. Note, however, that the proof is indirect. Only assuming that a solution is square-summable, one can then show that it does not decay, and hence cannot be square-summable. The Gordon lemma, on the other hand, which relies on local almost translation symmetries, does exclude the presence of decaying solutions in an unqualified way and hence it can sometimes be used to establish even stronger continuity properties of spectral measures. We will see instances of this in later sections.

The next result related to AMO-Problem~\ref{prob.amo2} was obtained as a consequence of the work of Gordon, Jitomirskaya, Last, and Simon \cite{GJLS97} on their version of Aubry duality stated in Theorem~\ref{t.gjlsdual}.

\begin{theorem}\label{t.gjlssccrit}
Assume $\alpha \in \T$ is irrational. Then, for Lebesgue almost every $\omega \in \T$, $H^{1,\alpha}_\omega$ has purely singular continuous spectrum.
\end{theorem}

Due to the gap in the proof of Theorem~\ref{t.gjlsdual} discussed after this theorem was stated above, it is important to note that in the meantime a different proof of Theorem~\ref{t.gjlssccrit} has been found by Avila \cite{A14b}, establishing an even stronger result. (This paper also helped in showing that Theorem~\ref{t.gjlsdual} is actually true as stated.)

\begin{theorem}\label{t.avsc}
Suppose $\alpha \in \T$ is irrational. Then, $H_\omega^{1,\alpha}$ has purely singular continuous spectrum for all but countably many $\omega \in \T$.
\end{theorem}

Removing the exclusion of the countable set of $\omega$'s in Theorem~\ref{t.avsc} is an interesting open problem.

The next major milestone was Jitomirskaya's \cite{J99}. Recognizing the necessary restrictions imposed by Theorems~\ref{t.godronavronsimon2} and \ref{t.jitosimon} on $\alpha$ and $\omega$, she established the expected spectral localization result outside of these exceptions.

\begin{theorem}\label{t.jitoloc}
Suppose $\lambda > 1$, $\alpha \in \T$ is Diophantine in the sense that there are constants $c > 0, r > 1$ such that
$$
| \sin ( 2 \pi n \alpha ) | > \frac{c}{|n|^r} \quad \text{ for every } n \in \Z \setminus \{0\},
$$
and $\omega \in \T$ is non-resonant in the sense that the relation
$$
\left| \sin \left( 2\pi \left( \omega + \frac{n}{2} \alpha \right) \right) \right| < \exp \left( -|n|^\frac{1}{2r} \right)
$$
holds for at most finitely many $n \in \Z$. Then, $H_\omega^{\lambda,\alpha}$ has pure point spectrum with exponentially decaying eigenfunctions.
\end{theorem}

Each of the conditions on $\alpha$ and $\omega$ in Theorem~\ref{t.jitoloc} holds on a set of full Lebesgue measure. In particular, part (c) of Theorem~\ref{t.amospectraltype} follows.

The proof of Theorem~\ref{t.amospectraltype} rests entirely on the positivity of the Lyapunov exponent. That is, nothing specific is assumed about $\lambda$, and the condition on $\lambda$ in Theorem~\ref{t.amospectraltype} is merely a consequence of the Herman estimate, Theorem~\ref{t.hermanestimate}. The paper \cite{J99} (along with its predecessors \cite{J94, J95}) therefore introduced the concept of ``nonperturbative localization'' in which the proof of localization uses the positivity of the Lyapunov exponent as input rather than a largeness assumption on the coupling constant. Another feature of a non-perturbative result is that the largeness condition on the coupling constant, which arises implicitly here as just explained, is independent of the frequency. This is in contrast to perturbative results where the largeness condition indeed does depend on the frequency; compare, for example, \cite{E97, FSW90, Si87}. The concept of non-perturbative localization was further explored in the more general setting of analytic sampling functions by Bourgain and Goldstein \cite{BG00}; compare Theorem~\ref{t.bourgaingold}.

The importance of \cite{J99} goes beyond merely validating part (c) of Theorem~\ref{t.amospectraltype}. Indeed, applying Aubry duality to Theorem~\ref{t.jitoloc} yields several nice consequences. The obvious one, stated in \cite{J99} as a consequence of part (b) of Theorem~\ref{t.gjlsdual}, is:

\begin{theorem}\label{t.jitoac}
Suppose $\lambda < 1$ and $\alpha \in \T$ is Diophantine in the sense of the previous theorem. Then, $H_\omega^{\lambda,\alpha}$ has purely absolutely continuous spectrum for Lebesgue almost every $\omega \in \T$.
\end{theorem}

Exponential localization at energy $E$ immediately implies boundedness of all solutions for the dual model at energy $E/\lambda$, from which purely absolutely continuous spectrum for almost every phase follows immediately by subordinacy theory. Thus, absolutely continuous spectrum for almost every phase is an immediate corollary of exponential localization for the dual model. The paper \cite{J99} referred to \cite{GJLS97} for this conclusion, but in view of the problems mentioned in the discussion of \cite{GJLS97}, the simple reasoning just described should be implemented instead.

Going beyond that, Avila and Jitomirskaya developed a quantitative formulation of Aubry duality in \cite{AJ10} and used it to prove the following stronger result.

\begin{theorem}\label{t.avjitoac}
Suppose $\lambda < 1$ and $\alpha \in \T$ is Diophantine in the sense above. Then, $H_\omega^{\lambda,\alpha}$ has purely absolutely continuous spectrum for every $\omega \in \T$.
\end{theorem}

Moreover, recall that for Liouville $\alpha$, the dual model in the supercitical regime has purely singular continuous spectrum, and Aubry duality cannot predict what should happen in the subcritical regime for such frequencies. If anything, one might be tempted to expect purely singular continuous spectrum as well, as the case of critical coupling shows that the dual of singular continuous may be singular continuous; see Theorem~\ref{t.gjlssccrit}. Note, however, the the singular continuous spectra at supercritical coupling and at critical coupling are different animals, as one comes with positive Lyapunov exponents and the other one comes with zero Lyapunov exponents. In any event, clarifying the spectral type in the subcritical regime for Liouville frequencies remained a challenge that put this issue on Barry Simon's list of Schr\"odinger operator problems for the 21st century \cite{S00} as well.

The first step toward a complete understanding of the spectral type in the subcritical regime was actually taken in the supercritical regime. Given $\alpha \in \T$ irrational with continued fraction approximants $\frac{p_k}{q_k}$, let
$$
\beta(\alpha) := \limsup_{k \to \infty} \frac{\log q_{k+1}}{q_k}.
$$
Note that $\beta(\alpha) = 0$ if $\alpha$ is Diophantine and $\beta(\alpha) = \infty$ if $\alpha$ is Liouville. Frequencies $\alpha$ with $0 < \beta(\alpha) < \infty$ are of a weak Liouville type. By carefully examining the argument of \cite{J99}, it follows that in Theorem~\ref{t.jitoloc} the Diophantine condition can be replaced by the weaker condition $\beta(\alpha) = 0$. Applying Aubry duality, one can replace the Diophantine condition by $\beta(\alpha) = 0$ in Theorem~\ref{t.jitoac} as well.

Avila and Damanik then extended Theorem~\ref{t.jitoac} to the regime where $\beta(\alpha) > 0$. They showed the following in \cite{AD08}.

\begin{theorem}\label{t.avdamac}
Suppose $\lambda < 1$ and $\alpha \in \T$ is irrational with $\beta(\alpha) > 0$. Then, $H_\omega^{\lambda,\alpha}$ has purely absolutely continuous spectrum for almost every $\omega \in \T$.
\end{theorem}

Their proof relies on Corollary~\ref{c.pureaccoro}. That is, they proved that the density of states measure is purely absolutely continuous when $\lambda \not= 1$ and $\beta(\alpha) > 0$.\footnote{Of course, it is purely singular when $\lambda = 1$ due to Theorems~\ref{t.lastzeromeasure} and \ref{t.avilakrikorianmeasure}, and its absolute continuity when $\lambda \not=1$ and $\beta(\alpha) = 0$ follows from purely absolutely continuous spectrum and an application of Corollary~\ref{c.pureaccoro} in the other direction.} This application shows the importance of Corollary~\ref{c.pureaccoro} which had been somewhat overlooked until then.

All remaining cases (namely, the exceptional frequencies in \cite{AJ10} and the exceptional phases in \cite{AD08}) were finally handled by Avila in \cite{A14a}, and hence part (a) of Theorem~\ref{t.amospectraltype} as stated followed:

\begin{theorem}\label{t.avac}
Suppose $\lambda < 1$ and $\alpha \in \T$. Then, $H_\omega^{\lambda,\alpha}$ has purely absolutely continuous spectrum for every $\omega \in \T$.
\end{theorem}

This completes our summary of the metal-insulator transition for the almost Mathieu operator at the critical coupling $\lambda = 1$. In fact there is a second spectral transition for frequencies $\alpha$ with $0 < \beta(\alpha) < \infty$. Namely, for such frequencies, beside the known transition from absolutely continuous to singular continuous at $\lambda = 1$, there is another transition from singular continuous to spectrally localized at $\lambda = e^{\beta(\alpha)}$;  see \cite{AJ09, J95b, J07} for the conjecture and partial results, as well as \cite{AYZ14} for the full result. Note that this beautifully interpolates between the cases $\beta(\alpha) = 0$ and $\beta(\alpha) = \infty$, where there is no second transition. Thus, one can always state that the (typical) spectral type is absolutely continuous between $0$ and $1$, singular continuous between $1$ and $e^{\beta(\alpha)}$, and pure point between $e^{\beta(\alpha)}$ and $\infty$.

\subsection{Cantor Spectrum}\label{ss.tenmartini}

The first result on Cantor spectrum for the almost Mathieu operator was established by Bellissard and Simon \cite{BS82}.

\begin{theorem}\label{t.bstenmartini}
The set $\{ (\lambda,\alpha) : \Sigma^{\lambda,\alpha} \text{ is a Cantor set} \}$ is residual.
\end{theorem}

The proof is quite soft and uses the Baire category theorem. The conditions on $\lambda$ and $\alpha$ are not explicit. While the proof is such that the $\alpha$'s in question will be well approximated by rational numbers, no explicit class of Liouville numbers can be singled out for which Cantor spectrum follows from this approach.

The need to vary $\lambda$ was eliminated in a work of Choi, Elliott, and Yui \cite{CEY90}. They proved a Cantor spectrum result for fixed $\lambda$ and an explicit generic set of frequencies.

\begin{theorem}\label{t.ceytenmartini}
Suppose that $\lambda > 0$ and $\alpha \in \T$ is a Liouville number in the sense of \eqref{e.liouville}. Then, $\Sigma^{\lambda,\alpha}$ is a Cantor set.
\end{theorem}

The condition in \cite{CEY90} is actually more general. The proof can also treat $\alpha$'s for which $\beta(\alpha)$ is finite, with the required bound on it depending on $\lambda$; see \cite[Remark~5.3]{CEY90} for the precise condition. The main advance in \cite{CEY90} concerns the gap structure of $\Sigma^{\lambda,\alpha}$ for rational $\alpha$. The authors identify all gaps of this set (relative to the value the integrated density of states takes on them) and prove a lower bound for their lengths. Together with a continuity result like Theorem~\ref{t.avmsthm2}, this proves the existence of gaps in the irrational spectrum, provided the approximation is good enough. Actually, the authors of \cite{CEY90} prove their own version of Theorem~\ref{t.avmsthm2}, and they obtain only $\frac13$-H\"older continuity. Using better continuity statements, one could improve the final conclusion somewhat, but one would always only cover a suitable class of Liouville numbers whose measure will not exceed zero.

As with the previous two AMO problems, the status of this problem around the turn of the century was such that some nice partial results were known, but to resolve the problem completely, one would have to invent a fundamentally different approach, as all the known approaches were understood to not be sufficient to cover the entire parameter space. Consequently, finding a complete solution to AMO-Problem~\ref{prob.amo3}  also became one of the problems on Barry Simon's list of Schr\"odinger operator problems for the 21st century \cite{S00}.

A major breakthrough was obtained in the paper \cite{P04} by Puig. He was able to connect spectral localization for some coupling constant $\lambda$ to the occurrence of Cantor spectrum for the dual coupling constant $\lambda^{-1}$ (which of course implies Cantor spectrum for $\lambda$ by applying Aubry duality again).

\begin{lemma}\label{l.puiglemma}
Suppose $\lambda > 1$, $\alpha \in \T$ is Diophantine, and $E$ is an eigenvalue of $H^{\lambda,\alpha}_0$ with an exponentially decaying eigenfunction. Then the dual energy $\lambda^{-1} E$ is an endpoint of a gap of the spectrum of $H^{\lambda^{-1},\alpha}_0$.
\end{lemma}

Inspecting the assumptions of Theorem~\ref{t.jitoloc}, we see that the phase zero is always non-resonant, and hence Theorem~\ref{t.jitoloc} implies that $H^{\lambda,\alpha}_0$ has pure point spectrum (with exponentially decaying eigenfunctions). Thus, the eigenvalues of $H^{\lambda,\alpha}_0$ must be dense in the spectrum of $H^{\lambda,\alpha}_0$, which is equal to the set $\Sigma^{\lambda,\alpha}$. Lemma~\ref{l.puiglemma} then shows that the endpoints of gaps of $\Sigma^{\lambda^{-1},\alpha}$ are dense in $\Sigma^{\lambda^{-1},\alpha}$! (Here we used that $\Sigma^{\lambda,\alpha} = \lambda \Sigma^{\lambda^{-1},\alpha}$; see Theorem~\ref{t.dualspectrum}.) Thus, $\Sigma^{\lambda^{-1},\alpha}$ is a Cantor set, and by $\Sigma^{\lambda,\alpha} = \lambda \Sigma^{\lambda^{-1},\alpha}$ again, $\Sigma^{\lambda,\alpha}$ is a Cantor set as well. In other words, given the results that were known at the time, Lemma~\ref{l.puiglemma} immediately implies the following result, also stated and derived by Puig in \cite{P04}, which resolves the Ten Martini Problem for almost all parameter values.

\begin{theorem}\label{t.puigtenmartini}
Suppose that $\lambda \in (0,\infty) \setminus \{ 1 \}$ and $\alpha \in \T$ is Diophantine. Then, $\Sigma^{\lambda,\alpha}$ is a Cantor set.
\end{theorem}

Since the proof of Lemma~\ref{l.puiglemma} is relatively easy and the result and its consequences are so important, let us give some details. We start with a simple Aubry duality statement. Consider the equations
\begin{equation}\label{f.amoeve}
u(n+1) + u(n-1) + 2 \lambda \cos(2 \pi n \alpha) u(n) = E u(n),
\end{equation}
\begin{equation}\label{f.damoeve}
u(n+1) + u(n-1) + 2 \lambda^{-1} \cos(2 \pi (\omega + n \alpha)) u(n) =
(\lambda^{-1} E) u(n).
\end{equation}
Then, the following pair of statements is not difficult to prove.

\begin{lemma}\label{l.puiglem2}
{\rm (a)} Suppose $u$ is an exponentially decaying solution of
\eqref{f.amoeve}. Consider its Fourier series
$$
\hat u(\omega) = \sum_{m \in \Z} u(m) e^{2 \pi i m \omega}.
$$
Then, $\hat u$ is real-analytic on $\T$, it extends analytically
to a strip, and the sequence $\tilde u(n) = \hat u(\omega + n \alpha)$ is a
solution of \eqref{f.damoeve}.

{\rm (b)} Conversely, suppose $u$ is a solution of
\eqref{f.damoeve} with $\omega = 0$ of the form $u(n) = g(n\alpha)$ for some
real-analytic function $g$ on $\T$. Consider the Fourier series
$$
g(\omega) = \sum_{n \in \Z} \hat g (n) e^{2 \pi i n \omega}.
$$
Then, the sequence $\{\hat g (n)\}$ is an exponentially decaying
solution of \eqref{f.amoeve}.
\end{lemma}

Next we use the information provided by the previous lemma to reduce the situation at hand to constant coefficients. Here is a general statement to this effect:

\begin{lemma}\label{l.puiglem3}
Let $\alpha \in \T$ be Diophantine and suppose $A : \T \to \mathrm{SL}(2,\R)$ is a real-analytic map, with analytic extension to the strip $|\Im \omega| < \delta$ for some $\delta > 0$. Assume that there is a non-vanishing real-analytic map $v : \T \to \R^2$ with analytic extension to the same strip $|\Im \omega| < \delta$ such that
$$
v(\omega + \alpha) = A(\omega) v(\omega) \quad \text{ for every } \omega \in \T.
$$
Then, there are a real number $c$ and a real-analytic map $B : \T \to  \mathrm{SL}(2,\R)$ with analytic extension to the strip $|\Im \omega| < \delta$ such that with
\begin{equation}\label{f.puigcdef}
C = \begin{pmatrix} 1 & c \\ 0 & 1 \end{pmatrix},
\end{equation}
we have
\begin{equation}\label{f.redtoconst}
B(\omega + \alpha)^{-1} A(\omega) B(\omega) = C \quad \text{ for every } \omega \in \T.
\end{equation}
\end{lemma}

Let us explain how this is shown. Since $v$ does not vanish, $d(\omega) = v_1(\omega)^2 + v_2(\omega)^2$ is strictly positive and hence we can define
$$
B_1(\omega) = \begin{pmatrix} v_1(\omega) & - \frac{v_2(\omega)}{d(\omega)} \\ v_2(\omega) & \frac{v_1(\omega)}{d(\omega)} \end{pmatrix} \in \mathrm{SL}(2,\R)
$$
for $\omega \in \T$. We have
\begin{equation}\label{f.puig1}
A(\omega) B_1 (\omega) = \begin{pmatrix} v_1(\omega + \alpha) & \ast \\ v_2(\omega + \alpha) & \ast \end{pmatrix} \in \mathrm{SL}(2,\R)
\end{equation}
and hence
$$
A(\omega) B_1(\omega) = B_1 (\omega + \alpha) \tilde C(\omega)
$$
with
$$
\tilde C(\omega) = \begin{pmatrix} 1 & \tilde c(\omega) \\ 0 & 1 \end{pmatrix},
$$
where $\tilde c : \T \to \R$ is analytic. Indeed, by \eqref{f.puig1} the first column of $\tilde C(\omega)$ is determined and then its $(2,2)$ entry must be one since $\tilde C(\omega) = B_1 (\omega + \alpha)^{-1} A(\omega) B_1(\omega) \in \mathrm{SL}(2,\R)$. Now let
$$
c = \int_\T \tilde c(\omega) \, d\omega.
$$
and define the matrix $C$ as in \eqref{f.puigcdef}.

We claim that we can find $b : \T \to \R$ analytic (with analytic extension to a strip) such that
\begin{equation}\label{f.cohomeq}
b(\omega + \alpha) - b(\omega) = \tilde c(\omega) - c \quad \text{ for every } \omega \in \T.
\end{equation}
Indeed, expand both sides of the hypothetical identity \eqref{f.cohomeq} in Fourier series:
$$
\sum_{k \in \Z} b_k e^{2 \pi i (\omega + \alpha) k} - \sum_{k \in \Z} b_k e^{2 \pi i \omega k} = \sum_{k \in \Z} \tilde c_k e^{2 \pi i \omega k} - c.
$$
Since we have $\tilde c_0 = c$, the $k=0$ terms disappear on both sides and hence all we need to do is to require
$$
b_k (e^{2 \pi i \alpha k} - 1) = \tilde c_k \quad \text{ for every } k \in \Z \setminus \{ 0 \}.
$$
In other words, if we set $b_0 = 0$ and
$$
b_k  = \frac{\tilde c_k}{e^{2 \pi i \alpha k} - 1} \quad \text{ for every } k \in \Z \setminus \{ 0 \},
$$
then
$$
b(\omega) = \sum_{k \in \Z} b_k e^{2 \pi i \omega k}
$$
satisfies \eqref{f.cohomeq}. Since $\tilde c (\cdot)$ has an analytic extension to a strip, the coefficients $\tilde c_k$ decay exponentially. On the other hand, the Diophantine condition which $\alpha$ satisfies ensures that the coefficients $b_k$ decay exponentially as well and hence $b(\cdot)$ is real-analytic with
an extension to the same open strip.

Setting
$$
B_2(\omega) = \begin{pmatrix} 1 & b(\omega) \\ 0 & 1 \end{pmatrix} \in \mathrm{SL}(2,\R),
$$
and using \eqref{f.cohomeq}, we find
$$
B_2(\omega + \alpha)^{-1} \tilde C(\omega) B_2 (\omega) = \begin{pmatrix} 1 & c \\ 0 & 1 \end{pmatrix} = C
$$
for every $\omega \in \T$. Thus, setting $B(\omega) = B_1 (\omega) B_2(\omega)$, we obtain \eqref{f.redtoconst}.

We can now prove Lemma~\ref{l.puiglemma}. Consider an eigenvalue $E$ of $H^{\lambda,\alpha}_0$ and a corresponding exponentially decaying eigenfunction. Then, Lemma~\ref{l.puiglem2} yields the
real-analytic function $\hat u$, which has an analytic extension to a strip, and a quasi-periodic solution of the dual difference equation at the dual energy. Using this as input to Lemma~\ref{l.puiglem3}, we then obtain that
$$
A(\omega) = \begin{pmatrix} \lambda^{-1} E - 2 \lambda^{-1} \cos(2 \pi \omega) & - 1 \\ 1 & 0 \end{pmatrix}
$$
may be analytically conjugated via $B(\cdot)$ to the constant
$$
C = \begin{pmatrix} 1 & c \\ 0 & 1 \end{pmatrix}.
$$

Let us show that $c \not= 0$. Assume to the contrary $c = 0$. Then, $A(\omega) = B(\omega + \alpha) B(\omega)^{-1}$ for every $\omega \in \T$ and therefore, all solutions of \eqref{f.damoeve} are analytically quasi-periodic! Indeed,
\begin{align*}
\begin{pmatrix} u(n) \\ u(n-1) \end{pmatrix} & = A(\omega + (n-1) \alpha) \begin{pmatrix} u(n-1) \\ u(n-2) \end{pmatrix} \\
& = \cdots \\
& = A(\omega + (n-1) \alpha) \times \cdots \times A(\omega) \begin{pmatrix} u(0) \\
u(-1) \end{pmatrix} \\
& = B(\omega + n \alpha) B(\omega)^{-1} \begin{pmatrix} u(0) \\ u(-1) \end{pmatrix},
\end{align*}
that is,
$$
u(n) = \left\langle \begin{pmatrix} 1 \\ 0 \end{pmatrix} , B(\omega + n \alpha) B(\omega)^{-1} \begin{pmatrix} u(0) \\ u(-1) \end{pmatrix} \right\rangle,
$$
and hence $u(n) = g(n\alpha)$ with a real-analytic function $g$ on $\T$. Now consider two linearly independent solutions of \eqref{f.damoeve} and associate with them via Lemma~\ref{l.puiglem2} the corresponding exponentially decaying solutions of the dual equation \eqref{f.amoeve}. They must be linearly independent too, which yields the desired contradiction since by constancy of the Wronskian there cannot be two linearly independent exponentially decaying solutions. This contradiction shows $c \not= 0$.

Let us now perturb the energy and consider
$$
\tilde A(\omega) = \begin{pmatrix} (\lambda^{-1} E + \lambda^{-1}
\delta) - 2 \lambda^{-1} \cos(2 \pi \omega) & - 1 \\ 1 & 0
\end{pmatrix} = A(\omega) + \begin{pmatrix} \lambda^{-1}
\delta & 0 \\ 0 & 0 \end{pmatrix}.
$$
One can show that there is $\delta_0 > 0$ such that
\begin{equation}\label{f.puigncg}
0 < | \delta | < \delta_0 \text{ and } \delta c < 0 \quad
\Rightarrow \quad \lambda^{-1} E + \lambda^{-1} \delta \not\in
\sigma(H^{\lambda^{-1},\alpha}_0),
\end{equation}
and hence $\lambda^{-1} E$ is a gap boundary, as claimed. Lemma~\ref{l.puiglemma}, and hence Theorem~\ref{t.puigtenmartini}, now follow.

\medskip

By the nature of Puig's approach leading to Theorem~\ref{t.puigtenmartini} via Lemma~\ref{l.puiglemma}, the critical coupling, $\lambda = 1$, has to be excluded. Note, however, that for critical coupling, the spectrum has zero Lebesgue measure as discussed above, and this implies that it cannot contain any intervals. Therefore, the zero-measure results of Last \cite{L94} and Avila-Krikorian \cite{AK06} described in the previous subsection imply Cantor spectrum in this case.

\begin{theorem}\label{t.lastavilakriktenmartini}
Suppose $\alpha \in \T$ is irrational. Then, $\Sigma^{1,\alpha}$ is a Cantor set.
\end{theorem}

This leaves non-critical couplings and frequencies that are neither Lioville nor Diophantine. Avila and Jitomirskaya managed in \cite{AJ09} to close this gap by working from both sides of the intermediate region and establish Theorem~\ref{t.amotenmartini} in the form stated.

\begin{theorem}\label{t.avjittenmartini}
For every $\lambda > 0$ and every irrational $\alpha \in \T$, $\Sigma^{\lambda,\alpha}$ is a Cantor set.
\end{theorem}

\section{The Fibonacci Hamiltonian}\label{s.8}

The Fibonacci Hamiltonian is the most extensively studied operator within the context of Schr\"odinger operators with subshift potentials, which will be discussed in the next section. As in the case of the almost Mathieu operator we devote a separate section to the Fibonacci Hamiltonian to acknowledge the multitude of additional results that are known for it beyond the results that hold for the general class of models.

The Fibonacci Hamiltonian was proposed in the early 1980's by Kohmoto, Kadanoff, and Tang \cite{KKT83} and Ostlund, Pandit, Rand, Schellnhuber, and Siggia \cite{OPRSS83}. After the discovery of quasicrystals by Shechtman \cite{SBGC84}, it became the central model for the study of electron transport in one-dimensional quasicrystals. Beyond its relevance to physics, this operator is also fascinating from a purely mathematical perspective.

The Fibonacci subshift can be generated in various equivalent ways. The two most popular ones are by means of a coding of some irrational rotation of the circle and by the fixed point of the Fibonacci substitution.

Throughout this section, let $\phi$ denote the golden ratio, that is,
\begin{equation}\label{e.goldenratio}
\phi = \frac{\sqrt{5}+1}{2}.
\end{equation}
The inverse of the golden ratio is then given by $1/\phi = \phi - 1 = \frac{\sqrt{5}-1}{2}$. Denote by $R_{1/\phi}$ the irrational rotation of the circle $\T = \R / \Z$ by $1/\phi$, $R_{1/\phi} (x) = x + 1/\phi \mod 1$. A coding of the rotation $R_{1/\phi}$ is obtained by sampling the iteration of this map with a finitely-valued observable. The specific choice leading to the object of interest is obtained by the partition $\T = [0, 1- 1/\phi) \sqcup [1-1/\phi, 1)$, and mapping $[0, 1- 1/\phi)$ to zero and $[1-1/\phi, 1)$ to $\lambda > 0$. That is, for an initial point $x \in \T$, consider the coding sequence
\begin{equation}\label{e.fib1def}
\lambda \chi_{[1-1/\phi, 1)}(n/\phi + x).
\end{equation}
This family (indexed by $x \in \T$) is sometimes considered to be the family of Fibonacci potentials. The parameter $\lambda$ is naturally called the \emph{coupling constant}, while $x \in \T$ is called the \emph{phase}. These potentials do have the general form studied in this paper, with the space and transformation as specified above, and Lebesgue measure on $\T$ as the unique invariant probability measure.

The description of the Fibonacci potentials via the Fibonacci substitution goes as follows. Consider the alphabet $\mathcal{A} = \{ a, b \}$, and the map $S : \mathcal{A} \to \mathcal{A}^*$ given by $S(a) = ab$, $S(b) = a$. Here, $\mathcal{A}^*$ denotes the set of finite words over $\mathcal{A}$. One can extend $S$ to $\mathcal{A}^*$ and also to $\mathcal{A}^{\Z_+}$ by concatenation. Note that the map $S : \mathcal{A}^{\Z_+} \to \mathcal{A}^{\Z_+}$ has a unique fixed point, namely $u = S(u) = abaababaabaab\ldots$. This sequence $u$ is called the \emph{Fibonacci} (\emph{substitution}) \emph{sequence}. The \emph{Fibonacci subshift} is then given by
$$
\Omega_\mathrm{Fib} = \{ \omega \in \mathcal{A}^{\Z} : \text{ every finite subword of } \omega \text{ occurs in } u \}.
$$
It is easy to see that $\Omega_\mathrm{Fib}$ is indeed a subshift, that is, it is $T$-invariant and closed. It is also not too difficult to see that $\Omega_\mathrm{Fib}$ is strictly ergodic. The sampling function that is usually considered is the locally constant function
$$
f(\omega) = \begin{cases} \lambda & \omega_0 = a, \\ 0 & \omega_0 = b . \end{cases}
$$
The resulting potentials $V_\omega(n) = f(T^n \omega)$ are almost precisely those given in \eqref{e.fib1def} (with $\lambda$ fixed and $x$ running through $\T$.) More precisely, $\{ V_\omega \}_{\omega \in \Omega_\mathrm{Fib}}$ consists of all the sequences in \eqref{e.fib1def} and an additional orbit, namely the $T$-orbit of the sequence $\lambda \chi_{(1-1/\phi, 1]}(n/\phi \mod 1)$. Moreover, Lebesgue measure on $\T$ pushes forward via $x \mapsto \lambda \chi_{[1-\alpha, 1)}(n\alpha + x \mod 1)$ to the unique $T$-invariant Borel probability measure on $\Omega_\mathrm{Fib}$, and $R_{1/\phi}$ and $T$ are (semi-)conjugate (the conjugacy holds on $\Omega_\mathrm{Fib}$ minus the exceptional orbit, but the exclusion of this orbit is responsible for the different topologies---note that $\T$ is connected, while $\Omega_\mathrm{Fib}$ is totally disconnected).

Thus the family of Fibonacci potentials can be described in two different ways, with two different sets of choices for space, transformation, sampling function, and invariant measure. Each of the two descriptions has its advantages, and one usually works with the more appropriate choice, depending on the context and the question being studied.

Since the subshift $\Omega_\mathrm{Fib}$ is minimal, Proposition~\ref{p.constantspec} is applicable and shows that there is a common spectrum, which will be denoted by $\Sigma_\lambda$. The density of states measure will be denoted by $\nu_\lambda$. By Theorem~\ref{t.dos.eq} and Theorem~\ref{t.boshernitzan}\footnote{In the special case at hand, the absence of non-uniform hyperbolicity was first shown by S\"ut\H{o} in \cite{S89}.}, it is the equilibrium measure of the set $\Sigma_\lambda$. We will also discuss the integrated density of states $N_\lambda$ and the transport exponents $\tilde \alpha_u^\pm$. For all of these quantities, the behavior in the small and large coupling regimes have been studied. In the following four subsections we describe the known results.

\subsection{Spectrum and Spectral Measures}\label{ss.fibspectrum}

S\"ut\H{o} \cite{S89} showed the following general result about the spectrum of the Fibonacci Hamiltonian:

\begin{theorem}\label{t.fibzmspec}
For every $\lambda > 0$, the spectrum $\Sigma_\lambda$ is a Cantor set of zero Lebesgue measure.
\end{theorem}

The proof of this theorem was made possible by Theorem~\ref{t.kotthmfv}. Indeed, S\"ut\H{o} proved his result shortly after \cite{K89} was released. Recall that as a consequence of Theorem~\ref{t.kotthmfv}, we have $\mathrm{Leb} \, (\mathcal{Z}) = 0$. S\"ut\H{o} then proceeded by proving that $\mathcal{Z}$ and the spectrum actually coincide for every value of the coupling constant. That is, for every energy in the spectrum, the norm of the transfer matrix grows subexponentially. (The converse holds in general, that is, subexponential growth implies that the energy in question must belong to the spectrum; see Theorem~\ref{t.spectrumandenergypartition}.) This result was later strengthened by Iochum and Testard, who showed in \cite{IT91} that for energies in the spectrum, the norm of the transfer matrix is actually polynomially bounded. This strengthening turned out to be quite important as it is essential for establishing continuity properties of spectral measures as well as estimates for the transport exponents.

Naturally, once it is known that the spectrum has zero Lebesgue measure, one is interested in its fractal dimension. The standard quantities of interest are the Hausdorff dimension and the (upper and lower) box counting dimension of the set. The proof of the zero-measure property of the spectrum given by S\"ut\H{o}, which relies on Kotani's general result \cite{K89}, does not shed any light on these dimensions and does not allow one to obtain quantitative estimates for them. Nevertheless, quantitative results have been obtained in recent years. All of these results are proved through a study of the trace map, which provides an alternative way of understanding the Cantor structure of the spectrum, which does in fact allow one to obtain quantitative statements. We will say more on this approach in Subsection~\ref{ss.fibproofs} below.

The first important result on the fractal dimension of the spectrum is the following.

\begin{theorem}\label{t.dimeqthm}
For every $\lambda > 0$, the box counting dimension of $\Sigma_\lambda$ exists and obeys $\dim_B \Sigma_\lambda = \dim_H \Sigma_\lambda$.
\end{theorem}

This is useful as it is usually easier to estimate the Hausdorff dimension from above and the box counting dimension from below. Knowing that they are equal allows one to estimate their common value, henceforth denoted by $\dim \Sigma_\lambda$, conveniently from above as well as from below.

Theorem~\ref{t.dimeqthm} was shown for $\lambda \ge 16$ by Damanik, Embree, Gorodetski, and Tcheremchantsev in \cite{DEGT08}. For $\lambda > 0$ sufficiently small, it was shown by Damanik and Gorodetski in \cite{DG09}. A proof that works for all $\lambda > 0$ was given by Damanik, Gorodetski, and Yessen in \cite{DGY14}.

For sufficiently large values of $\lambda$, upper and lower bounds for the dimension were obtained in \cite{DEGT08}. In particular, Damanik et al.\ were able to prove the following theorem identifying the large coupling asymptotics of the dimension.

\begin{theorem}\label{t.fibspecdimlargelambda}
We have
$$
\lim_{\lambda \to \infty} \dim \Sigma_\lambda \cdot \log \lambda = \log (1 + \sqrt{2})  \approx 1.83156 \log \phi.
$$
\end{theorem}

The following result from \cite{DG11} addresses the small coupling asymptotics of the dimension of the spectrum.

\begin{theorem}\label{t.fibspecdimsmalllambda}
There are constants $c_1, c_2 > 0$ such that for $\lambda > 0$ sufficiently small, we have
$$
1 - c_1 \lambda \le \dim \Sigma_\lambda \le 1 - c_2 \lambda.
$$
\end{theorem}

The spectral measures associated with the Fibonacci Hamiltonian are only partly understood. On the one hand, the qualitative behavior is completely known:

\begin{theorem}
For every $\lambda > 0$ and every $\omega \in \Omega_\mathrm{Fib}$, all spectral measures are purely singular continuous.
\end{theorem}

The absence of eigenvalues was shown in this generality by Damanik and Lenz in \cite{DL99a} with the help of the Gordon two-block criterion, Lemma~\ref{l.twoblockgordon}. The required estimates for transfer matrix traces may be established via a study of the trace map, and this was accomplished by S\"ut\H{o} in \cite{S87}. Earlier partial results regarding the absence of eigenvalues for the Fibonacci Hamiltonian can be found in \cite{DP86, HKS95, K96, S87}. In addition, by Theorem~\ref{t.fibzmspec} the spectral measures are supported by a set of zero Lebesgue measure and hence must be purely singular.

On the other hand, the quantitative behavior is not well understood. That is, we don't have a very good handle on the optimal $\alpha > 0$ such that a given spectral measure is $\alpha$-continuous, or at least has a non-trivial $\alpha$-continuous component. Recall from Subsection~\ref{ss.hdpsm} that answers to these questions are desirable. Nevertheless, we do know that such $\alpha$'s exist:

\begin{theorem}\label{t.fibalphacont}
For every $\lambda > 0$, there is $\alpha > 0$ such that for every $\omega \in \Omega_\mathrm{Fib}$, all spectral measures are $\alpha$-continuous.
\end{theorem}

This was shown by Damanik, Killip, and Lenz in \cite{DKL00}, again using Lemma~\ref{l.twoblockgordon}. This time, the lemma needs to be applied to a number of consecutive sites in order to reproduce multiples of local $\ell^2$ norms. This mass-reproduction technique was originally developed by Damanik \cite{D98} in the case of zero phase; see also \cite{JL00}. The main question concerns the optimization of the value of $\alpha$. The optimization of the $\alpha$ that results from current technology can be found in \cite{DG11}, but it is quite clear that this result does not describe the optimal value. As a consequence, the estimates we can deduce from these spectral continuity results for the transport exponents are likely far from optimal. Identifying the actual dimensionality properties of spectral measures associated with the Fibonacci Hamiltonian remains an important open problem.

\subsection{The Density of States Measure}\label{ss.fibdos}

As pointed out above, the density of states measure $\nu_\lambda$ is the equilibrium measure of the set $\Sigma_\lambda$. The following theorem summarizes some properties of this measure, including its exact-dimensionality and the large and small coupling asymptotics of its dimension.

\begin{theorem}\label{t.fibdsmeas}
For every $\lambda > 0$, there is $d_\lambda \in (0,1)$ so that the density of states measure $\nu_\lambda$ is of exact dimension $d_\lambda$, that is, for $\nu_\lambda$-almost every $E \in \R$, we have
$$
\lim_{\varepsilon \downarrow 0} \frac{\log \nu_\lambda(E - \varepsilon , E + \varepsilon)}{\log \varepsilon} = d_\lambda.
$$
Moreover, $d_\lambda$ is an analytic function of $\lambda$, and we have
$$
\lim_{\lambda \to 0} d_\lambda = 1
$$
and
$$
\lim_{\lambda \to \infty} d_\lambda \cdot \log \lambda = \frac{5 + \sqrt{5}}{4} \log \phi \approx 1.80902 \log \phi.
$$
\end{theorem}

The exact-dimensionality at small coupling and the small-coupling asymptotics were shown by Damanik and Gorodetski in \cite{DG12}. The extension of the exact-dimensionality result to all couplings and the large coupling asymptotics are contained in the paper \cite{DGY14} by Damanik, Gorodetski, and Yessen.

\subsection{The Optimal H\"older Exponent of the Integrated Density of States}\label{ss.fibidshoeldexp}

The integrated density of states of the Fibonacci Hamiltonian depends on the coupling constant $\lambda$ and will be denoted by $N_\lambda$ in this subsection. This function is always H\"older continuous:

\begin{theorem}
For every $\lambda > 0$, there are $C_\lambda < \infty$ and $\gamma_\lambda > 0$ such that
$$
| N_\lambda(E_1) - N_\lambda(E_2) | \le C_\lambda |E_1 - E_2|^{\gamma_\lambda}
$$
for every $E_1,E_2$ with $|E_1 - E_2|$ small enough.
\end{theorem}

This result was stated in \cite{DG13}. It follows quickly, however, from the uniform H\"older continuity of spectral measures, as established in \cite{DKL00}.

The supremum of all possible exponents $\gamma_\lambda$ in this statement may be called the optimal H\"older exponent of $N_\lambda$. The asymptotic behavior of the optimal H\"older exponent in the regimes of small and large coupling was studied in \cite{DG13}, where the following results were obtained.

\begin{theorem}\label{t.fibidslargelambda}
The optimal H\"older exponent of $N_\lambda$ is asymptotically $\frac{3 \log \phi}{2 \log \lambda}$ in the large coupling regime. More precisely,

{\rm (a)} Suppose $\lambda > 4$. Then, for every
$$
\gamma < \frac{3\log(\phi)}{2\log(2\lambda + 22)},
$$
there is some $\delta > 0$ such that
$$
| N_\lambda(E_1) - N_\lambda(E_2) | \le |E_1 - E_2|^{\gamma}
$$
for every $E_1,E_2$ with $|E_1 - E_2| < \delta$.

{\rm (b)} Suppose $\lambda \ge 8$. Then, for every
$$
\tilde \gamma > \frac{3\log(\phi)}{2\log \left(\frac{1}{2} \left( (\lambda - 4) + \sqrt{(\lambda - 4)^2 - 12} \right)\right)}
$$
and every $0 < \delta < 1$, there are $E_1,E_2$ with $0< |E_1 - E_2| < \delta$ such that
$$
| N_\lambda(E_1) - N_\lambda(E_2) | \ge |E_1 - E_2|^{\tilde \gamma}.
$$
\end{theorem}

\begin{theorem}
The optimal H\"older exponent of $N_\lambda$ converges to $\frac{1}{2}$ as $\lambda\to 0$, and is strictly less than $\frac{1}{2}$ for small $\lambda>0$. More precisely,

{\rm (a)} For any $\gamma \in (0, \frac{1}{2})$, there exists $\lambda_0 > 0$ such that for any $\lambda \in (0, \lambda_0)$, there exists $\delta>0$ such that
$$
|N_{\lambda}(E_1)-N_{\lambda}(E_2)|\le |E_1-E_2|^{\gamma}
$$
for every $E_1, E_2$ with $|E_1-E_2|<\delta$.

{\rm (b)} For any sufficiently small $\lambda > 0$, there exists $\tilde \gamma = \tilde \gamma (\lambda) < \frac{1}{2}$ such that for every $\delta > 0$, there are $E_1, E_2$ with $0<|E_1-E_2|<\delta$ and
$$
|N_{\lambda}(E_1)-N_{\lambda}(E_2)|\ge |E_1-E_2|^{\tilde\gamma}.
$$
\end{theorem}

\subsection{Transport Exponents}\label{ss.fibtransport}

The transport exponents associated with the Fibonacci Hamiltonian have been studied extensively; see \cite{BLS11, D98, D05, DG14, DGY14, DKL00, DST, DT03, DT05, DT07, DT08, JL00, KKL03} for a partial list of papers containing results about them. For simplicity, let us focus in our discussion here on the time-averaged upper and lower transport exponents associated with the initial state $\psi_0$ given by $\delta_0$, for which the results are easy to state. We refer the reader to the original papers for the estimates that are known for $\tilde \beta^\pm(p)$. Most of the known results concern the time-averaged exponents; see, however, \cite{DT08}.

Notice that, contrary to the quantities considered in the previous subsections, the transport exponents formally depend on $\omega \in \Omega_\mathrm{Fib}$. However, the following result was shown by Damanik, Gorodetski, and Yessen in \cite{DGY14}.

\begin{theorem}\label{t.equaltransportexponents}
For every $\lambda > 0$, $\tilde \alpha^+_u$ and $\tilde \alpha^-_u$ are equal and independent of $\omega \in \T$.
\end{theorem}

Naturally, one is interested in upper and lower bounds for the transport exponents. Given what we have already discussed, we can note that there is always some form of transport, that is, for all parameter values, all transport exponents are positive.

\begin{theorem}
For every $\lambda > 0$, we have
$$
\tilde \alpha_l^\pm > 0.
$$
\end{theorem}

Indeed, this follows from the uniform lower bounds for the upper Hausdorff dimension of spectral measures, see Theorem~\ref{t.fibalphacont}, and the Guarneri-Combes-Last estimate \eqref{e.gclbound2}.

The Fibonacci model had long been the primary candidate for a physically relevant model that displays anomalous transport. This was finally rigorously established by Damanik and Tcheremchantsev in \cite{DT07} where it was shown that all transport exponents are strictly less than one (and eventually also less than $1/2$) for sufficiently large coupling. In fact, the upper bound established in \cite{DT07} turned out to be asymptotically sharp and the method introduced there allowed the authors to prove a corresponding lower bound in a similar way in \cite{DT08}. Combining the two results, the large-coupling asymptotics of $\tilde \alpha_u^\pm$ could be identified.\footnote{To be precise, the paper \cite{DT07} studies the case of a single $\omega \in \Omega_\mathrm{Fib}$ (the so-called zero-phase potential $\lambda \chi_{[1-1/\phi, 1)}(n/\phi \mod 1)$), while the paper \cite{DT08} considers the case of general $\omega \in \Omega_\mathrm{Fib}$. With the help of \cite{D05}, the analysis of \cite{DT07} can be extended to general $\omega \in \Omega_\mathrm{Fib}$. This idea is also behind the proof of Theorem~\ref{t.equaltransportexponents}, which was obtained later, showing complete independence of $\omega$.}

\begin{theorem}\label{t.alphafiblargelambda}
We have
$$
\lim_{\lambda \to \infty} \tilde \alpha_u^\pm \cdot \log \lambda = 2 \log \phi.
$$
\end{theorem}

The asymptotics in the small coupling regime were studied in \cite{DG14}, where the following result was obtained.

\begin{theorem}
There is a constant $c > 0$ such that for $\lambda > 0$ sufficiently small, we have
$$
1 - c\lambda^2 \le \tilde \alpha_u^\pm \le 1.
$$
\end{theorem}

\subsection{Strict Inequalities}\label{ss.fibstrictinequ}

The discussion above shows that much of the recent work on the Fibonacci Hamiltonian has focused on the following four quantities: the upper transport exponents $\tilde \alpha^\pm_u(\lambda)$, the dimension of the spectrum $\dim_H \Sigma_\lambda$, the dimension of the density of states measure $\dim_H \nu_\lambda$, and the optimal H\"older exponent $\gamma_\lambda$. By general principles, we have
$$%\begin{equation}\label{e.weakinequalities}
\gamma_\lambda \le \dim_H \nu_\lambda \le \dim_H \Sigma_\lambda.
$$%\end{equation}
This is obvious since $\Sigma_\lambda$ supports the measure $\nu_\lambda$, and the almost everywhere scaling exponent of $\nu_\lambda$ is at least as big as one that works at every point. On the other hand, there is no inequality that relates $\tilde \alpha^\pm_u(\lambda)$ to one of the other three quantities, which holds for general operators.\footnote{For example, for the Fibonacci Hamiltonian, $\tilde \alpha^\pm_u(\lambda)$ is strictly larger than the other three quantities, while for random potentials, $\tilde \alpha^\pm_u(\lambda)$ is strictly smaller than each of them.}

The large or small coupling asymptotics for these four quantities, see Theorem~\ref{t.fibspecdimlargelambda} (resp., Theorem~\ref{t.fibspecdimsmalllambda}), Theorem~\ref{t.fibdsmeas}, Theorem~\ref{t.fibidslargelambda}, and Theorem~\ref{t.alphafiblargelambda}, show that they in fact obey strict inequalities for $\lambda$ sufficiently large or sufficiently close to zero (the strict inequality between $\dim_H \nu_\lambda$ and $\dim_H \Sigma_\lambda$ for $\lambda > 0$ sufficiently small was shown in \cite{DG12} as well). In fact, Damanik, Gorodetski, and Yessen showed in \cite{DGY14} that the strict inequalities hold for all $\lambda > 0$:

\begin{theorem}\label{t.strictinequalities}
For every $\lambda > 0$, we have
\begin{equation}\label{e.inequalities}
\gamma_\lambda < \dim_H \nu_\lambda < \dim_H \Sigma_\lambda < \tilde \alpha^\pm_u(\lambda).
\end{equation}
\end{theorem}

This result is a consequence of the thermodynamic formalism and formulas, established in \cite{DGY14}, relating the four quantities to suitable dynamical quantities associated with the trace map.

The particular inequality $\dim_H \nu_\lambda < \dim_H \Sigma_\lambda$ in \eqref{e.inequalities} was conjectured by Barry Simon, based on an analogy with work of Makarov and Volberg \cite{Mak, MV, V}; see \cite{DG12} for a more detailed discussion.\footnote{The conjecture does not appear anywhere in print, but it was popularized by Barry Simon in many talks given by him in the past four years.} This inequality was shown in \cite{DG12} for $\lambda > 0$ sufficiently small, and hence the conjecture had been partially established there. It was established for all values of the coupling constant in \cite{DGY14}.

The inequality
\begin{equation}\label{e.dimspectranspexp}
\dim_H \Sigma_\lambda < \tilde \alpha^\pm_u(\lambda)
\end{equation}
in \eqref{e.inequalities} is related to a question of Yoram Last. He asked in \cite{l} whether in general $\dim_B^+ \Sigma_\lambda$ bounds $\tilde \alpha^\pm_u(\lambda)$ from above and conjectured that the answer is no. The inequality \eqref{e.dimspectranspexp} confirms this (recall from Theorem~\ref{t.dimeqthm} that $\dim_B^+ \Sigma_\lambda = \dim_B^- \Sigma_\lambda = \dim_H \Sigma_\lambda$).

\subsection{Some Comments on the Proofs}\label{ss.fibproofs}

It is an important fact that the restriction of the zero-phase potential $\lambda \chi_{[1-1/\phi, 1)}(n/\phi \mod 1)$ to $\Z_+$ is precisely the (image under $f$ of the) Fibonacci sequence $u$. In particular, this restriction is $S$-invariant. This has the following consequence. Denote the Fibonacci numbers by $\{F_k\}$, that is, $F_0 = F_1 = 1$ and $F_{k+1} = F_k + F_{k-1}$ for $k \ge 1$. Then the fact that the restriction of the potential for zero phase to the right half-line is invariant under the Fibonacci substitution implies that the matrices
$$
M_{-1}(E) = \begin{pmatrix} 1 & -\lambda \\ 0 & 1 \end{pmatrix} , \quad M_0(E) = \begin{pmatrix} E & -1 \\ 1 & 0 \end{pmatrix}
$$
and
$$
M_k(E) = T_{\lambda,0}(F_k,E) \times \cdots \times T_{\lambda,0}(1,E) \quad \text{ for } k \ge 1
$$
obey the recursive relations
$$
M_{k+1}(E) = M_{k-1}(E) M_k(E)
$$
for $k \ge 0$. Passing to the variables
$$
x_k(E) = \frac12 \mathrm{Tr} M_k(E),
$$
this in turn implies
\begin{equation}\label{e.tracerec}
x_{k+1}(E) = 2 x_k(E) x_{k-1}(E) - x_{k-2}(E)
\end{equation}
for $k \ge 1$, with $x_{-1}(E)=1$, $x_0(E)=E/2$, and $x_1=(E-\lambda)/2$. The recursion relation \eqref{e.tracerec} exhibits a conserved quantity; namely, we have
\begin{equation}\label{e.traceinvariant}
x_{k+1}(E)^2+x_k(E)^2+x_{k-1}(E)^2 - 2 x_{k+1}(E) x_k(E) x_{k-1}(E) -1 = \frac{\lambda^2}{4}
\end{equation}
for every $k \ge 0$.

Given these observations, it is then convenient to introduce the \textit{trace map}
\begin{equation}\label{e.tracemap}
T : \Bbb{R}^3 \to \Bbb{R}^3, \; T(x,y,E)=(2xy-z,x,y).
\end{equation}
The following function,
$$
G(x,y,E) = x^2+y^2+z^2-2xyz-1,
$$
is invariant under the action of $T$,\footnote{The function $G(x,y,E)$ is usually called the \textit{Fricke character} or \textit{Fricke-Vogt invariant}.} and hence $T$ preserves the family of cubic surfaces\footnote{The surface $S_0$ is called the \textit{Cayley cubic}.}
\begin{equation} \label{Slam}
S_\lambda = \left\{(x,y,E)\in \Bbb{R}^3 : x^2 + y^2 + z^2 - 2xyz = 1 + \frac{\lambda^2}{4} \right\}.
\end{equation}
It is therefore natural to consider the restriction $T_{\lambda}$ of the trace map $T$ to the invariant surface $S_\lambda$. That is, $T_{\lambda}:S_\lambda \to S_\lambda$, $T_{\lambda}=T|_{S_\lambda}$. We denote by $\Lambda_{\lambda}$ the set of points in $S_\lambda$ whose full orbits under $T_{\lambda}$ are bounded (it is known that $\Lambda_\lambda$ is equal to the non-wandering set of $T_\lambda$).

Denote by $\ell_\lambda$ the line
$$
\ell_\lambda = \left\{ \left(\frac{E-\lambda}{2}, \frac{E}{2}, 1 \right) : E \in \Bbb{R} \right\}.
$$
It is easy to check that $\ell_\lambda \subset S_\lambda$.

S\"ut\H{o} proved the following central result in \cite{S87}.

\begin{theorem}\label{t.spectrum}
For every $\lambda > 0$, an energy $E \in \R$ belongs to the spectrum $\Sigma_\lambda$ if and only if the positive semiorbit of the point $\left(\frac{E-\lambda}{2}, \frac{E}{2}, 1 \right)$ under the iterates of the trace map $T$ is bounded.
\end{theorem}

This connection shows that spectral properties of the Fibonacci Hamiltonian can be studied via an analysis of the dynamics of the trace map.

Another very important ingredient is the following. For every $\lambda > 0$, $\Lambda_\lambda$ is a locally maximal compact transitive hyperbolic set of $T_{\lambda} : S_\lambda \to S_\lambda$; see \cite{C09, Cas, DG09}. This fact allows one to use powerful tools from hyperbolic dynamics in exploring the connection between the operator and the trace map. To fully exploit this, one needs that the stable manifolds of points in $\Lambda_\lambda$ intersect the line of initial conditions, $\ell_\lambda$, transversally. This crucial fact was shown for $\lambda$ sufficiently large in \cite{Cas}, for $\lambda$ sufficiently small in \cite{DG09}, and in complete generality in \cite{DGY14}:

\begin{theorem}\label{t.transversal}
For every $\lambda > 0$, $\ell_\lambda$ intersects $W^s(\Lambda_\lambda)$ transversally.
\end{theorem}

Moreover, the measure of maximal entropy $\mu_\lambda$ for $T_{\lambda} : \Lambda_\lambda \to \Lambda_\lambda$ is related to the density of states measures $\nu_\lambda$. Namely, we have \cite{DG12, DG14}:

\begin{theorem}
For every $\lambda > 0$, the following holds. Consider a normalized restriction of the measure of maximal entropy for the trace map to an element of a Markov partition. The projection of this measure to $\ell_\lambda$ along the stable manifolds of the hyperbolic set $\Lambda_\lambda$ is equal to the normalized restriction of the push-forward of the measure $\nu_\lambda$ under $E \mapsto \left(\frac{E-\lambda}{2}, \frac{E}{2}, 1 \right)$ to the image of the projection.
\end{theorem}

Equipped with these results, the spectral analysis of the Fibonacci Hamiltonian can be completely reduced to a study of the dynamics of the Fibonacci trace map.

\section{Subshifts Over Finite Alphabets}\label{s.7}

In this section we discuss Schr\"odinger operators with dynamically defined potentials taking finitely many values. The central example is the Fibonacci Hamiltonian, which was discussed in Section~\ref{s.8}. Much of the recent advances in the study of the Fibonacci Hamiltonian were made possible through the use of sophisticated tools and results from partially hyperbolic dynamics, applied to the specific case of the trace map. Here, on the other hand, we will describe those results that can be shown without resorting to the special properties of the Fibonacci Hamiltonian. Primarily, we will discuss how one establishes one of the key features of these models, zero-measure Cantor spectrum, through a proof of the absence of non-uniform hyperbolicity for the associated Schr\"odinger cocycles. Then we will also discuss some specific classes of examples, each of which contains the Fibonacci case, and each of which generalized one particular feature of the Fibonacci Hamiltonian, namely the self-similarity or the minimal combinatorial complexity. Thus, we will discuss substitution potentials and Sturmian potentials.

For a more comprehensive survey of the spectral properties of operators associated with subshifts over a finite alphabet, we refer the reader to \cite{D07a}. We also want to mention the survey \cite{DEG}, which is written from the perspective of quasicrystal models.

Let us present the models of interest in this section. Given some finite set $\mathcal{A}$, called the \emph{alphabet}, and equipped with the discrete topology, consider the product space $\mathcal{A}^\Z$, equipped with the product topology. The shift transformation $T : \mathcal{A}^\Z \to \mathcal{A}^\Z$ acts as $[T \omega]_n = \omega_{n+1}$. Any $T$-invariant closed (and hence compact) set $\Omega \subseteq \mathcal{A}^\Z$ is called a \emph{subshift}. A subshift $\Omega$ is called \emph{minimal} if the topological dynamical system $(\Omega,T)$ is minimal, that is, every orbit $O(\omega) = \{ T^n \omega : n \in \Z \}$ is dense in $\Omega$. A subshift $\Omega$ is called \emph{uniquely ergodic} if there is a unique $T$-invariant Borel probability measure on $\Omega$. A subshift that is both minimal and uniquely ergodic is called \emph{strictly ergodic}.

Given a subshift $\Omega$ and a continuous sampling function $f : \Omega \to \R$, we consider as usual the potentials $\{ V_\omega \}_{\omega \in \Omega}$ as defined in \eqref{e.ergpotential}. There are certain sampling functions of special interest. A function $f : \Omega \to \R$ is called \emph{locally constant} if it depends only on finitely many entries. More precisely, there is $N \ge 0$ and a function $g : \mathcal{A}^{2N + 1} \to \R$ such that for every $\omega \in \Omega$, we have $f(\omega) = g(\omega_{-N} \ldots \omega_0 \ldots \omega_N)$. In fact, most papers in the literature study the case of locally constant sampling functions with $N = 0$. In other words, the alphabet consists of finitely many real numbers, and the subshift elements themselves serve as the potentials of the Schr\"odinger operators in question.

Note that in the case of a locally constant sampling function, the potentials take on only finitely many values. In particular, Kotani's theorem becomes relevant, so that we know from the outset that the Lyapunov exponent is almost everywhere positive, provided the potentials are aperiodic. In particular, all spectral measures are purely singular.

If the subshift is strictly ergodic, one usually wants to go beyond this initial observation and prove at least the following two statements: The spectrum coincides with the set $\mathcal{Z}$, and hence has zero Lebesgue measure, and the operators $H_\omega$ have no eigenvalues, so that the spectral measures are in fact purely singular continuous.

Let us stress that not a single counterexample is known, that is, within the context of strictly ergodic subshifts with locally constant sampling functions, no example is known where we do not have zero-measure Cantor spectrum and purely singular continuous spectral measures. However, the paper \cite{ADZ} constructs minimal subshifts and locally constant sampling functions so that the spectrum has positive Lebesgue measure.

\subsection{Absence of Non-Uniform Hyperbolicity}\label{ss.uniformity}

The absence of non-uniform hyperbolicity holds in great generality. That is, for a large class of strictly ergodic subshifts and all locally constant sampling functions, the set $\mathcal{NUH}$ defined in \eqref{e.nuhdef} is empty. This has the immediate consequence that $\Sigma = \mathcal{Z}$, and hence in the aperiodic case, Kotani's theorem implies zero-measure Cantor spectrum. A sufficient condition for $\mathcal{NUH} = \emptyset$ to hold is the so-called Boshernitzan condition. It holds for large classes of strictly ergodic subshifts. Let us recall the definition and some of the subshifts that may be treated in this way.

Let $\Omega$ be a strictly ergodic subshift with unique $T$-invariant measure $\mu$. It satisfies the \emph{Boshernitzan condition} (B) if
\begin{equation}\label{boshcond}
\limsup_{n \to \infty} \left( \min_{w \in \mathcal{W}_\Omega(n)} n \cdot \mu \left( [w] \right) \right) > 0.
\end{equation}
Here, $\mathcal{W}_\Omega(n)$ denotes the set of words of length $n$ that occur in elements of $\Omega$, and $[w]$ denotes the cylinder set associated with a finite word $w$, that is,
$$
[w] = \{ \omega \in \Omega : \omega_1 \ldots \omega_{|w|} = w \}.
$$

The following result was shown in \cite{DL06a} (see \cite{L02} for an important precursor).

\begin{theorem}\label{t.boshernitzan}
If the subshift $\Omega$ is strictly ergodic and satisfies {\rm (B)}, and the sampling function $f : \Omega \to \R$ is locally constant, then $\mathcal{NUH} = \emptyset$.
\end{theorem}

The idea of the proof is the following. When studying convergence of $\frac1n \log\|A^n_E(\omega)\|$, unique ergodicity allows one to bound the quantities uniformly from above. That is, the $\limsup$, equal to the Lyapunov exponent, is uniform in general, and hence the only way that uniformity can fail is that the $\liminf$ is not uniform. In this situation, one has a subsequence of words of increasing length for which the associated quantity is bounded away from the uniform $\limsup$. But if (B) holds, these outliers occur very often and this in turn implies that the average (given by the Lyapunov exponent) must be strictly smaller than the uniform $\limsup$, which is a contradiction.

Combining Theorem~\ref{t.boshernitzan} with Theorems~\ref{t.spectrumandenergypartition} and \ref{t.dos.eq}, we obtain the following statement.

\begin{coro}
If the subshift $\Omega$ is strictly ergodic and satisfies {\rm (B)}, and the sampling function $f : \Omega \to \R$ is locally constant, then $\Sigma = \mathcal{Z}$ and the density of states measure is the equilibrium measure of the spectrum.
\end{coro}

Combining this in turn with Theorem~\ref{t.kotthmfv}, we obtain:

\begin{coro}\label{c.dlcorobosh}
Suppose the subshift $\Omega$ is strictly ergodic and satisfies {\rm (B)}, the sampling function $f : \Omega \to \R$ is locally constant, and the resulting potentials $V_\omega$ are aperiodic. Then, $\Sigma$ is a Cantor set of zero Lebesgue measure.
\end{coro}

This shows that the property of having zero-measure spectrum, first established in the Fibonacci case, actually holds for a large class of subshifts. Many concrete examples are described in \cite{DL06b}. Indeed, at the time \cite{DL06a, DL06b} were written, all known subshift models with zero-measure spectrum did obey condition (B), and hence Corollary~\ref{c.dlcorobosh} served as a unifying result that also provided many new examples. More recently, however, Liu and Qu have found subshift models, generated by Toeplitz sequences, with zero-measure spectrum for which (B) fails \cite{LQ11, LQ12}.

\subsection{Substitution Subshifts}\label{ss.substitutions}

Recall that the Fibonacci model may be generated by the substitution $S(a) = ab$, $S(b) = a$. Replacing this particular substitution by a more general one, and pursuing the same construction, one can generate general substitution subshifts and associated Schr\"odinger operators. That is, fix some alphabet $\mathcal{A} = \{ a, b \}$, and consider a map $S : \mathcal{A} \to \mathcal{A}^*$. One can again extend $S$ to $\mathcal{A}^*$ and also to $\mathcal{A}^{\Z_+}$ by concatenation. Assume that the map $S : \mathcal{A}^{\Z_+} \to \mathcal{A}^{\Z_+}$ has a fixed point, $u = S(u)$. This sequence $u$ is called a \emph{substitution sequence} associated with $S$. The \emph{subshift} generated by $u$ is then given by
$$
\Omega_u = \{ \omega \in \mathcal{A}^{\Z} : \text{ every finite subword of } \omega \text{ occurs in } u \}.
$$
Now one can define Schr\"odinger operators in the usual way, namely by considering the shift transformation $T : \Omega_u \to \Omega_u$ and by choosing a sampling function $f : \Omega_u \to \R$. The sampling function should be at least continuous, but is usually assumed to be locally constant. In fact, the most common choice is that where $f$ depends only on $\omega_0$, that is, one replaces symbols from $\mathcal{A}$ by real numbers. Let us call such sampling functions \emph{letter-to-letter}.

The substitution $S$ is called \emph{primitive} if there is $k \in \Z_+$ such that for every $a \in \mathcal{A}$, $S^k(a)$ contains all symbols from $\mathcal{A}$. For example, for the Fibonacci substitution, we can choose $k = 2$ since $a \mapsto ab \mapsto aba$ and $b \mapsto a \mapsto ab$. Primitivity of the substitution $S$ has many nice consequences and much of the literature on Schr\"odinger operators with substitution potentials focuses on this special case (see, however, \cite{DL02, LO03} for some exceptions to this rule). One such consequence is that for a primitive substitution $S$, the associated subshift does not depend on the choice of the fixed point, and may therefore be denoted by $\Omega_S$. For example, the \emph{Thue-Morse substitution} $S(a) = ab$, $S(b) = ba$ has two fixed points, one obtained by iterating $S$ on $a$, and the other by iterating $S$ on $b$. Another general consequence of primitivity of $S$ is that $(\Omega_S,T)$ is strictly ergodic and satisfies (B).

Recall that one very useful way of studying the Fibonacci Hamiltonian proceeds via the trace map. The existence of the trace map is a direct consequence of the self-similarity of the underlying substitution sequence with respect to the substitution rule. Thus, for any family of Schr\"odinger operators generated by a substitution as described above, we have an associated trace map. This makes it in principle possible to relate spectral properties of these operators to dynamical properties of these maps. Trace maps will always be polynomial maps of some $\R^d$, but the size of $d$ and the degree of the polynomials may be arbitrarily large. They depend on the size of the alphabet and the complexity of the words $\{ S(a) : a \in \mathcal{A} \}$. Thus it may be difficult to actually implement this strategy.

Nevertheless, one can show the following in full generality via the trace map approach.

\begin{theorem}\label{t.zmspecsubst}
Suppose that $S$ is primitive, $f$ is letter-to-letter, and the resulting potentials $V_\omega$ are aperiodic. Then, $\Sigma$ has zero Lebesgue measure.
\end{theorem}

Under additional assumptions on $S$, this result was shown by Bovier and Ghez in \cite{BG93}. Their paper was a generalization of earlier work on the Fibonacci case \cite{S89}, the invertible/Sturmian case \cite{BIST89}, the Thue-Morse case \cite{B90}, and the period doubling case \cite{BBG91}. The result as stated is due to Liu, Tan, Wen, and Wu \cite{LTWW02}.

Whenever Theorem~\ref{t.zmspecsubst} applies, it follows that all spectral measures are purely singular. It is not clear whether the absence of a pure point components follows in the same generality. There are, however, many results that exclude point spectrum for Schr\"odinger operators generated by primitive substitutions; see \cite{D98b, D98c, D00b, D01, DGR01, DKL00, DL99a, DL03, DP86, DP91, HKS95, S87} for a partial list. It is an interesting open question whether purely singular continuous spectrum holds throughout the class of primitive substitution Hamiltonians or whether there is a counterexample.

Going beyond these qualitative results in the Fibonacci case required a deeper analysis of the trace map dynamics, in particular using tools from hyperbolic dynamics. Extending this to a more general class of substitution models is in general quite challenging, but this has been successfully implemented for invertible primitive substitutions over a two-letter alphabet by Girand \cite{G14} and Mei \cite{M14}.

\subsection{Sturmian Subshifts}\label{ss.sturmian}

Recall that the other way of generating the Fibonacci model relies on the coding of a circle rotation; see \eqref{e.fib1def}. This suggests a natural generalization. Namely, replace the inverse of the golden ratio by a general irrational number $\alpha \in \T$ and consider
\begin{equation}\label{e.sturmiandef}
\lambda \chi_{[1-\alpha, 1)}(n\alpha + x)
\end{equation}
for $x \in \T$. As in the Fibonacci case one can now either consider the family of potentials indexed by $x \in \T$ or generate a subshift $\Omega_\alpha \subset \{ 0,\lambda \}^\Z$ by taking the orbit closure of the sequence \eqref{e.sturmiandef} (which is independent of $x$). In either case one obtains a family of Schr\"odinger operators that depends on the coupling constant $\lambda > 0$ and the frequency $\alpha \in (0,1) \setminus \Q$. The ($x$/$\omega$-independent) spectrum will be denoted by $\Sigma_{\lambda,\alpha}$, and the density of states measure and the integrated density of states will be denoted by $\nu_{\lambda,\alpha}$ and $N_{\lambda,\alpha}$, respectively.

Naturally, many results will depend on the continued fraction of $\alpha$,
\begin{equation}\label{e.alphacontfrac}
\alpha = \cfrac{1}{a_1 + \cfrac{1}{a_2 + \cfrac{1}{a_3 + \cdots}}},
\end{equation}
with uniquely determined $a_k \in \Z_+$. The associated rational approximants $\frac{p_k}{q_k}$ are defined by
\begin{align*}
p_0 & = 0, \quad p_1 = 1, \quad p_{k+1} = a_{k+1} p_k + p_{k-1}, \\
q_0 & = 1, \quad q_1 = a_1, \quad q_{k+1} = a_{k+1} q_k + q_{k-1}.
\end{align*}

\subsubsection{Spectrum and Spectral Measures}

First, we have the following generalization of Theorem~\ref{t.fibzmspec}, as shown by Bellissard, Iochum, Scoppola, and Testard \cite{BIST89}.

\begin{theorem}\label{t.zmspecsturmian}
For every $\lambda > 0$ and every $\alpha \in (0,1) \setminus \Q$, the spectrum $\Sigma_{\lambda,\alpha}$ is a Cantor set of zero Lebesgue measure.
\end{theorem}

As in the Fibonacci case, this result was made possible by Kotani's short paper \cite{K89}, and \cite{BIST89} was indeed completed shortly after \cite{K89} had been released. The proof of Theorem~\ref{t.zmspecsturmian} follows the same strategy as that of Theorem~\ref{t.fibzmspec}. Namely, a trace map analysis shows that the transfer matrices have subexponentially growing norms for energies in the spectrum and then the result follows from Theorem~\ref{t.kotthmfv}. More precisely, there isn't a single trace map in this case, but rather a sequence of maps, which is determined by the sequence of partial quotients $\{a_k\}_{k \in \Z_+}$ in \eqref{e.alphacontfrac}. The implementation of this strategy is more difficult in this general case. Nevertheless, the authors of \cite{BIST89} succeeded in complete generality. Moreover, the improved power-law estimate on transfer matrix norms for energies in the spectrum obtained by Iochum and Testard in the Fibonacci case \cite{IT91} was generalized to the Sturmian case by Iochum, Raymond, and Testard in \cite{IRT92} who showed that for frequencies with bounded density, that is, those for which the (upper) density
\begin{equation}\label{e.boundeddensity}
d^*(\alpha) := \limsup_{K \to \infty} \frac1K \sum_{k=1}^K a_k
\end{equation}
is finite, one again has a power-law estimate on transfer matrix norms for energies in the spectrum.

The other result that extends in full generality from the Fibonacci case to the Sturmian case is the following.

\begin{theorem}
For every $\lambda > 0$, every $\alpha \in (0,1) \setminus \Q$, and every $\omega \in \Omega_\alpha$, all spectral measures are purely singular continuous.
\end{theorem}

In this form the result was shown by Damanik, Killip, and Lenz in \cite{DKL00}; earlier partial results can be found, for example, in \cite{BIST89, DL99a, DP86, HKS95, K96}. Again, the central tool is the Gordon two-block criterion, Lemma~\ref{l.twoblockgordon}, along with the trace estimates established in \cite{BIST89}.

Given that the Lebesgue measure of the spectrum is zero, it is again natural to study the fractal dimension of the spectrum. In general, the Hausdorff dimension of $\Sigma_{\lambda,\alpha}$ will not coincide with its box counting dimension. An in-depth study of these fractal dimensions of the spectrum at large coupling has been carried out by Liu, Peyri\`ere, Qu, Wen (in various combinations of co-authorship) in \cite{FLW11, LPW07, LQW14, LW04}. This study has culminated in the following result.

\begin{theorem}\label{t.sturmianspec}
Let $\alpha \in (0,1) \setminus \Q$ be given. Then, for every $\lambda > 24$, we have
$$
\dim_H \Sigma_{\lambda,\alpha} = s_*(\lambda) \quad \text{and} \quad \dim_B^+ \Sigma_{\lambda,\alpha} = s^*(\lambda),
$$
where $s_*(\lambda)$ and $s^*(\lambda)$ are the lower and upper pre-dimensions, respectively. We have $s_*(\lambda), s^*(\lambda) > 0$, as well as
$$
s_*(\lambda) < 1 \quad \Leftrightarrow \quad K_*(\alpha) := \liminf_{k \to \infty} (a_1 \cdots a_k)^{1/k} < \infty
$$
and
$$
s^*(\lambda) < 1 \quad \Leftrightarrow \quad K^*(\alpha) := \limsup_{k \to \infty} (a_1 \cdots a_k)^{1/k} < \infty.
$$
Moreover, if $K_*(\alpha) < \infty$ {\rm (}resp., $K^*(\alpha) < \infty${\rm )}, then
$$
L_*(\alpha) := \lim_{\lambda \to \infty} \dim_H \Sigma_{\lambda,\alpha} \cdot \log \lambda \quad \left( \text{resp., } L^*(\alpha) := \lim_{\lambda \to \infty} \dim_B^+ \Sigma_{\lambda,\alpha} \cdot \log \lambda \right)
$$
exists.
\end{theorem}

The lower and upper pre-dimensions $s_*(\lambda)$ and $s^*(\lambda)$ are defined via the band structure of the canonical periodic approximants.  This approach was pioneered by Raymond in \cite{R97}, and his work is the basis for the papers \cite{FLW11, LPW07, LQW14, LW04} mentioned above.

The quantities $L_*(\alpha), L^*(\alpha)$ can be described rather explicitly, and one recovers in particular the value $\log (1 + \sqrt{2})$ in the Fibonacci case.

Theorem~\ref{t.sturmianspec} shows that $\dim_H \Sigma_{\lambda,\alpha}$ and $\dim_B^+ \Sigma_{\lambda,\alpha}$ need not coincide, simply choose an $\alpha$ with $\liminf_{k \to \infty} (a_1 \cdots a_k)^{1/k} < \infty$ and $\limsup_{k \to \infty} (a_1 \cdots a_k)^{1/k} = \infty$.

\bigskip

At small coupling, the Fibonacci result has been extended to a much smaller set of frequencies. Namely, Mei \cite{M14} showed the following:

\begin{theorem}
Suppose the continued fraction expansion of $\alpha$ is eventually periodic, that is, there are $k_0,p \in \Z_+$ such that $a_{k+p} = a_k$ for every $k \ge k_0$. Then the box counting dimension of $\Sigma_{\lambda,\alpha}$ exists and obeys $\dim_B \Sigma_{\lambda,\alpha} = \dim_H \Sigma_{\lambda,\alpha}$. Moreover, there are constants $c_1, c_2 > 0$ such that for $\lambda > 0$ sufficiently small, we have
$$
1 - c_1 \lambda \le \dim \Sigma_{\lambda,\alpha} \le 1 - c_2 \lambda.
$$
\end{theorem}

The statement about $\dim_B \Sigma_{\lambda,\alpha} = \dim_H \Sigma_{\lambda,\alpha}$ is a corollary of Cantat's work \cite{C09}, while the proof of the small coupling asymptotics extends the proof given in the Fibonacci case by Damanik and Gorodetski \cite{DG11}.

Recall that for Sturmian models, the appropriate trace map point of view is given by a sequence of maps that corresponds to the sequence of continued fraction coefficients. Thus, if the continued fraction expansion of $\alpha$ is eventually periodic, one can combine the maps over one period, and then the dynamics of this particular map governs the evolution of the traces. Suppose one can show that this map is hyperbolic (on its nonwandering set). Then the methods and results from hyperbolic dynamics apply, as they did in the Fibonacci case. This is the reason why those results whose proofs rely on the hyperbolic theory in a crucial way may be extended to frequencies $\alpha$ with eventually periodic continued fraction expansion.

\subsubsection{The Density of States Measure}

The results for the density of states measure in the Fibonacci (cf.\ Theorem~\ref{t.fibdsmeas}) case were extended by Mei \cite{M14} in the weak coupling regime:

\begin{theorem}
Suppose the continued fraction expansion of $\alpha$ is eventually periodic. Then, there exists $\lambda_0 \in (0,\infty]$ such that for every $\lambda \in (0,\lambda_0)$, there is $d_{\lambda,\alpha} \in (0,1)$ so that the density of states measure $\nu_{\lambda,\alpha}$ is of exact dimension $d_{\lambda,\alpha}$, that is, for $\nu_{\lambda,\alpha}$-almost every $E \in \R$, we have
$$
\lim_{\varepsilon \downarrow 0} \frac{\log \nu_{\lambda,\alpha}(E - \varepsilon , E + \varepsilon)}{\log \varepsilon} = d_{\lambda,\alpha}.
$$
Moreover, $d_{\lambda,\alpha}$ is an analytic function of $\lambda$, and we have $d_\lambda < \dim \Sigma_{\lambda,\alpha}$ for $\lambda > 0$ sufficiently small, and
$$
\lim_{\lambda \to 0} d_\lambda = 1.
$$
\end{theorem}

This is the extension of the main result obtained by Damanik and Gorodetski in \cite{DG12} for the Fibonacci case. The additional results from \cite{DGY14} that are stated in Theorem~\ref{t.fibdsmeas} should extend to $\alpha$'s with eventually periodic continued fraction expansion in a similar way. For a partial result in this direction, the large coupling asymptotics of the dimension of the density of states measure was identified for frequencies of constant type by Qu in \cite{Q14}.

\subsubsection{The Optimal H\"older Exponent of the Integrated Density of States}

Theorem~\ref{t.fibidslargelambda} establishes results about the optimal H\"older exponent of the integrated density of states in the Fibonacci case. This result was extended by Munger in \cite{M14b} to some Sturmian models:

\begin{theorem}\label{t.hcont-bb}
Suppose that $\lambda > 24$ and the continued fraction expansion of $\alpha$ is constant, that is, $a_k = a$ for every $k \in \Z_+$. Then for every
$$
\gamma < \begin{cases}
\frac{2 \log \alpha}{- a \log (\lambda+5) - 3 \log(a+2)} & a > 3 \\
\frac{\log \alpha}{-\log (\lambda+5) - 3 \log (a+2)} & a = 2,3 \\
\frac{3 \log \alpha}{-2 \log(27 (\lambda+5))} & a = 1
\end{cases},
$$
there is a $\delta > 0$ such that the density of states measure $N_{\lambda,\alpha}$ obeys
$$
|N_{\lambda,\alpha}(x) - N_{\lambda,\alpha}(y)| \le |x-y|^\gamma
$$
for all $x$, $y$ with $|x-y| < \delta$.
\end{theorem}

\begin{theorem}\label{t.opt-bb}
Suppose that $\lambda > 24$ and the continued fraction expansion of $\alpha$ is constant, that is, $a_k = a$ for every $k \in \Z_+$. Then for every
$$
\tilde \gamma > \begin{cases}
\frac{2 \log \alpha}{ - a \log (\lambda-8) - \log (a) + a \log 3} & a > 2 \\
\frac{\log \alpha}{-\log (\lambda-8) + \log(a) - \log 3} & a = 2 \\
\frac{3 \log \alpha}{-2 \log (\lambda-8) - 2 \log 3} & a = 1
\end{cases},
$$
and any $0 < \delta < 1$, there are $x$ and $y$ with $0 < |x-y| < \delta$ such that $|N_{\lambda,\alpha}(x) - N_{\lambda,\alpha}(y)| \ge |x-y|^{\tilde\gamma}$.
\end{theorem}

For constant continued fraction coefficients, this identifies the asymptotic behavior of optimal H\"older exponent as $\lambda \to \infty$; compare also \cite{Q14}. More generally, the qualitative behavior is determined by the upper density \eqref{e.boundeddensity} and the lower density
$$
d_*(\alpha) := \liminf_{K \to \infty} \frac1K \sum_{k=1}^K a_k < \infty
$$
of $\alpha$. Namely, still assuming $\lambda > 24$, Munger \cite{M14b} has also shown that $N_{\lambda,\alpha}$ is H\"older continuous if $d^*(\alpha)$ is finite, and it is not H\"older continuous if $d_*(\alpha)$ is infinite.

\subsubsection{Transport Exponents}

Damanik and Tcheremchantsev \cite{DT05} established the following lower bounds for the transport exponents in the Sturmian case.

\begin{theorem}
Suppose $\lambda > 0$, $\alpha \in (0,1) \setminus \Q$ with $\max_k a_k \le C$, the phase is zero, and the initial state is $\delta_0$. With
$$
\zeta = c \, d^*(\alpha) \, \log (2 + \sqrt{8 + \lambda^2})
$$
{\rm (}$c$ is some universal constant{\rm )} and
$$
\kappa = \frac{\log (\sqrt{17} / 4)}{(C+1)^5},
$$
we have
$$
\tilde \beta^-(p) \ge \begin{cases} \frac{p+2\kappa}{(p+1)(\zeta + \kappa + 1/2)} & p \le 2 \zeta + 1, \\
\frac{1}{\zeta + 1} & p>2 \zeta + 1. \end{cases}
$$
In particular,
$$
\tilde \alpha_l^- \ge \frac{2\kappa}{\zeta + \kappa + 1/2}.
$$
\end{theorem}

Generalizing the work \cite{DT07} of Damanik and Tcheremchantsev in the Fibonacci case, Marin \cite{M10} showed the following upper bound in the Sturmian case:

\begin{theorem}
Assume that $\lambda > 20$, $\alpha \in (0,1) \setminus \Q$ is such that
$$
D := \limsup_{k \to \infty} \frac{\log q_k}{k}
$$
is finite, the phase is zero, and the initial state is $\delta_0$. Then,
$$
\tilde \alpha_u^+ \le \frac{2D}{\log \frac{\lambda - 8}{3}}.
$$
Moreover, if $a_k \ge 2$ for every $k \in \Z_+$, then
$$
\tilde \alpha_u^+ \le \frac{D}{\log \frac{\lambda - 8}{3}}.
$$
\end{theorem}

The analog of the lower bounds in the Fibonacci case from \cite{DT08}, which turned out to be asymptotically optimal in the large coupling regime, has not yet been worked out.

\section{Miscellanea}

\subsection{Potentials Generated by the Skew-Shift}

Schr\"odinger operators with potentials generated by the skew-shift are perhaps the most challenging among the ones discussed in this survey. They are the least ``random'' among those models for which ``random phenomena'' are expected to occur. More precisely, the potentials in question are quite close to one-frequency quasi-periodic potentials, and this is in fact reflected on a technical level in several ways. On the other hand, they are expected to display the full range of phenomena known to occur for random potentials, that is, positive Lyapunov exponents and spectral and dynamical localization without any largeness assumptions. In particular, if one introduces a coupling constant $\lambda$, then all these statements will hold for all $\lambda > 0$ (and they are of course most surprising for small $\lambda$, where these statements do not hold for one-frequency quasi-periodic potentials).

Let us begin with some remarks. First of all, for these statements to even be potentially true, the sampling function $f : \T^2 \to \R$ will have to depend on the second coordinate. Clearly, any sampling function that only depends on the first coordinate gives rise to a one-frequency quasi-periodic potential and hence the expected statements are know to be false. For this reason, some authors isolate the case where $f$ only depends on the second component as the special case of primary interest. Concretely, this means that one studies sampling functions of the form $f(\omega_1,\omega_2) = g(\omega_2)$ with some $h : \T \to \R$ and hence potentials, which take for $n \ge 0$ the following form,
$$
V_{(\omega_1,\omega_2)}(n) = h \left( \omega_2 + n \omega_1 + \frac{n(n-1)}{2} \alpha \right).
$$
Notice that if we replace $\alpha$ by $2 \alpha$ and then consider $( \omega_1, \omega_2 ) = ( \alpha, 0 )$, this potential takes the form
\begin{equation}\label{e.nsquaredpotential}
V(n) = h ( n^2 \alpha ), \quad n \ge 0.
\end{equation}
In this form, the potential strongly resembles a one-frequency potential (the argument $n \alpha$ has been replaced by $n^2 \alpha$), and this suggests moreover that we consider the one-parameter family of potentials
\begin{equation}\label{e.ntothegammapotential}
V(n) = h ( n^\gamma \alpha ), \quad n \ge 0,
\end{equation}
in which the cases $\gamma = 1$ and $\gamma = 2$ correspond to the quasi-periodic case and the skew-shift case, respectively. Note that, since we had to impose the restriction $n \ge 0$ above, it is natural to consider half-line operators when studying potentials given by \eqref{e.ntothegammapotential}.

The second remark concerns the regularity of the sampling function. Recall from Theorem~\ref{t.abdthm} that for generic $f \in C(\T^2,\R)$, the spectrum of the skew-shift Schr\"odinger operator is a Cantor set. This is a marked non-random-type statement, as random models can never have Cantor spectrum; see Theorem~\ref{t.randomspectrum}. Moreover, the generic spectral type also differs from what is expected to sufficiently nice sampling functions. Indeed, as pointed out earlier, the method developed by Boshernitzan and Damanik in \cite{BD08} may be applied to the skew-shift. They showed the following result for the skew-shift model in \cite{BD08}.

\begin{theorem}\label{t.skewshiftgenericcontinuousspectrum}
Suppose $\alpha \in \T$ is irrational with unbounded partial quotients. Then, there is a residual set $\mathcal{F}_\mathrm{c} \subseteq C(\T^2,\R)$ such that for every $f \in \mathcal{F}_\mathrm{c}$ and Lebesgue almost every $\omega \in \T^2$, the skew-shift Schr\"odinger operator has purely continuous spectrum.
\end{theorem}

Recall once again that Lebesgue almost all $\alpha \in \T$ satisfy the assumption of this theorem. Combining Theorem~\ref{t.skewshiftgenericcontinuousspectrum} with the general Theorem~\ref{t.avdamthm1}, we see that in the conclusion we may actually claim purely singular continuous spectrum. Thus, again the generic spectral type is singular continuous.

Nevertheless, the conjectures for sufficiently regular sampling functions are intriguing. For definiteness, let us state the following expected result for a specific potential of the form \eqref{e.nsquaredpotential}, which is currently believed to be extremely difficult to prove; compare, for example, \cite[Chapter~15]{B05b}.

\begin{ssprob}
Consider the skew-shift Schr\"odinger operator with the sampling function $f(\omega_1,\omega_2) = g(\omega_2) = 2 \lambda \cos (2 \pi \omega_2)$. Show the following statements for every $\lambda > 0$:
\begin{itemize}

\item The spectrum is an interval.

\item The Lyapunov exponent is positive at all energies.

\item The operator is spectrally and dynamically localized.

\end{itemize}
\end{ssprob}

More generally, it is expected that these properties continue to hold if $2 \lambda \cos (2 \pi \omega_2)$ is replaced by a sufficiently nice (trigonometric polynomial?, real-analytic?) non-constant function. Given that the problem above is essentially wide open for $2 \lambda \cos (2 \pi \omega_2)$, it is perhaps too soon to speculate about threshold regularity questions.

\medskip

Against this backdrop, let us now turn to the known results that are relevant in this context. The first positive result, relative to the conjectures described above, was obtained by Bourgain, Goldstein, and Schlag in \cite{BGS01}. This paper was a follow-up to work by these authors in the quasi-periodic setting \cite{BG00, GS01} and extended their approach to localization to the skew-shift case.

\begin{theorem}\label{t.bgsssloc}
Fix a non-constant real-analytic function $g$ on $\T^2$ and $\varepsilon > 0$. Then there exist $A_\varepsilon \subset \T$, whose complement has measure less than $\varepsilon$,
and $\lambda_0 = \lambda_0(\varepsilon,g) > 0$ so that for every $\alpha \in A_\varepsilon$ and every $\lambda \ge \lambda_0$, the skew-shift model with sampling function $f = \lambda g$ and base rotation by $\alpha$ has the following properties:
\begin{itemize}

\item The Lyapunov exponent is positive at all energies.

\item The integrated density of states is continuous with modulus of continuity
$$
m(t) = \exp\left( -c| \log t |^{\frac{1}{24} - } \right).
$$

\item The operators $\{ H_\omega \}$ are spectrally localized.

\end{itemize}
\end{theorem}

This result proves several of the expected properties. Unfortunately, it does not establish the absence of gaps in the spectrum and it does not say anything about the weak-coupling behavior, where the skew-shift model is expected to behave differently from quasi-periodic ones and where the conjectures are especially interesting. The latter point is not surprising as the method of proving Theorem~\ref{t.bgsssloc} is an adaptation of the proof developed in the quasi-periodic setting, and hence it cannot be expected to prove more than what is known (and true) in that context.

In other words, in order to establish the two properties that are expected for the skew-shift model, but are known to fail for one-frequency quasi-periodic models, namely absence of gaps in the spectrum and weak-coupling localization, one will need new ideas. At present, establishing the conjectures in full generality appears to be completely out of reach, and the Skew-Shift-Problem stated above is likely very hard.

Nevertheless, quite significant partial progress has been made by Bourgain \cite{B02a} and Kr\"uger \cite{K09, K12, K12b, K14}. Some results establish partial statements for the skew-shift model, while others establish full statements for models that are close to (but different from) the skew-shift model.

Bourgain has shown the following result of the first kind \cite{B02a}.

\begin{theorem}\label{t.ssbourgainpointspec}
Consider the skew-shift Schr\"odinger operator with the sampling function $f(\omega_1,\omega_2) = g(\omega_2) = 2 \lambda \cos (2 \pi \omega_2)$. Then, there exists an explicit $\omega \in \T^2$ such that for every $\lambda > 0$, there is a positive measure set of $\alpha$'s for which the operator has some point spectrum of positive Lebesgue measure.
\end{theorem}

Note here that the point spectrum is the closure of the set of eigenvalues. It is claimed that this closure has positive Lebesgue measure for $\omega$ and $\alpha$ suitably chosen.

The next result, obtained by Kr\"uger in \cite{K09}, is of the second kind.

\begin{theorem}\label{t.kruegerintervalspec}
If $h \in C(\T,\R)$, $\alpha \in \T \setminus \{ 0 \}$ and $\gamma \in (0,\infty) \setminus \Z_+$, then the Schr\"odinger operator in $\ell^2(\Z_+)$ with potential \eqref{e.ntothegammapotential} has spectrum given by the interval $[\min f -2 , \max f + 2]$.
\end{theorem}

Note that in analogy to the random case, the spectrum is given by the sum of the spectra of the Laplacian and the potential; compare Theorem~\ref{t.randomspectrum}. That is, the shape of the spectrum is truly pseudorandom in this case. Unfortunately, Theorem~\ref{t.kruegerintervalspec} excludes integer $\gamma$'s and the proof in fact does not and will not cover those $\gamma$'s. In any event, Theorem~\ref{t.kruegerintervalspec} is interesting given its generality and its relatively simple proof (relying on not quite as simple work of Boshernitzan \cite{B94}). It does show that the expected result for the shape of the spectrum of the skew-shift model holds ``nearby'' after an arbitrarily small perturbation of $\gamma$.\footnote{If the reader is concerned about us passing to half-line operators here, this is not a real issue. The skew-shift model, which is a whole-line operator, may be restricted to any half-line, and the spectrum of the whole-line model will coincide with the essential spectrum of the half-line model. Since the spectrum is purely essential in the context of Theorem~\ref{t.kruegerintervalspec}, the two spectra will in fact coincide.}

For the genuine skew-shift model, Kr\"uger has shown a weaker result, which is still very interesting in that it shows behavior different from the quasi-periodic case. Namely, in \cite{K12} he established the presence of intervals in the spectrum for sufficiently large coupling. The formulation of the result becomes nicer if one moves the coupling constant over to the Laplacian. That is, one simply multiplies the operator $\Delta + \lambda V$ by $\lambda^{-1}$, and thus considers $\lambda^{-1} \Delta + V$. This merely scales the spectrum but leaves all other characteristics unchanged.

\begin{theorem}\label{t.kruegerintervalspec2}
Suppose $\alpha \in \T$ is Diophantine, $h : \T \to \R$ is real-analytic, and $\delta > 0$. Then there exists $\lambda_0 > 0$ such that for every $\lambda > \lambda_0$, we have that the spectrum of the skew-shift operator with base frequency $\alpha$ and sampling function $f(\omega_1,\omega_2) = \lambda h(\omega_2)$, times $\lambda^{-1}$, contains the set
$$
\mathcal{E}_{h,\delta} = \{ E \in \R : \exists \, \omega \in \T \text{ such that } h(\omega) = E \text{ and } |h'(\omega)| \ge \delta \}.
$$
\end{theorem}

In other words, the spectrum of $\lambda^{-1} \Delta + V$ contains large intervals for $\lambda$ sufficiently large and it approximates the range of $h$ as $\lambda$ is sent to infinity.

Finally, there is also a partial result on the positivity of the Lyapunov exponent at small coupling. Kr\"uger shows the following in \cite{K14}.

\begin{theorem}\label{t.kruegerpositivele}
Suppose $\alpha \in \T$ is Diophantine in the sense that $\kappa = \inf_{q \in \Z_+} q^2 \|q \alpha\| > 0$, and the sampling function is given by $f = \lambda g$ with $\lambda > 0$ and $g(\omega_1,\omega_2) = \cos (2\pi \omega_2) - \cos (2\pi(\omega_1 + \omega_2))$. Then there exists $\varepsilon = \varepsilon(\kappa,\lambda) > 0$ such that the Lyapunov exponent satisfies $L(E) \ge \frac14 \log ( 1 + \lambda^2 )$ for $|E| \le \varepsilon$. Moreover, $\Sigma \cap [-\varepsilon , \varepsilon] \not= \emptyset$.
\end{theorem}

The result in \cite{K14} actually holds for more general skew-shift transformation. Taken together, Theorems~\ref{t.ssbourgainpointspec} and \ref{t.kruegerintervalspec2}--\ref{t.kruegerpositivele} demonstrate that in all three respects (namely, positive Lyapunov exponents, existence of eigenvalues, and existence of intervals in the spectrum), the skew-shift model does behave differently from one-frequency quasi-periodic operators. All three properties have been partially established, but proving them in full generality will require new ideas. Let us also mention Kr\"uger's work on the eigenvalue spacings of the finite-volume restrictions of skew-shift operators \cite{K12b}, which shows that from this perspective as well, it behaves more like a random model than a quasi-periodic model.

\subsection{Avila's Disproof of the Schr\"odinger and Kotani-Last Conjectures}

Avila obtained the following result in \cite{A14f}.

\begin{theorem}
There exist a uniquely ergodic map $T : \Omega \to \Omega$, a sampling function $f : \Omega \to \R$, and a set $S \subset \R$ of positive Lebesgue measure such that for $\mu$-almost every $\omega$, $S$ is contained in the essential support of the absolutely continuous spectrum of $H_\omega$, and for every $E \in S$ and $\mu$-almost every $\omega \in \Omega$, any non-trivial eigenfunction of \eqref{e.eve} is unbounded.
\end{theorem}

This provides a counterexample to the so-called Schr\"odinger conjecture, which had asked whether for Lebesgue almost all energies in the essential support of the absolutely continuous part of a Schr\"odinger operator, all solutions of the generalized eigenvalue equation are bounded; compare the discussion following Corollary~\ref{c.boundedsolac}.

\medskip

In the same paper, \cite{A14f}, Avila also proved the following result.

\begin{theorem}\label{t.klcounterexample}
There exist a weakly mixing uniquely ergodic map $T : \Omega \to \Omega$ and a non-constant sampling function $f : \Omega \to \R$ such that $H_\omega$ has non-empty absolutely continuous spectrum for every $\omega \in \Omega$.
\end{theorem}

Recall from our discussion in Subsection~\ref{ss.kotani} that the Kotani-Last conjecture asked whether $\Sigma_\mathrm{ac} \not= \emptyset$ implies the almost periodicity of the potentials. Theorem~\ref{t.klcounterexample} provides a counterexample to this conjecture. Note that the proof does not yield purely absolutely continuous spectrum (whereas Avila does obtain this stronger statement in the continuum setting), so that at this point it is not clear whether the implication ``purely a.c.\ spectrum $\Rightarrow$ almost periodicity of the potentials'' also fails in the discrete case. However, the existence of the counterexample in the continuum case (as shown in \cite{A14f, DY14, YZ16}) sheds a lot of doubt on this conjectural implication.

\section*{Acknowledgments}

I am deeply grateful to all my collaborators for the joy of collaboration, for teaching me a wealth of mathematics, and for their patience with my limitations. I also wish to thank Manfred Einsiedler for the invitation to write this survey and Jake Fillman, Michael Goldstein, Anton Gorodetski, Svetlana Jitomirskaya, Milivoje Lukic, Christoph Marx, and Barry Simon for useful comments.

\end{document}